\documentclass[reqno,11pt]{amsart}
\usepackage{amsmath,amsthm,amsfonts,color,graphicx,bbm,mathabx}
\graphicspath{ {code/} }
\usepackage[latin1]{inputenc}
\usepackage[makeroom]{cancel}
\usepackage{pdfsync,subfigure,multirow,float}
\usepackage{filecontents,diagbox,bbold}

\definecolor{lightGray}{RGB}{235,235,235}
\definecolor{orange}{RGB}{255,128,0}
\definecolor{ucib}{RGB}{0,36,105}
\definecolor{mygreen}{RGB}{0,128,0}
\definecolor{lightBlue}{RGB}{102,153,204}

\newtheorem{thm}{Theorem}[section]

\newtheorem{cor}[thm]{Corollary}
\newtheorem{prop}[thm]{Proposition}

\newtheorem{rems}[thm]{Remarks}
\newtheorem{rem}[thm]{Remark}

\newtheorem{deff}[thm]{Definition}
\newtheorem{exa}[thm]{Example}

\DeclareMathAlphabet{\mathpzc}{OT1}{pzc}{m}{it}

\numberwithin{equation}{section}

\begin{document}
\bibliographystyle{plain}

\title{The Curve Shortening Flow for Curves of Finite Total (Absolute)
  Curvature}

\author{Patrick Guidotti}
\address{University of California, Irvine\\
Department of Mathematics\\
340 Rowland Hall\\
Irvine, CA 92697-3875\\ USA }
\email{gpatrick@math.uci.edu}

\begin{abstract}
We revisit the well-known Curve Shortening Flow for immersed curves in
the $d$-dimensional Euclidean space. We exploit a fundamental
structure of the problem to derive a new global construction of a
solution, that is, a construction that is valid for all times and is
insensitive to singularities. The construction is characterized by
discretization in time and the approximant, while still exhibiting
the possibile formation of finitely many singularities at a finite set of
singular times, exists globally and is well behaved and simpler to
analyze than a solution of the CSF. A solution of the latter is
obtained in the limit. Estimates for 
a natural (geometric) norm involving length and total absolute
curvature allow passage to the limit. Many classical qualitative
results about the flow can be recovered by exploiting the simplicity
of the approximant and new ones can be proved. The construction also
suggests a numerical procedure for the computation of the flow which
proves very effective as demonstrated by a series of numerical
experiments scattered throughout the paper.
\end{abstract}

\subjclass[1991]{}

\maketitle

\section{Introduction}
The so-called Curve Shortening Flow is a special case of a geometric
evolution equation known as the Mean Curvature Flow whereby a spatial
curve is moved in normal direction with a speed given by its
curvature. There is a vast literature concerning the MCF even in
its simplest CSF form. In this paper we mainly consider the CSF for immersed
initial curves even though some of the results are interesting even in
the embedded case. The contributions of this paper are a new
construction of the solution for the flow in $\mathbb{R}^d$ ($d\geq
2$) by a semi-discrete approximation procedure which generates
approximating flows that can
undergo singularity formation but are always globally defined in
time. This is to be compared to the common perception since Grayson's
work \cite{Gray89} that ``The real
difficulty lies in showing that the flow [CSF] is complete. One discovers
quickly that this is equivalent to showing that the curvature remains
bounded until the entire curve shrinks to a point.''  We argue that
singularity formation is a natural phenomenon for the flow (even in
its linearized form, as we shall see) and that, in a generalized sense,
curvature does in fact not blow up until the extinction time of
the flow (when the curve shrinks to a point) but merely occasionally
concentrates. The issue of curvature blow up can be circumvented by
resorting to various concepts of weak solution. This typically results
in the loss of uniqueness (as is the case for Brakke's construction
\cite{Brakke78} applied to the CSF). The approximating flows
constructed here are easier to analyze and allow us to prove
qualitative and asymptotic properties of solutions and to gain 
insight into special solutions. Furthermore uniqueness holds in as far
it holds for classical solutions of the CSF since the solutions
obtained are classical away from singularities and continuous through
the singularities. Last but not least the construction of the
solutions suggests a natural discretization which turns out to yield
an effective algorithm for the numerical computation of the flow.

Next we summarize important contributions to the understanding of the
CSF. The embedded case was thoroughly analyzed in the ground-breaking
contributions of Gage \cite{Gage83,Gage84} and, later, of Gage and Hamilton
\cite{GH86}, in the convex case, 
and of Grayson \cite{Gray87}, in the general case. In the case of
immersed curves, seminal contributions are due to Angenent
\cite{Ang90,Ang91}. In these papers the author obtains a variety of
results concerning solutions of a class of geometric (curve) evolutions on
two dimensional surfaces of which the curve shortening flow is an
example. Abresch and Langer \cite{AL86} characterize all closed curves that
evolve homotetically and show that they serve as the model for the
asymptotics in the singularity for curves which shrink to a point
in way such that the (normalized) curvature converges in
$\operatorname{L}^1$. Along these lines Huisken \cite{Hui90} and
Angenent \cite{Ang912} show that planar curves which develop type-I
singularities are asymptotic to an Abresch-Langer curve. A result of
Altschuler \cite{Alt91} shows that space curves exhibit planar
asymptotic behavior (hence of Abresch-Langer type) in the case of
type-I singularities. For Type-II singularities, a blow-up of the
solution is shown to be asymptotic to the Gream Reaper.
Halldorsson \cite{Hall12} characterizes all self-similar solutions of
the planar CSF. For completeness we also mention the work of
Altschuler and Grayson \cite{AG92} that approximates the planar CSF by
a regularized flow for curves in $\mathbb{R}^3$ that exhibits no
singularities and exist globally in time.
\section{The Equation}
In order to formulate the problem mathematically in a way
that is convenient for the purposes of this paper, we fix some
notation. The set $\mathcal{C}=\mathcal{C}(\mathbb{R}^d)$ of all
closed curves in $\mathbb{R}^d$ ($d\in \mathbb{N}$) is defined by
\begin{equation}\label{spaceOfCurves}
  \mathcal{C}(\mathbb{R}^d)=\Big\{ X:[0,1]\to \mathbb{R}^d\,\big |\,
  X\in \operatorname{W}^{1,\infty}_\pi \bigl(
  [0,1],\mathbb{R}^d\bigr), X_r\in \operatorname{BV}_\pi\bigl(
  [0,1],\mathbb{R}^d\bigr)\Big\}, 
\end{equation}
where the subscript $\pi$ indicates that the parametrizations are
periodic. By BV we denote the space functions of bounded variation. We
shall also use the notation SBV for the space of BV functions for
which the singular part of their derivative consists only of
vector-valued Dirac masses. We can assume without loss of generality
that all parametrizations have common domain of definition
$[0,1)=:I\widehat{=}\mathbb{S}^1$. Notice that  elements of
$\mathcal{C}(\mathbb{R}^d)$ are not necessarily immersed. We therefore
denoted by $\mathcal{IC}(\mathbb{R}^d)$ the subset of immersed
curves. By immersion we mean that $X$, while allowed to exhibit
singularities (discontinuities of the tangent vector), is still
an everywhere local injection. For any $X\in \mathcal{C}$, $X(I)$ is a
rectifiable curve and its arc-length is given by 
$$
s=\overline{\varphi}_X(r)=\int_0^r|X'(\rho)|\, d \rho ,\: r\in I.
$$
Its length $\overline{\varphi}_X(1)$ is denoted by $L(X)$. Any curve in
$\mathcal{C}(\mathbb{R}^d)$ can be parametrized {\em normally},
i.e. it possesses a parametrization $Y\in \mathcal{C}(\mathbb{R}^d)$
for which $|Y'|\equiv L(X)$. Given a parametrization $X$, we denote its
normal reparametrization $X\circ\varphi_X^{-1}$ by
$R(X)$. Notice that
$$
R: \mathcal{C}(\mathbb{R}^d)\to \mathcal{C}(\mathbb{R}^d),\: X\to
X\circ\varphi_X^{-1},
$$
is a nonlinear operator and that we define
$$
\varphi_X=\frac{1}{L(X)}\overline{\varphi}_X:I\to I.
$$
We identify parametrizations that only differ by a continuous,
piecewise smooth change of variables $\varphi:[0,1]\to[0,1]$ and
allow for orientation reversing change of 
variables as they do not alter the evolution of the curve under the
Curve Shortening Flow. While this removes some ambiguity, there are
still distinct immersed, non-embedded, curves that share the same
trace set (see e.g. Figure \ref{fig:infinityLoopsCurves}). In this
sense, an immersed curve is a an equivalence class of parametrizations
and, given $X$, we sometimes denote the corresponding curve by $[X]$
but, more often, abuse notation and simply speak of the curve $X$
conflating the parametrization at hand with the equivalence class it
represents. We observe, as pointed out e.g. in \cite{Ang90}, that a curve's normal 
parametrization is unique up to a rigid transformation of
$\mathbb{S}^1$, so that, given a normal parmetrization $X$, the only
other equivalent normal parametrizations $Y$ are given by
$$
Y(s)=X(s+\theta) \text{ or by }Y(s)=X(-s+\theta),\: s\in I,
$$
for some $\theta\in I$ and where addition is interpreted modulo 1. The
tangent line (not vector) to a curve $[X]$ at one of its
points is either well-defined or experiences a discontinuity. The
tangent vector depends on the direction of parametrization. If $X'\in
\operatorname{SBV}_\pi(I)$, it is well-defined or undergoes a jump for
any chosen representative $X$ of a curve $[X]$. In particular
$$
\lim_{r\to r_0\pm}X'(r)
$$
exist for every $r_0\in I$. We say that a curve is {\em essentially}
parametrized by $X$ if there is no interval $J\subset I$ on which $X$
is constant. We denote the space of (equivalence classes of)
[immersed] curves described above by $\mathcal{C}_e(\mathbb{R}^d)$
[$\mathcal{IC}_e(\mathbb{R}^d)$] and notice that each equivalence
class has a normally and essentially parametrized representative. On
$\mathcal{C}_e(\mathbb{R}^d)$, we can use the norm given by
 \begin{equation}\label{norm}
   \| X \| _{\mathcal{C}}=\| X\|_\infty+\| X_r\|_{\operatorname{L}^1}+
  \| X_{rr}\|_{\mathcal{M}},\: X\in \mathcal{C}_e(\mathbb{R}^d),
\end{equation}
where the last term denotes the $|X_{rr}|(I)$ for the total variation
measure $|X_{rr}|$ of the vector valued measure $X_{rr}$. One has of course
to use the corresponding quotient norm for the equivalence classes of
parametrizations but we shall see later that a more geometric choice
of norm can be made that is independent of parametrization and thus
behaves like a quotient norm.
\begin{exa}
A simple example of a non-immersed curve of the kind considered that
we will revisit later is a doubly covered segment. Given two points $X^i\in
\mathbb{R}^d$ ($i=1,2$), a normal and essential parametrization of the
doubly covered segment connecting them is given by
\begin{equation}\label{dcs}
  X(s)=\begin{cases}
    (1-2s)P_0+2sP_1,&s\in[0,\frac{1}{2}),\\
    (2-2s)P_1+(2s-1)P_0,& s\in[\frac{1}{2},1).
    \end{cases}
\end{equation}
Notice that $|X_s(s)|=2|P_1-P_0|=L(X)$ and that
$$
X_s=2(P_1-P_0)\chi_{[0,\frac{1}{2})}+2(P_0-P_1)\chi_{(\frac{1}{2},1)},\:
X_{ss}=4(P_1-P_0)\bigl[ \delta_0-\delta_\frac{1}{2}\bigr].
$$
\end{exa}
For a parametrized curve $X$, its curvature (vector) is
defined by
$$
k=\frac{1}{L^2(X)}\partial_{ss} R(X).
$$
This curvature coincides with the usual curvature vector of a curve in
the smooth case and $R$ becomes the identity for normally parametrized
curves.  We have that
$$
k=k^r+k^s,
$$
where $k^r\in \operatorname{L}^1_\pi(I)$ is the regular part and
$k^s\in \mathcal{M}_\pi(I)$ is the
singular part consisting of vector-valued Dirac measures, if $X_s\in
\operatorname{SBV}_\pi$. The case
of countably many singularities is not excluded but, in this paper,
only the finite case will play a role. Then
$R(X)_{s}=\frac{L(X)}{|X_r|}X_r$ only has a finite number of jump
discontinuities.
\begin{exa}\label{dscc}
For the doubly covered segment we have that
$$
k(X)=\frac{P_1-P_0}{|P_1-P_0|^2}\bigl[ \delta_0-\delta_\frac{1}{2}\bigr].
$$
\end{exa}
The smooth {\em mean curvature flow} for a time dependent family of
curves $X=X(t,\cdot)$, $t\geq 0$, is given by
$$
V(X)=\kappa(X),
$$
where $V(X)$ is the so-called normal velocity of $X$ and
$\kappa(X)=\nu(X)\cdot k(X)$ (see \cite{Ang90}). In this formulation
the sign of $\kappa$ depends on the 
choice of the normal (in relation to the orientation of the curve) and
on the orientation chosen for the ambient space. Working with the
curvature vector, as we will, avoids this issue.
A solution of the CSF is obtained if solutions of $X_t=k(X)$ can be
produced. Notice that
$$
R:\mathcal{C} _e(\mathbb{R}^d)\to \mathcal{C}_e(\mathbb{R}^d),\:
X\to X\circ\varphi_X^{-1}, 
$$
is well-defined since $\varphi_X$ is strictly increasing for $X\in
\mathcal{C} _e(\mathbb{R}^d)$ even if $|X_r|$ can have some
zeros. From now on, we shall refer to the PDE system
\begin{equation}\label{csf}
\begin{cases} \partial _t X=\frac{1}{L^2(X)}\partial _{ss} R(X),&
  t>0,\\
 X(0)=X^0\in \mathcal{C}_e,&\end{cases}
\end{equation}
as the {\em curve shortening flow (CSF)} for the initial curve $X^0$. It is
apparent that the nonlinear nature of \eqref{csf} only stems from the
normal reparametrization operator $R$ and the length operator $L$. It
can be verified that 
$$
R(\partial_tX)\neq \partial_tR(X),
$$
in general and, therefore, no simple equation can be derived if $X$ is
replaced by $R(X)$. It is, however, an important observation that
$$
 R(X_t)\cdot \nu(X)=R(X)_t\cdot \nu(X).
$$
Here and above $X_t$ is a tangent vector(field) at $X$ to the manifold of
immersed curves and $R(X_t)=X_t \circ\varphi _X^{-1}$ is the natural
extension of $R$ to the tangent space at $X$. It follows that
$[R(X)]=[X]$ can be ``updated'' using the alternative equation 
$$
R(X)_t=\frac{1}{L^2\bigl(R(X)\bigr)}\partial _{ss} R(X),
$$
while still yielding a solution of the mean curvature flow, which only
involves the normal velocity $V(X)$ as observed above. This
fact is the essential motivation for the construction of the solution
and for the discretization presented later as well as the reason for the
efficacy and simplicity of the numerical scheme. We
reiterate, however, that it is not possible to simply replace the CSF
\eqref{csf} by the equation
$$
 Y_t=\frac{1}{L^2(Y)}\partial_{ss}Y
$$
for $Y=R(X)$, because this equation does not preserve the normalized
length of a parametrization. 

Notice that $L$ is a nonlinear function of
$X$ but $L(X)=L(R(X))=|\partial_sR(X)|$ only depends on
time.
\begin{rem}\label{square}
We point out that Definition \eqref{csf} of the CSF is more general
than the usual differential geometric definition. It is solely based
on the existence of a generalized curvature vector. The latter is
sometimes defined even when the actual curvature is not. If a curve is
smooth, then we have the relation $\kappa(X)=k(X)\cdot \nu(X)$. In the
non-smooth case, take for instance a square of unit length and
normal-length parametrized by
$$
X(s)=\begin{cases}
  (s,0),&s\in[0,\frac{1}{4}),\\
  (\frac{1}{4},s-\frac{1}{4}),& s\in[\frac{1}{4},\frac{1}{2}),\\
  (\frac{3}{4}-s,\frac{1}{4}) ,& s\in[\frac{1}{2},\frac{3}{4}),\\
  (0,1-s),& s\in[\frac{3}{4},1).
\end{cases}
$$
Then
$X_{ss}=(1,1)\delta_0+(-1,1)\delta_{\frac{1}{4}}+(-1,-1)\delta_{\frac{1}{2}}+(1,-1)\delta_{\frac{3}{4}}$
is a measure and the unit normal is merely $L^\infty$ and these can therefore not
be multiplied together. Likewise the Euclidean norm of $X_{ss}$ cannot be
computed.
\end{rem}
Motivated by the above example, we investigate curvature and total
curvature in the context of curves with singularities, starting in the
plane. For smooth curves, it is known that the total curvature $\kappa_{tot}$ of a curve
parametrized by $X$ can be computed by
$$
\kappa_{tot}(X)=\int_0^{L(X)}|\partial_{\bar s\bar s}X|\,
d\bar\sigma=\frac{1}{L(X)}\int_0^1|\partial_{ss}X|\, d \sigma,
$$
where $\bar s$ is the arc-length and $s$ the normalized length,
i.e. $\bar s= L(X)s$. In the smooth case this is equivalent to
$$
\kappa_{tot}(X)=\int _0^L| \partial_{\bar s}\theta|\, d\bar \sigma =\int_0^1|\partial
_s \theta |\, d \sigma,
$$
where $\theta(\bar s)=\arccos(\partial_{\bar s}X\cdot
e_1)=\arccos(\frac{1}{L}\partial _sX\cdot e_1)$ measures the angle 
the tangent vector makes with the vector $e_1$ along the curve. This
is an immediate consequence of writing
$$
 \partial_{\bar s}X=\bigl( \cos(\theta(\bar s)),\sin(\theta(\bar
 s))\bigr), \: \partial_{s}X=L\bigl( \cos(\theta( s)),\sin(\theta(s))\bigr),
$$
and computing $|\partial_{\bar s\bar s}X|=|\partial_{\bar
  s}\theta(\bar s)|$.
In particular it always holds that
\begin{equation}\label{ktotid}
\int_0^L|\partial_{\bar s\bar s}X|\, d\bar \sigma =\int _0^L |\partial_{\bar
  s}\theta |\, d\bar \sigma.
\end{equation}
In the singular case, it would seem natural to replace the length of
the curvature vector $|\partial_{\bar s\bar s}X|$, now a Radon
measure, with its variation (measure) $|\partial_{\bar s\bar s}X|$ and
total curvature by
$$
|\partial_{\bar s\bar s}X|\bigl([0,L)\bigr)=\frac{1}{L}|\partial_{ss}X|(I).
$$
It, however, turns out that identity \eqref{ktotid} is no longer valid for
singular curves. Take for instance the square discussed in Remark
\ref{square} to see that 
$$
\partial_{\bar s}X=e_1\chi_{[0,\frac{1}{4})}+e_2\chi_{[\frac{1}{4},\frac{1}{2})}-
e_1\chi_{[\frac{1}{2},\frac{3}{4})}-e_2\chi_{[\frac{3}{4},1)}.
$$
Then
$$
\partial_{\bar s\bar s}X=(e_1+e_2)\delta_0+(e_2-e_1)
\delta_{\frac{1}{4}}+(-e_1-e_2)\delta_{\frac{1}{2}}+(e_1-e_2)\delta_{\frac{3}{4}},
$$
which has a compelling geometric interpretation, and
$$
|\partial_{\bar s\bar s}X|=\sqrt{2}\sum_{i=0}^3
\delta_{\frac{i}{4}},\: |\partial_{\bar s\bar s}X|(I)=4\sqrt{2}.
$$
On the other hand,
\begin{equation*}
\theta=0\chi_{[0,\frac{1}{4})}+\frac{\pi}{2}\chi_{[\frac{1}{4},\frac{1}{2})}+
        \pi\chi_{[\frac{1}{2},\frac{3}{4})}+\frac{3\pi}{2}\chi_{[\frac{3}{4},1)},\:
\partial_{\bar s}\theta =\frac{\pi}{2}\sum_{i=0}^3 \delta_{\frac{i}{4}},
\end{equation*}
so that
$$
|\partial_{\bar s}\theta |(I)=2\pi,
$$
which is indeed the total curvature of the square. We conclude that
$\partial_{\bar s\bar s}X$ and $\partial_{\bar s}\theta$, while
essentially containing the same information and being useful
quantities, are not related in the
obvious way they are for smooth curves. This is due to the fact that
nonlinear operations are involved in the transition between these
quantities that are not well-defined for measures. For the singular
curves considered in this paper, we shall from now on define the
curvature vector by $\frac{1}{L^2}\partial_{ss}X$ or $\partial_{\bar
  s\bar s}X$ and total curvature by
\begin{equation}\label{ktot}
\kappa_{tot}(X)=\frac{1}{L(X)}|\partial_{ss}X|(I).
\end{equation}
At least in $d=2$, one could use the alternative definition
$$
\kappa_{tot}(X)=\int _I \big | (d\theta)^r\big | +\big | (d\theta)^s\big | (I) 
$$
where $(d \theta)^r$ and $(d \theta)^s$ are the regular and singular
(jump) part of $d \theta$, respecitvely. 
The difference between $|\partial_{ss}X|(I)$ and
$|\partial_s\theta|(I)$ for singular curves in $\mathbb{R}^d$ has been
discussed before by \cite{Sull08}. The
unit tangent to a curve traces a curve in 
$\mathbb{S}^{d-1}$, the so-called tantrix. Its length as a
(discontinuous) curve in $\mathbb{S}^{d-1}$ corresponds to
$|\partial_s\theta|(I)$, whereas its length as a curve in
$\mathbb{R}^d$ yields $|\partial_{ss}X|(I)$. The former has the
advantage of coinciding with the limit of the total curvature of
approximating polygons (defined as their total turning angle). Notice
that, for smooth curves, both quantities coincide since,
infinitesimally, the ambient space distance and the unit sphere
distance coincide and the two diverge only when singularities are
present. Clearly these total curvatures are either both finite or both
infinite and there is a simple relation between them as explained in
\cite[Section 3]{Sull08}. In the same paper, it is 
also shown, following \cite{Mil50}, that the total
curvature of a curve $\Gamma$ in $\mathbb{R}^d$ can be defined by
$$
\kappa_{tot}(\Gamma)=\sup_{\mathbb{P}<\Gamma}\kappa_{tot}(\mathbb{P}),
$$
where the supremum is taken over all possible polygons generated by
selecting an (ordered) sequence $\mathbb{P}=(X_0,\dots,X_{m-1})$ of
points from the curve $\Gamma$, and
$$
\kappa_{tot}(\mathbb{P})=\sum_{i=0}^{m-1}\theta_i,
$$
where $\theta_i$ is the turning angle at $X_i$, i.e. the angle between
the unit vector in the direction determined by $X_{i-1}$ and $X_i$, and
that in the direction determined by $X_i$ and $X_{(i+1)\mod m}$. We
shall make use of this fact to compute a numerical approximation of
this quantity later in the paper.

With total absolute curvature in hand, we can introduce the following more
geometric norm on the space $\mathcal{C}(\mathbb{R}^d)$ defined by
 \begin{equation}\label{gnorm}
   \| X \| _{\mathcal{G}}=\| X\|_\infty+L(X)+\kappa_{tot}(X)
 \end{equation}
for the curve with parametrization
$X\in\mathcal{C}_e(\mathbb{R}^d)$. Notice that its value is 
independent of the chosen parametrization within the same equivalence 
class. In particular $\| X\|_{\mathcal{G}}=\| R(X)\|_{\mathcal{G}}=\|
R(X)\|_{\mathcal{C}}$ and so this quantity bounds the quotient norm of $X$.
\section{Steady States}
Among the shapes that evolve only by rescaling there may be some that
remain normally parametrized if they are initially normally
parametrized, i.e. curves for which it holds that $R\bigl((X(t)\bigr)=X(t)$
whenever $R(X^0)=X^0$. Such curves satisfy the quasi-linear equation
$$
\partial_tX=\frac{1}{L^2(X)}\partial_{ss}X.
$$
In this case, we have that $X(t)=r(t)X^0$, where $r(0)=1$, so that the
equation becomes
$$
\dot r(t)X^0=\frac{1}{L^2_0r(t)}\partial_{ss}X^0,
$$
where $L_0=L(X^0)$. By separation of variables, it must hold that
\begin{equation}\label{noReparam}
  \begin{cases}
  \partial_{ss}X^0=\mu X^0,&\\
  L_0^2\dot r(t)r(t)=\mu,& r(0)=1,
\end{cases}
\end{equation}
for some $\mu\in \mathbb{R}$. Then $\mu=-4\pi^2n^2$ for $n\in
\mathbb{N}\cup\{0\}$ and
$$
X^0=a\cos(2\pi n s)+b\sin(2\pi n s),\: r(t)=\sqrt{1+2\frac{\mu}{L_0^2}t},
$$
for some $a,b\in \mathbb{R}^d$. If $n=0$, then the only (connected)
solution would correspond to a constant point (where $r(0)=0$) but this solution is not considered as it
will always be assumed that $L_0>0$, when it needs to hold that
$r(0)=1$. For $n\in \mathbb{N}$, it must hold that
$$
\big | -2\pi na\sin(2\pi n s)+2\pi nb\cos(2\pi n
s)\big|^2=L_0^2,\: s\in I,
$$
which implies that
$$
4\pi^2n^2\bigl[ |a|^2\sin^2(2\pi n s)+|b|^2\cos^2(2\pi n
s)\bigr]-8\pi^2 n^2 a\cdot b \sin(2\pi n s)\cos(2\pi n s)\equiv L_0^2.
$$
This entails that
$$
 |a|=|b|=\frac{L_0}{2\pi n}\text{ and that }a\cdot b=0.
$$
The curve is therefore a circle traversed $n$-times that lies in the
plane spanned by $a$ and $b$. If that circle is normalized to have
length 1, then the curve has length $2\pi n$ and $a,b$ are
orthonormal. Moreover, this $n$-fold circle shrinks to a point when
$$
\sqrt{1-2\frac{4\pi^2n^2}{4\pi^2n^2}t}=\sqrt{1-2t}=0,
$$
i.e. when $t=\frac{1}{2}$. It is known that simple closed curves evolve
reducing the enclosed area at the constant rate $2\pi$. Thus the
multiply traversed circle decreases its length at a speed that is
higher so as to make the extinction time the same for all $n\in
\mathbb{N}$, in this simple case of a non-simple curve.
\begin{rem}
There are self-similar solutions, even when considering closed curves
only, which evolve by rescaling but for which $R \bigl( X(t)\bigr)
\neq X(t)$ even if $R(X^0)=X^0$. All planar curves evolving in this
way have been identified in \cite{Hall12} by the use an ODE
system. They include the curves already known from \cite{AL86} which
are closed curves of finite length that shrink homotetically. These
are characterized by the fact that the parametrization does not remain
normal along \eqref{csf} and hence
$$
X(t,s)=r(t)X^0 \bigl( \varphi(t,s)\bigr),\: s\in I,\: t\geq 0, 
$$
where $\varphi(0,s)=s$ for $s\in I$. The CSF in this cases reduces to
$$
\dot r
X^0\circ\varphi\cdot\nu^0\circ\varphi=\frac{1}{r}k^0\circ\varphi\cdot
\nu^0\circ\varphi,
$$
or, equivalently, to
$$
\dot r rX^0\cdot \nu^0=k^0\cdot \nu^0.
$$
By rescaling $X^0$ if necessary, we can assume that $\dot r r=-1$ and
thus homotetically evolving curves must satisfy
$$
X^0\cdot \nu^0+k^0\cdot\nu^0=0.
$$
From here an ODE can be derived for the curvature $\kappa^0$ which
allows for the classification of all possible solutions. A thorough
discussion is found e.g. in \cite[Appendix E]{Man11}.
\end{rem}
\begin{rem}
Assuming smoothness, the CSF flow is equivalent to the system
$$\begin{cases}
  X_t=\Sigma_*k(X),& k(X)=\frac{1}{L^2}R(X)_{ss},\\
  \sigma _t= \sigma \bigl[\int_0^1 k(X)\cdot k(X)-\Sigma_*k(X)\cdot
  \Sigma_*k(X)\bigr],&s=\int_0^r \sigma(t,\rho)\, 
  d\rho =:\Sigma(r),\\
  L_t=-L\int _0^1 k(X)\cdot k(X),&
\end{cases}
$$
where $X(0)=X^0$, $L(0)=L(X^0)$, and $\sigma(0,r)=|X^0_r(r)|$, $r\in
I$. This system is obtained using the well-known equations that can be
derived for the CSF as they can be found in \cite[Section 3]{GH86},
for instance. It shows that a solution is continuously undergoing
reparametrization unless its curvature is constant.
\end{rem}
\section{Construction of a Solution}
In order to obtain a solution to the CSF for any given initial datum
in $\mathcal{C}_e(\mathbb{R}^d)$ and in order to derive an effective
numerical approach to its solution (see later) we will heavily rely on
formulation \eqref{csf}.

The idea is to linearize \eqref{csf} by intertwining
reparametrizations and small time interval evolutions by the linear
heat equation while keeping track of singularity formation. We define
recursively
\begin{equation}\label{solApprox}
  \begin{cases}
    X^n(t)=e^{\frac{t-kh}{L^2(X^n(kh))}\partial_{ss}}R(X^n(kh)),&
    t\in\bigl[kh,(k+1)h\bigr),\: k\geq 0,\\
    X^n(0)=X_0=R(X_0),&
  \end{cases}
\end{equation}
where $X^n(kh)$ has been determined in the previous step. We point out
a last time that this evolution is an evolution for the immersed curve
$[X^0]$ and that the use of the operator $R$ preserves the equivalence
class, i.e. the curve. Here we
switched the notation for the initial parametrized 
curve to $X_0$,  we set $h=\frac{1}{n}$ for $n\in \mathbb{N}$ and
$$
L(X^n(kh))=\int _0^1 |\partial_rX^n(kh)|\,
d\rho=\big|\partial_sR\bigl(X^n(kh)\bigr)\big|, 
$$
During this evolution it is possible to encounter (singular) times
$t_{sing}$ at which $|\partial_sX^n(t_{sing})|=0$ at one or several parameter
values. At such times the evolving curve experiences one or more singularities.  
In any time interval of length $h$ (and overall) this happens at most a finite
number of times as explained below. For future use, we also define
$$
Y^n(t)=
\begin{cases}
  X^n(t),&t\not\in h \mathbb{N}\\
  e^{\frac{h}{L^2(X^n(kh))}\partial_{ss}}R(X^n(t-h)),&t\in h \mathbb{N}.
\end{cases}
$$
Letting $\Gamma^n(t)$ denote image set of
$X^n(t)$, $Y^n(t)$ is a smooth (analytic) parametrization of it for
each time $t>0$.

We interpret $X^n$ (or $\Gamma^n$) as the $n$-th interpolant
(approximant) to the CSF and collect some important properties it
enjoys in the next sections. In particular, \eqref{solApprox} defines
an orbit of curves that depends continuously on time as does its
length. The tangent vector(field) to the orbit is piecewise smooth and
bounded (in the appropriate sense). The approximant does converge (as
$n\to\infty$) to a globally defined limit $X$, which becomes constant
(the curve shrinks to a point) in finite time and satisfies the CSF
equation up to a finite number of singular times (that include the
extinction time).
\begin{rem}\label{CoC}
The function $Y^n(t)$ satisfies a linear heat equation for all $t\neq kh$,
$k\in \mathbb{N}$ and $Y^n(kh)$ is the end value of a linear heat
equation. It follows that $V+MY^n(t)$ also satisfies a linear heat
equation for any $V\in \mathbb{R}^d$ and $M\in \mathbb{R}^{d\times
  d}$. In particular any rigid change of coordinates in the ambient
space yields parametrizations that satisfy the same linear
heat equation. 
\end{rem}
\subsection{Extremity Points and Singularities}
\begin{deff}
A point $P$ on a curve smoothly parametrized by $X$ is called an
extremity point if $\partial_sX^i(r_0)=0$ for at least one $i\in\{
1,\dots, d\}$ and $X(r_0)=P$. Singular points are special extremity
points where all components of $\partial _sX$ vanish
simultaneously. The set of extremity points is denoted by $\mathcal{E}
(X)$. We say that a curve $\Gamma$ has finitely many extremity points
if there is a coordinate system (an origin and an orthonormal
basis) in which it admits a parametrization $X$ with
$\mathcal{E}(X)<\infty$. 
\end{deff}
\begin{rems}
{\bf (a)} We include inflection points of components in the definition of
extremity point for simplicity, but, minima and maxima are the
interesting critical points for the evolution. This is related to the
fact that diffusion instantenously removes inflection points.\\
{\bf (b)} A singularity point along a smoothly and essentially
parametrized curve is a special extremity point where all components
of the tangent vector vanish simultaneously.\\
{\bf (c)} If a curve is parametrized in a non-smooth way, then
extremity points also include local maxima and minima at which the
parametrization is not differentiable.\\
{\bf (d)} The only way to parametrize about a singular point in a
smooth manner is by ``stalling'' at the singularity. In this way all
components of the parametrization must vanish at the singular
points. Inflection points are generated in all components for which
the singular point is not a local extremum.
\end{rems}

\begin{prop}
Given $X_0$ with finitely many extremity points, the number of
extremity points of $X^n(t)$ is a non-increasing function of time,
regardless of $n\in \mathbb{N}$.
\end{prop}
\begin{proof}
Under scalar diffusion maxima decrease and minima increase, while
inflection points disappear instantaneously. As maxima and minima may
merge during the evolution, the number of extremity points will not
increase, as diffusion does not allow for the generation of new
extrema. Inflection points can appear but only when extrema merge,
like, e.g., when a minimum and a maximum coalesce. Non-smooth
extremity points are
instantaneously regularized, possibly turning into smooth extremal
points. Notice that the operator $R$ does not affect smooth extremity points,
which are preserved after reparametrization (their location in 
parameter space can of course change in the process). Extrema that do
occur at singular points are also preserved by $R$, while inflection
points that occur at a singular point may disappear when $R$ is
applied.
\end{proof}
\begin{rem}
If the initial curve is not smooth and possesses finitely many
extremity points, we can assume without loss of generality that it be
smooth and exhibit a finite number of extremity points by the
regularizing effect of diffusion. As diffusion does not generate
oscillations, a small time application of it removes singularities
without affecting the (macroscopic) oscillatory properties of the
initial components since we are assuming that they have at most
finitely many extremity points.
\end{rem}
An important class of curves with finitely many extremity points
(especially with respect to numerical approximation) is that of closed
polygons $\mathbb{P}=(P_0,P_1,\dots,P_{n-1})$ for $n\in \mathbb{N}$
which we identify with the closed curve
$$
\bigl( 1-(nt-i) \bigr) P_i+(nt-i)P_{(i+1)\operatorname{mod}n},\text{
  if }t\in [\frac{i}{n},\frac{i+1}{n}),\: i=1,\dots, n-1.
$$
\begin{prop}
It holds that $|\mathcal{E}(\mathbb{P})|<\infty$ for any closed polygon
$\mathbb{P}$.
\end{prop}
\begin{proof}
A polygon comprises a finite number of segments, each determining a
direction $u_i=\overrightarrow{P_iP_{(i+1)\operatorname{mod}n}}$
and an orthogonal hyperplane $\mathcal{P}_i^\perp$, $i=1,\dots,n$ . It suffices
to choose a coordinate system, the basis vectors of which do not
lie in any of these orthogonal planes.
\end{proof}
\subsection{Flow Singularities}
\begin{figure}
\begin{center}
  \includegraphics[scale=0.4]{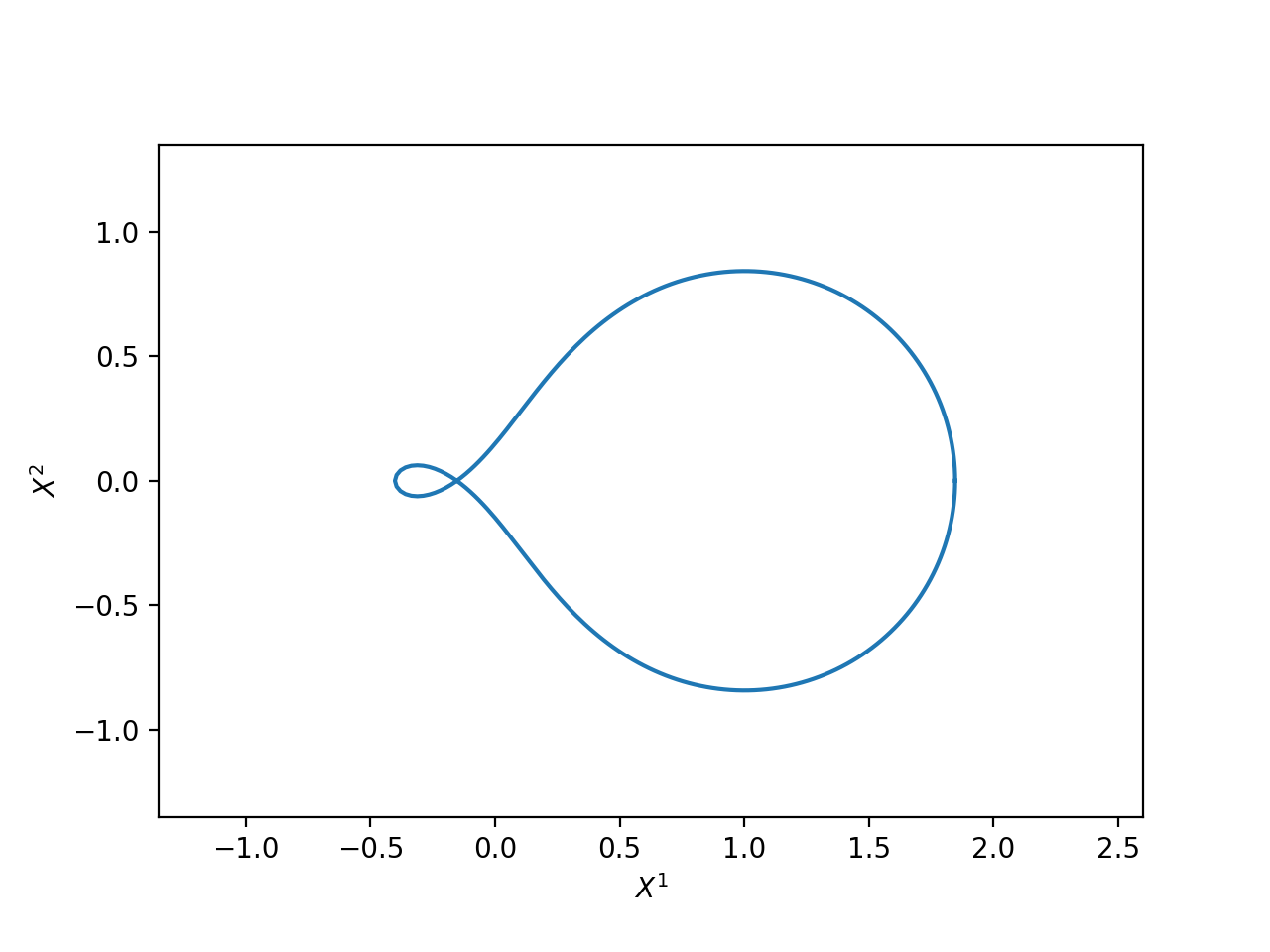}
  \includegraphics[scale=0.4]{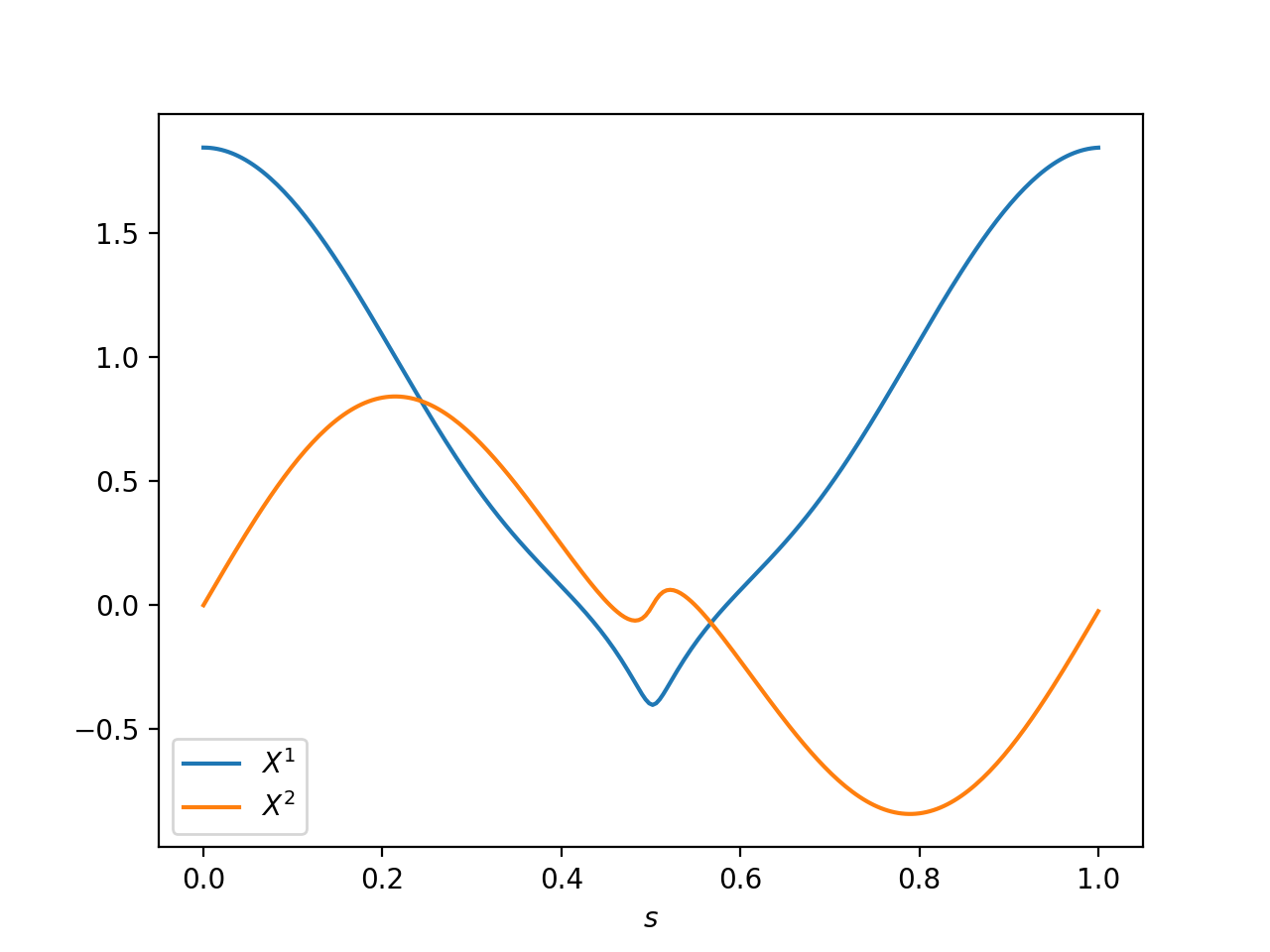}
  \includegraphics[scale=0.4]{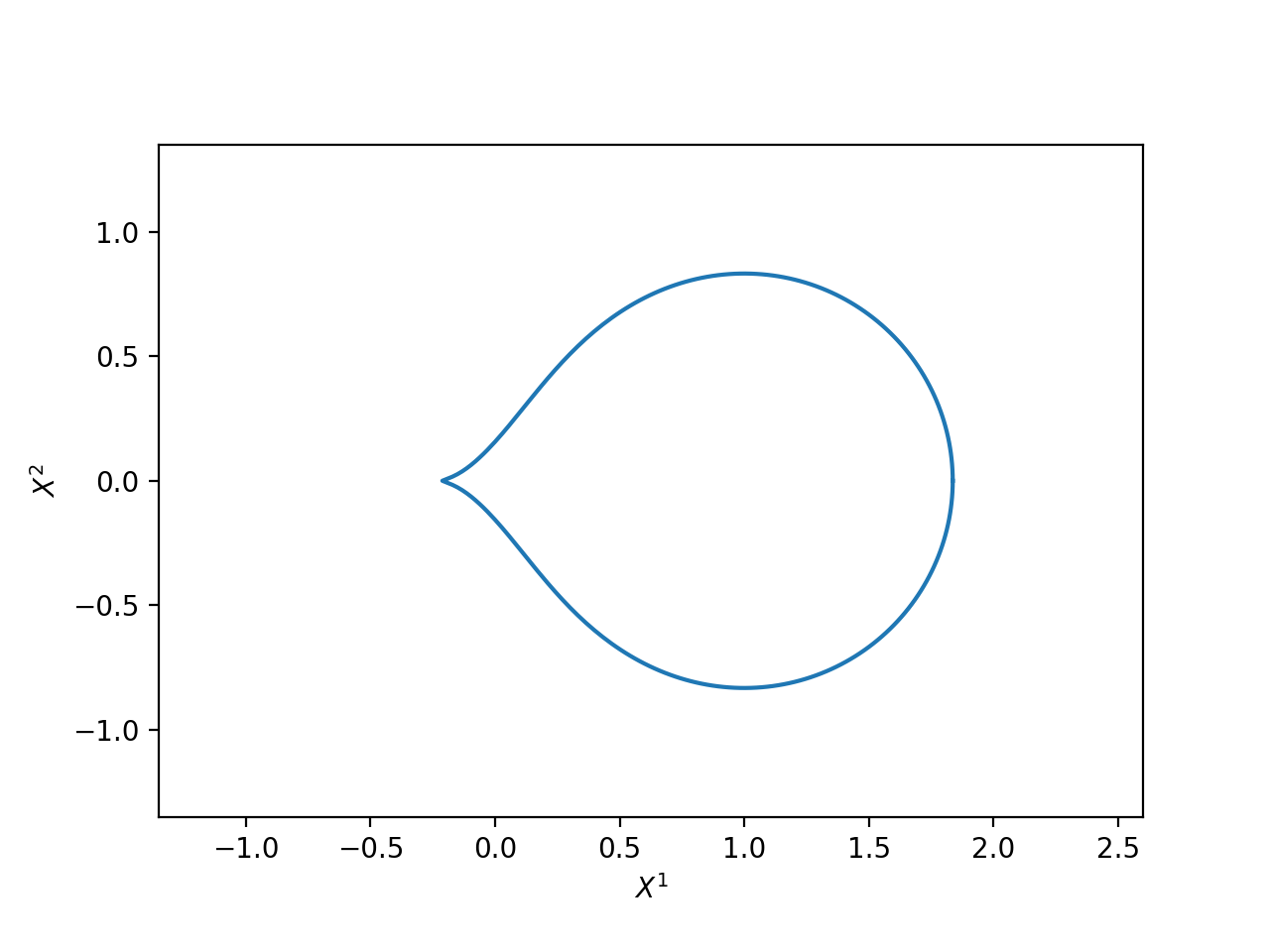}
  \includegraphics[scale=0.4]{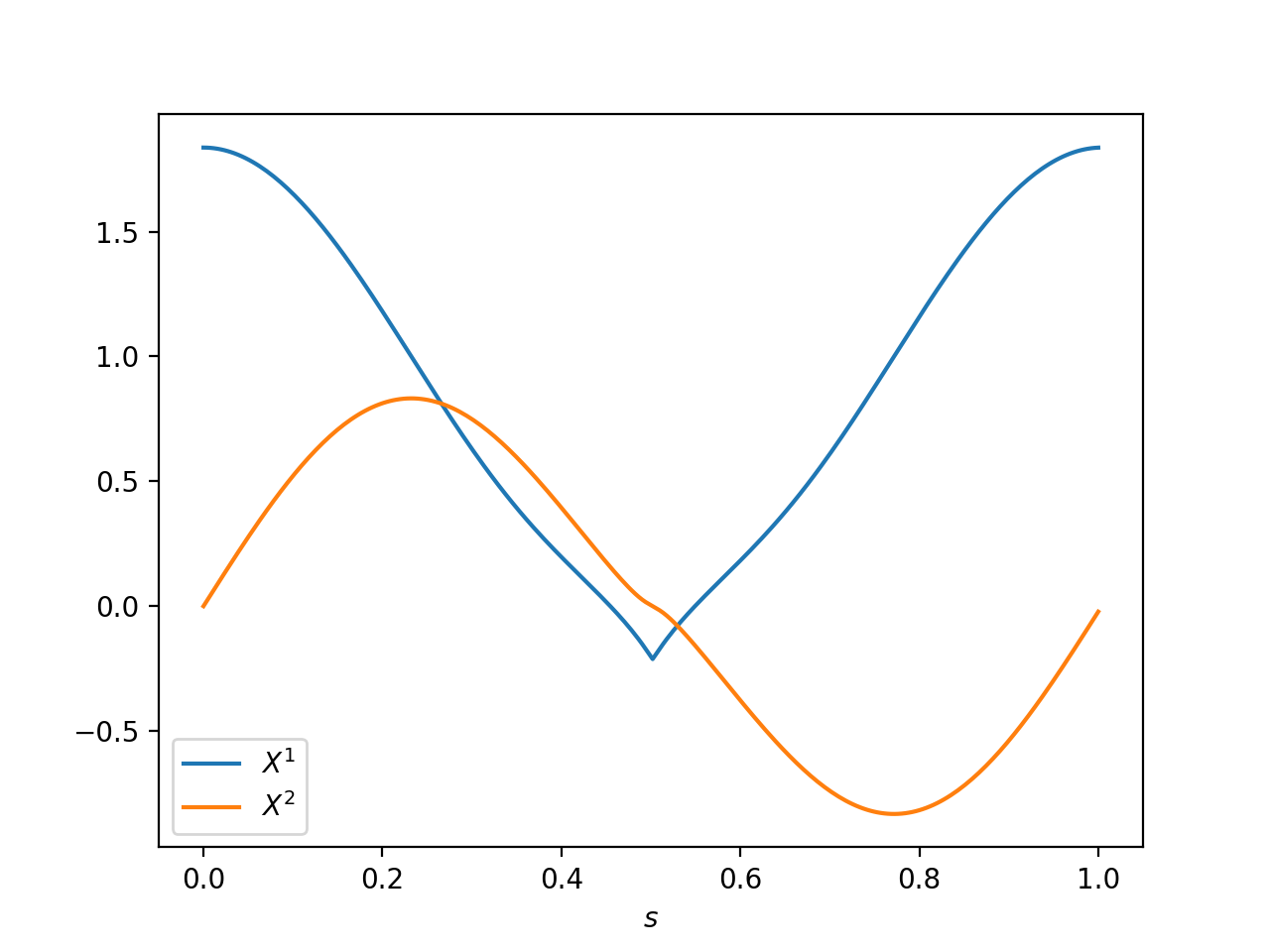}
  \includegraphics[scale=0.4]{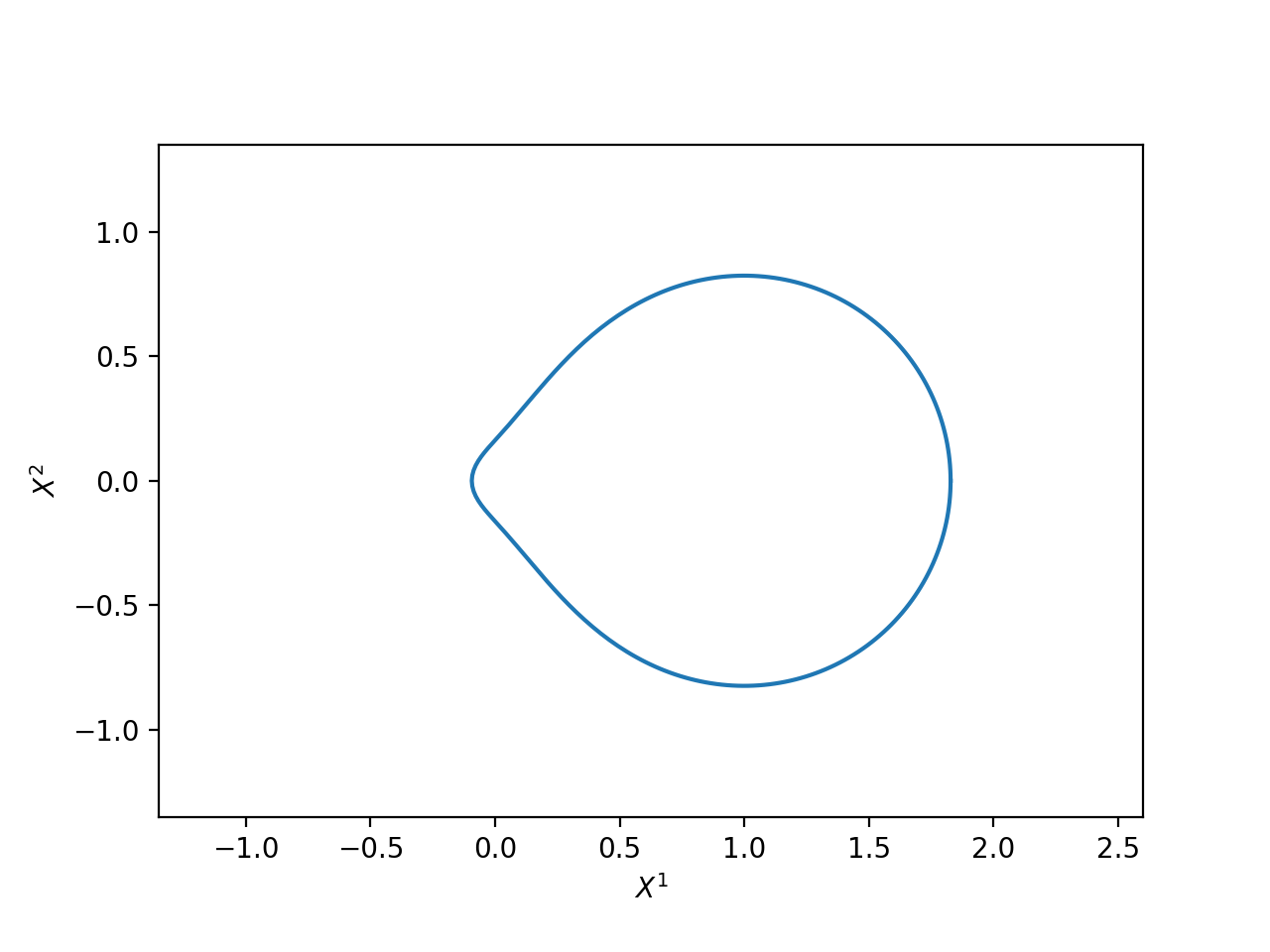}
  \includegraphics[scale=0.4]{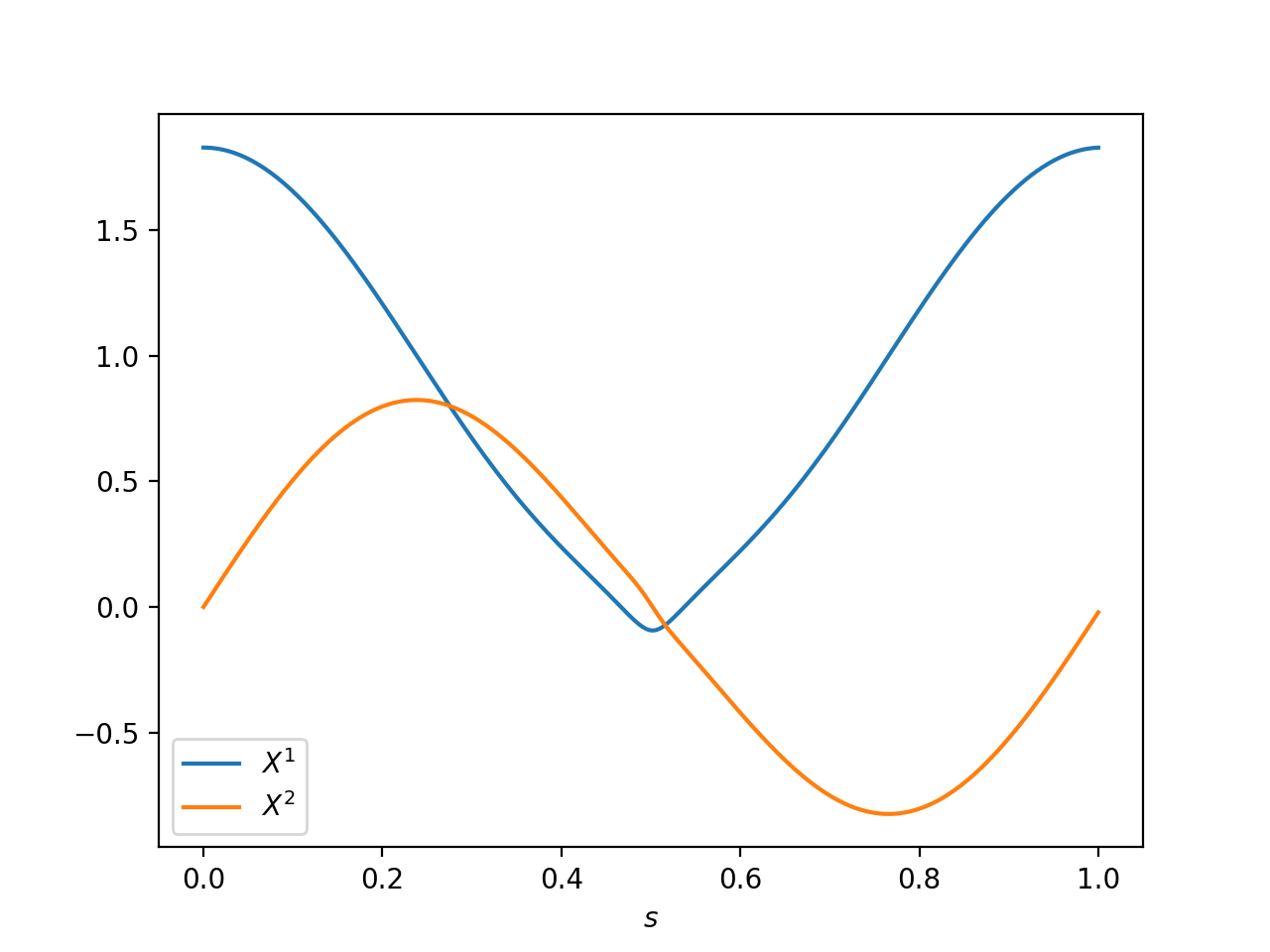}
\caption{Shedding a loop and corresponding singularity. In the first
  row the evolving curve and its (almost normal) coordinate functions
  are depicted. The same are shown very close to the singular time in
  the second row, and after the singularity is regularized by the flow
  in the third row.}
\label{fig:loopShedding}
\end{center}
\end{figure}

During the flow $X^n$ (arbitrary $n\in \mathbb{N}$), an
initial curve is immediately regularized but can later develop
singularities. They, however, must be very special in nature.
\begin{prop}\label{cusps}
Regardless of $n\in \mathbb{N} $, the only possible singularities in
the evolving curve $X^n$ are cusps.
\end{prop}
\begin{proof}
At a singular time $t^n_{sing}$, the function $Z=Y_s(t^n_{sing},\cdot)$ is
analytic. Take a singular point of the corresponding curve. It can be
assumed without loss of generality that the singular point occurs at
the parameter value $s=0$. Then, in a neighborhood, we have that
$$
Z(s)=\sum_{n\geq 1}a_ns^n,
$$
for coefficients $a_n\in \mathbb{R}^d$. Let $p\in \mathbb{N}$ be the
smallest integer for which $a_n\neq 0$. A rigid transformation
$T$ can be found such that $Ta_p=|a_p|e_1$ and thus
$$
 TZ(s)=|a_p|e_1s^p+\sum_{n\geq p+1}Ta_ns^n.
$$
Then we have that
$$
\lim_{s\to 0\pm}\frac{TZ(s)}{|TZ(s)|}=\lim_{s\to
  0\pm}\frac{s^p}{|s|^p}e_1=\begin{cases}
  \pm e_1,&p\text{ odd},\\ e_1,&p\text{ even}.
  \end{cases}
$$
As it is assumed that $s=0$ corresponds to a singular point of the
curve $Z$, $p$ must be odd and the tangent points in direction
$\pm T^{-1}e_1$ on either side of the singular point. This makes the
singular point a cusp. In order to determine the exact type of cusp
one can determine the smallest integer $q>p$ for which $(Ta_n)'\neq
0$, where $a'=(a^2,\dots,a^d)$ for $a\in \mathbb{R}^d$, and a rigid
transformation $\widetilde{T}$ of $\mathbb{R}^{d-1}$ such
that $\widetilde{T}(Ta_q)'=|a_q'|e_2'=|a_q'|e_1\in
\mathbb{R}^{d-1}$. Then the asymptotic properties of the cusp are
determined by the parametrization
$$
\bigl(|a_p|s^p+\sum_{n=p}^q (Ta_n)^1s^n,|a_q'|s^q,0,\dots, 0\bigr),
\: s\simeq 0,
$$
in the new coordinates determined by $T$ and
$\operatorname{id}_\mathbb{R}\otimes
\widetilde{T}$.
\end{proof}
\begin{cor}
Curves with kinks (changes of direction of less than $\pi$ degrees) do
not admit analytic parametrizations.
\end{cor}
\begin{rem}
It cannot be excluded, in general, that a cusp singularity as described
in the above proof be degenerate. This happens when only one component
of the parametrization is non-zero (in the proper coordinates). This
means that the cusp does not open. This is the case, for instance, for
the two cusps of a doubly covered segment. Such a degenerate loop can,
however, not appear on an immersed curve.
\end{rem}
\begin{rem}
It is an interesting observation worth further investigation that, for
the kind of curves considered here, the tangent line (as opposed to
the tangent vector) is well-defined at each point. Indeed a jump of
$\pi$ degrees in the direction of the tangent vector does not
represent a discontinuity in the tangent line (only a discontinuity of
its orientation).
\end{rem}
\begin{prop}
If $X_0$ has at most finitely many extremity points, i.e. if
$\mathcal{E}(X_0)<\infty$, then the flow $X^n$ experiences
singularities in at most finitely many singular times, i.e. times
$t_{sing}$  at which $\Gamma^n(t_{sing})$ has at most finitely many cuspidal
points.
\end{prop}
\begin{proof}
We can assume
that $X_0$ be smooth and normally parametrized.  The initial curve has
at most finitely many extremity points, 
i.e. zeros in any of the components of $\partial_sX_0$. Singularities are
formed when $|\partial _s Y^n(t,\cdot)|$ vanishes at some parameter
value $s_0$ and (singular) time $t_{sing}$. We already observed that
inflection points immediately disappear and cannot therefore be
involved in the lead-up to a singularity. They can appear at a
singularity but not coalesce with other inflection points or
extrema. We therefore need to follow maxima and minima of the
component functions, of which $X_0$ has a finite number and which are,
hence, isolated. As component extrema move around during the evolution
and cannot increase in number, zeros of the $j$th component
$\partial_sY^{n,j}$ will occasionally occur at a common parameter value
$s_0\in[0,1)$ for all $j=1,\dots, d$. The singularity formed will be a
cusp by Proposition \ref{cusps}. At such singular points $\partial_sR \bigl(
Y^n(t_{sing})\bigr)(s_0)$ will experience a U-turn. In the proper
coordinates, a (non-degenerate) cusp corresponds to a minimum in the first component
and an inflection point in the second (see Proof of \ref{cusps} where
it is shown that singularities are asymptotically in a plane). The latter will instantaneously
disappear and $\Gamma^n(t_{sing}+dt)$ will be a smooth curve with at least
one less critical point in its components compared to
$\Gamma^n(t_{sing})$. As the inflection point can 
only have been generated by the merger of extrema, the loss of
extremity points as compared to $\Gamma^n(t_{sing}-dt)$ is at least two. It
may happen that a singularity appears at $t\in h \mathbb{N}$, in
which case the discussion above still applies. There are therefore at
most finitely many singular times during the evolution at each of
which at most a finite number of cusps are formed. Notice that
degenerate cusps cannot occur since they require vanishing of all
derivatives coming into the vertex of the cusp (in the appropriate
coordinates), which is only possible for trivial analytic
parametrizations. Here we make use of Remark \ref{CoC}.
\end{proof}
\begin{rem}
The formation of (non degenerate) cusps, which amount to an extremum
in a component and an inflection point in another (in the proper
coordinates), stems from the merger of two (or more) extrema which
give rise to an inflection point. When a maximum and a minimum merge
(at a singular time), this corresponds to the evolving curve
``shedding'' a loop. This phenomenon is depicted in Figure \ref{fig:loopShedding}.
\end{rem}
\begin{rem}
The above result shows that inflection points (in any of the
components of $X^n$) appear in tandem with singularities, where maxima
and minima merge, and instantaneously disappear. Observe that the main
driver of the evolution is diffusion
(albeit with varying diffusivity) and that the operator $R$ preserves
component minima and maxima in size (while, of course,
distorting the shape of the graphs). This means that the
flow tries to remove component minima and maxima over time. This can
happen without the formation of any singulariy or, in the case of
coalescence, with the formation of one or more singularities. It
follows that there will be an asymptotic number of extremity points as
the extinction time (for the limiting flow) is approached. It can
be seen that this number is 6 for infinity like shapes, 4n for n-fold
circular points, and 8 for an 
Abresch-Langer curve with rotation index 2 and closing up after 3
periods of its curvature function (and hence a perturbation of a
doubly covered circle). Figure \ref{fig:convolutedCurve} depicts the
evolution of a convoluted curve. It is apparent how the number of
extremity points decreases over time with and without singularity
formation. 
\end{rem}

\subsection{The Behavior of Length}
\begin{prop}\label{lendgthDown}
The length is a continuous non-increasing function of time
$$
L\bigl(X^n(t)\bigr)\leq L \bigl( X^n(\tau)\bigr) \leq L \bigl( X^n(0)\bigr)
=L(X_0),\: t\geq \tau\geq 0,\: n\in \mathbb{N}.
$$
\end{prop}
\begin{proof}
In the construction of $X^n$, $X^n(t)$ is obtained either by heat
flowing a previous parametrization (with constant but changing
time-scale/diffusivity) or by reparametrization. As reparametrization
does not change the curve nor its length, it is enough to show that
the heat flow of a parametrization with arbitrary diffusivity
$\alpha>0$ exhibits non-increasing length. In other words, let $Y$
satisfy $Y_t=\alpha Y_{ss}$ and consider $L \bigl( Y(t)\bigr)=\int_0^1
|Y_s(t,\sigma)|\, d \sigma$. Then
\begin{align}\label{lengthDown}
\frac{d}{dt} L \bigl( Y(t)\bigr)&=\int_0^1
\frac{Y_s}{|Y_s|}\cdot Y_{st}\ d \sigma=\alpha \int_0^1
\frac{Y_s\cdot Y_{sss}}{|Y_s|}\, d \sigma\notag\\
&=-\alpha \int_0^1 \bigl( \frac{Y_s}{|Y_s|}\bigr) _s\cdot Y_{ss}\, d
\sigma =-\alpha \int _0^1 \bigl( \frac{Y_s}{|Y_s|}\bigr) _s \cdot
\bigl( Y_{ss}-\frac{Y_{ss}\cdot Y_s}{|Y_s|}\frac{Y_s}{|Y_s|}\bigr)\, d
\sigma \notag\\ &= -\alpha \int _0^1 k(Y)\cdot k(Y)|Y_s|^3\, d \sigma\leq 0,                                                 
\end{align}
where $k(Y)$ is the curvature vector of the curve parametrized by
$Y$. This holds everywhere with the exception of the finitely many
singular times away from which continuity follows. At a singular time $t_{sing}$,
$Y^n(t_{sing})$ is a smooth parametrization of a curve with at most finitely
many cuspidal points. If $t_{sing}\not\in k \mathbb{N}$, then
$$
Y^n\in \operatorname{C}\bigl( [t_{sing}-\delta,
t_{sing}+\delta],\operatorname{C}^1_\pi(I)\bigr), 
$$
for small enough $\delta>0$ as a solution of
$Y_t=\frac{1}{L(X^n(kh))^2}Y_{ss}$ with
$Y(t_{sing}-\delta)=Y^n(t_{sing}-\delta)$ and $t_{sing}\pm d\in\bigl(
kh,(k+1)h\bigr)$ and $\int_0^1|Y^n_s(\cdot,\sigma)|\, d\sigma $ is
therefore continuous at $t$. If, on the other hand, $t_{sing}=(k+1)h$ for
some $k=0,1,\dots$, then left continuity follows as above since
$$
L(X^n(t))=L(Y^n(t))=L \bigl(
e^{\frac{t-kh}{L^2(X^n(kh)}\partial_{ss}}X^n(kh)\bigr), t\in \bigl[ kh,(k+1)h\bigr),
$$
and
$$
X^n((k+1)h)=R \bigl( Y^n((k+1)h)\bigr),\: L \bigl( X^n((k+1)h)\bigr)=L
\bigl( Y^n((k+1)h)\bigr).
$$
For $t>(k+1)h$, we use the fact that $R(Y^n((k+1)h)_s=X^n_s((k+1)h)$ is a continuous
function except at mostly finite many points where it has jump
discontinuities. In particular, it holds that  $X^n_s((k+1)h)\in
\operatorname{L}^\infty_\pi(I)\subset \operatorname{L}^1_\pi(I)$ and
then
\begin{equation*}
\partial_se^{\frac{t-(k+1)h}{L^2(X^n(kh)}\partial_{ss}}
X^n((k+1)h)=e^{ \frac{t-(k+1)h}{L^2(X^n(kh)}\partial_{ss}}
X^n_s((k+1)h)\in\operatorname{C} \Bigl( \bigl[
(k+1)h,(k+2)h\bigr),\operatorname{L}^1_\pi(I)\Bigr),
\end{equation*}
from which the claim follows.
\end{proof}
\begin{rem}
Figure \ref{fig:l+tac} depicts the behavior of
length for the flow depicted in Figure \ref{fig:convolutedCurve} where
the singular times are clearly visible even in the numerical
approximation. The first two occur where the curve looses loops and
the last at the extinction time.
\end{rem}
\begin{prop}
Setting $L_k=L\bigl( X^n(kh)\bigr)$, $k\geq 1$, it holds that
\begin{equation}\label{lenRecursion}\begin{cases}
  L_{k+1}\leq e^{-\frac{h}{L_k^2}4\pi^2}L_k, &k\geq 0,\\
  L_0=L(X_0).
\end{cases}\end{equation}
\end{prop}
\begin{proof}
We compute
\begin{align*}
L_{k+1}^2&=\Bigl(\int_0^1 \big |Y^n_s \bigl( (k+1)h,\sigma)\bigr) \big |\, d
\sigma\Bigr)^2\leq \int_0^1 \big |Y^n_s \bigl( (k+1)h,\sigma)\bigr) \big
|^2\, d \sigma\\ &= \int _0^1 \sum_{j=1}^d \big |Y^{n,j}_s \bigl(
(k+1)h,\sigma\bigr)|^2\, d \sigma = \sum_{j=1}^d \|
Y^{n,j}_s((k+1)h)\|_2^2\\ &\leq e^{-2h\frac{4\pi^2}{L_k^2}}\sum_{j=1}^d
\| X^{n,j}_s(kh)\|_2^2=e^{-2h\frac{4\pi^2}{L_k^2}}\sum_{j=1}^d
\int_0^1 \big |(X^{n,j}_s(kh,\sigma)\big |^2\, d \sigma              
  \\ &= e^{-2h\frac{4\pi^2}{L_k^2}}\int_0^1 \sum_{j=1}^d \big
|(X^{n,j}_s(kh,\sigma)\big |^2\, d \sigma= e^{-2h\frac{4\pi^2}{L_k^2}}L_k^2,
\end{align*}
using the facts that $Y^n_s(t)$ solves $Y_t=\frac{1}{L_k^2}Y_{ss}$ for
$t\in\bigl( kh,(k+1)h\bigr)$ with
$Y(kh)=R(Y^n(kh))_s=X^n_s(kh)$ and that $|X_s^n(kh)|=L_k$ as $X^n(kh)$
is normally parametrized. The claim follows.
\end{proof}
\begin{thm}\label{extComparison}
For any given initial curve $X_0\in \mathcal{C}(\mathbb{R}^d)$ and any
$n\in \mathbb{N}$, it holds that
$$
 L(X^n(t))>0 \text{ for }t\in\bigl[0,\infty\bigr),
$$
for the evolving family of curves $X^n$ with $X^n(0)=X_0$. For each
$n\in \mathbb{N}$, the decrease in length is faster or equal to that
of a circle of radius $\frac{L(X_0)}{2\pi}$.
\end{thm}
\begin{proof}
It follows from careful inspection of the string of inequalities
leading to the recursive estimate \eqref{lenRecursion}, and it will be
computed again in Subsection \ref{circlesSection}, that identity
actually holds when the initial curve is a circle. Taking its radius to be
$\frac{L_0}{2\pi}$, a comparison is obtained for any curve with initial
length $L_0$. The length is hence strictly decreasing and tends to
zero exponentially since the decreasing sequence satisfies
$$
L_{k+1}\leq L(t)\leq L_k,\: t\in \bigl[ kh,(k+1)h\bigr), 
$$
and $(L_k)_{k\in\mathbb{N}}$ decreases exponentially since
$$
\frac{L_k}{L_0}=\prod _{i=1}^{k}\frac{L_i}{L_{i-1}}\leq 
\bigl( e^{-h\frac{4\pi^2}{L_0^2}}\bigr)^k,\: k\in \mathbb{N}.
$$
The length can, however, not vanish in a finite number of steps since
the evolution is driven by linear diffusion on each time interval
between $hk$ and $h(k+1)$ with a non trivial (non constant) initial
datum. 
\end{proof}
The approximant solution has a infinite extinction time
$t^n_e(X_0)=\infty$, regardless of $n\in \mathbb{N}$, and it converges
to a point under the effect of diffusion, which drives each component
of the position vector to a constant. 
\begin{rem}\label{zeros}
Given an initial datum $X_0\in \mathcal{C}(\mathbb{R}^d)$ and $0<t<h$,
$X_s^n(t)$ is an analytic function of the space variable (regardless of $n\in
\mathbb{N}$). It follows that each of its
components is either constant or has at most finitely many zeros. In
the former case, it must have been a constant component for $X_0$ also
and it will remain constant during the whole evolution. We can
therefore assume without loss of generality that the total number of
zeros of all components of $X^n_s(h)$ is finite and we may as well
assume that this is the case for the initial datum itself.
\end{rem}
\begin{rem}
It follows from the construction procedure that the length
reduction is increasingly pronounced as the total length gets smaller
and smaller. This is due to the fact the diffusion coefficient is
quadratically inversely proportional to the length of the curve. As
the approximation parameter $n$ increases, this coefficient 
is updated more often and leads to a faster decrease of length. The
approximating solution is therefore expected to experience singular
points later than any of its limiting curves. We will show rigorously
later that the extinction time $t_e(X_0)<\infty$ for any initial
curve (of finite length as considered here). It, however, cannot be
excluded in general that $t_e(X_0)=0$, as will be shown by example
below. This can happen when the initial curve is not immersed.
\end{rem}
\begin{rem}
The initial datum can be assumed to be essentially parametrized and
$X^n(t)$ is always essentially parametrized as either an analytic
function or a normal length re\-pa\-ra\-me\-tri\-za\-tion of it.
\end{rem}

\subsection{The Behavior of Total Curvature}
Next we study the behavior of total curvature for the approximant and
for smooth solutions of the CSF. For a smooth curve the total absolute
curvature is given by
$$
\kappa_{tot}=\int_0^L |k|\, d\bar s=\frac{1}{L}\int_0^1|k|\, ds,
$$
where $\bar s$ is arc- and $s$ is normal length.
\begin{prop}\label{tcDownApproxFlow}
Away from singular times, the total curvature along the
approximating flow is non-increasing, i.e.
$$
\kappa_{tot}\bigl(X^n(t)\bigr)\leq
\kappa_{tot}\bigl(X^n(\tau)\bigr)\leq \kappa_{tot}\bigl(X^n(0)\bigr)=
\kappa_{tot}(X_0),\: 0\leq\tau \leq t.
$$
\end{prop}
\begin{proof}
The total curvature of a curve is a geometric property that does not
depend on the parametrization and it 
is therefore enough to show that the curve obtained by linear
diffusion of a curve's parametrization does not increase total
curvature. As the diffusivity is kept piecewise constant for the
approximant, we can assume that the diffusivity is 1 and consider a
solution of $Y_t=Y_{ss}$ to show the claim. It holds that 
$$
\kappa_{tot}\bigl( Y)=\int_0^1 \big | \bigl( \frac{Y_s}{|Y_s|}\bigr)
_s\big |\, ds =\int _0^1 |T_s|\, ds,
$$
where $T$ is unit tangent. Notice that, away from the at most finitely
many singularities, $|Y_s|\neq 0$ and $T$ is analytic. As $|T_s|=0$ at a
finite number of parameter values we cannot work with $|T_s|$ directly
but replace this quantity by $a_\varepsilon
=\sqrt{|T_s|^2+\varepsilon}$ and compute
$$
\frac{d}{dt} \int_0^1\sqrt{|T_s|^2+\varepsilon}\, ds=\int _0^1
\frac{T_s}{a_\varepsilon}\cdot\frac{d}{dt} T_s\, ds=-\int_0^1 \bigl(
\frac{T_s}{a_\varepsilon}\bigr) _s\cdot \frac{d}{dt} T\,
ds,
$$
and we proceed with the computation of the various terms
involved. First notice that
$$
\frac{d}{dt} T=\frac{Y_{st}}{|Y_s|}-T\cdot
\frac{Y_{st}}{|Y_s|}T=\frac{Y_{sss}}{|Y_s|}-T\cdot
\frac{Y_{sss}}{|Y_s|}T=:\frac{Y_{sss}^\perp}{|Y_s|},
$$
where the perpendicularity superscript indicates that we are taking
the orthogonal projection onto the space perpendicular to $T$. Next it
can be verified that
$$
\frac{d}{ds} \bigl(\frac{T_s}{a_\varepsilon}\bigr)=
\frac{T_{ss}}{a_\varepsilon}-\frac{1}{a_\varepsilon^3}T_s\cdot T_{ss}T_s,
$$
taking into consideration the definition of $a_\varepsilon$. It holds that
\begin{align*}
  T_{ss}&=\Bigl( \frac{Y_{ss}}{|Y_s|}-T\cdot \frac{Y_{ss}}{|Y_s|} T
  \Bigr)_s=\frac{Y_{sss}}{|Y_s|}-T\cdot\frac{Y_{ss}}{|Y_s|}\frac{Y_{ss}}{|Y_s|}-T_s
  \cdot\frac{Y_{ss}}{|Y_s|} T-T\cdot \frac{Y_{ss}}{|Y_s|}T_s+\\
 &\phantom{=}-T\cdot \Bigl( \frac{Y_{sss}}{|Y_s|}-T\cdot
   \frac{Y_{ss}}{|Y_s|}\frac{Y_{ss}}{|Y_s|}\Bigr) T\\
  &=\frac{Y_{sss}^\perp}{|Y_s|}-T\cdot\frac{Y_{ss}}{|Y_s|}\frac{Y_{ss}^\perp}{|Y_s|}
 -T_s\cdot\frac{Y_{ss}}{|Y_s|} T-T\cdot \frac{Y_{ss}}{|Y_s|}T_s
 \\&=\frac{Y_{sss}^\perp}{|Y_s|}-2T\cdot \frac{Y_{ss}}{|Y_s|}T_s -T_s\cdot\frac{Y_{ss}}{|Y_s|} T.
\end{align*}
Then we have that
\begin{align*}
\bigl(\frac{T_s}{a_\varepsilon}\bigr)_s&=\frac{1}{a_\varepsilon}\Big\{
T_{ss}-\frac{T_s}{a_\varepsilon}\cdot T_{ss}\frac{T_s}{a_\varepsilon}\Big\}
=\frac{1}{a_\varepsilon}\Big\{\frac{Y_{sss}^\perp}{|Y_s|}-2T\cdot
\frac{Y_{ss}}{|Y_s|}T_s-T_s\cdot\frac{Y_{ss}}{|Y_s|} T+\\
&\phantom{=}-\frac{T_s}{a_\varepsilon}\cdot\Bigl(\frac{Y_{sss}^\perp}{|Y_s|}-2T\cdot
\frac{Y_{ss}}{|Y_s|}T_s-T_s\cdot\frac{Y_{ss}}{|Y_s|} T\Bigr) \frac{T_s}{a_\varepsilon}\Big\}\\
&=\frac{1}{a_\varepsilon}\Big\{\frac{Y_{sss}^\perp}{|Y_s|}-\frac{T_s}{a_\varepsilon}\cdot
\frac{Y_{sss}^\perp}{|Y_s|}\frac{T_s}{a_\varepsilon}\Big\}-2\frac{T\cdot Y_{ss}}{|Y_s|}
 \frac{T_s}{a_\varepsilon}-\frac{T_s}{a_\varepsilon}\cdot \frac{Y_{ss}}{|Y_s|}T\\
&\phantom{=}+2 \frac{T_s}{a_\varepsilon}\cdot \frac{T_s}{a_\varepsilon}
\frac{T\cdot Y_{ss}}{|Y_s|}\frac{T_s}{a_\varepsilon}+\frac{T_s}{a_\varepsilon}
\cdot T \frac{T_s}{a_\varepsilon}\cdot \frac{Y_{ss}}{|Y_s|}\frac{T_s}{a_\varepsilon}\\
&=\frac{1}{a_\varepsilon}\Big\{\frac{Y_{sss}^\perp}{|Y_s|}-\frac{T_s}{a_\varepsilon}\cdot
\frac{Y_{sss}^\perp}{|Y_s|}\frac{T_s}{a_\varepsilon}\Big\}-2\frac{T\cdot Y_{ss}}{|Y_s|}
 \frac{T_s}{a_\varepsilon}\Bigl( 1-\frac{T_s}{a_\varepsilon}\cdot \frac{T_s}{a_\varepsilon}\Bigr)
-\frac{T_s}{a_\varepsilon}\cdot \frac{Y_{ss}}{|Y_s|}T
\end{align*}
Combining everything we compute
\begin{equation*}
\bigl( \frac{T_s}{a_\varepsilon}\bigr) _s\cdot \frac{d}{dt}
  T=\frac{1}{a_\varepsilon}\Big\{\frac{|Y_{sss}^\perp|^2}{|Y_s|^2}-
  \bigl(\frac{T_s}{a_\varepsilon}\cdot\frac{Y_{sss}^\perp}{|Y_s|}\bigr)^2\Big\}
 -2\frac{T\cdot Y_{ss}}{|Y_s|}\frac{T_s}{a_\varepsilon}\cdot\frac{Y_{sss}^\perp}{|Y_s|}
\Bigl( 1-\frac{T_s}{a_\varepsilon} \cdot \frac{T_s}{a_\varepsilon}\Bigr)=I-II
\end{equation*}
In order gain insight into the two terms above we write
$u=\frac{Y_{ss}}{|Y_s|}$ and $v=\frac{Y_{sss}}{|Y_s|}$ and work in an
orthonormal basis starting with $e_1=T$, $e_2=\frac{T_s}{|T_s|}$. Notice that
such a basis exists almost everywhere since $T_s$ can only vanish
finitely many times. Then term $I$ reads
$$
\frac{1}{a_\varepsilon}\bigl(
v_2^2+\cdots +v_d^2-\frac{a_0^2}{a_\varepsilon^2}v_2^2\bigr)\geq 0
$$
since $|T_s|=a_0\leq a_\varepsilon$. For $II$ we have
$$
II=-2u_1\frac{a_0}{a_\varepsilon}v_2\bigl(1-\frac{a_0^2}{a_\varepsilon^2}\bigr)
\longrightarrow 0\: (\varepsilon\to 0).
$$
By Lebesgue's dominated convergence theorem we finally arrive at
$$
\int _0^1 |T_s|(t)\, ds-\int_0^1|T_s|(\tau)\, ds\leq 0\text{ for
}t\geq \tau,
$$
where $t$ and $\tau$ are in between singular times. Indeed
$$
\int _0^1 \sqrt{|T_s|^2(t)+\varepsilon}\, ds-\int_0^1\sqrt{|T_s|^2(\tau)+\varepsilon}\, ds
=-\int_\tau^t\int_0^1 \Big\{ I+II \Big\}\leq -\int_\tau^t\int_0^1 II,
$$
and the claim follows from the almost everywhere convergence and
everywhere boundedness of $II$, as well as the convergence of the
terms on the left as $\varepsilon\to 0$.
\end{proof}
\begin{prop}\label{tcDown}
The total curvature is non-increasing for any smooth solution of the
CSF in any space dimension.
\end{prop}
\begin{proof}
We borrow the idea of the proof from \cite{Alt91} where the case of 3D
space curves is considered. For this proof we denote the arclength
parameter by $s$ instead of $\overline{s}$. If $X$ solves the CSF and
is smooth, then 
$$
\kappa(X)=|k(X)|=|X_{ss}|.
$$
First we derive an equation for $|X_{ss}|^2$ and then deal with the
total curvature. It holds that
\begin{align*}
\frac{d}{dt} |X_{ss}|^2&=2 \frac{d}{dt} \frac{d}{ds} \frac{d}{ds}
  X\cdot X_{ss}=2\Bigl[ \frac{d}{ds} \frac{d}{dt} \frac{d}{ds} X
  +|X_{ss}|^2X_{ss}\Bigr]\cdot X_{ss}\\&=2\Bigl[
  \frac{d^2}{ds^2}\frac{d}{dt} X+\frac{d}{ds} \bigl( |X_{ss}|^2X_s\bigr)
  +|X_{ss}|^2X_{ss}\Bigr]\cdot X_{ss}\\&=
2 \frac{d^2}{ds^2}X_{ss}\cdot X_{ss}+2|X_{ss}|^4+4(X_{sss}\cdot X_{ss})
  (X_s\cdot X_{ss})\\&=2 \frac{d^2}{ds^2}X_{ss}\cdot X_{ss}+2|X_{ss}|^4
\end{align*}
since $\frac{d}{dt} \frac{d}{ds} =\frac{d}{ds} \frac{d}{dt} +\kappa^2
\frac{d}{ds}$ and $X_{ss}\cdot X_s\equiv 0$. Next we use that
$$
 \frac{d^2}{ds^2}\bigl( X_{ss}\cdot X_{ss}\bigr)=2\frac{d^2}{ds^2}
 X_{ss}\cdot X_{ss}+2X_{sss}\cdot X_{sss},
$$
to see that
$$
\frac{d}{dt} |X_{ss}|^2=\frac{d^2}{ds^2}|X_{ss}|^2-2X_{sss}\cdot
X_{sss}+2|X_{ss}|^4.
$$
In order to deal with $|X_{ss}|$ we use a ``regularization'' that
Altschuler attributes to Hamilton in \cite{Alt91} and derive first an
equation for $\kappa_\varepsilon=\sqrt{|X_{ss}|^2+\varepsilon}$. It
holds that
$$
\frac{d}{dt} \kappa_\varepsilon =\frac{1}{\kappa_\varepsilon
}\frac{d}{dt} |X_{ss}|^2=\frac{1}{\kappa_\varepsilon}\Bigl[
\frac{1}{2}\frac{d^2}{ds^2}|X_{ss}|^2-X_{sss}\cdot X_{sss}+|X_{ss}|^4\Bigr],
$$
and we compute
\begin{align*}
\frac{d^2}{ds^2}\kappa_\varepsilon &=\frac{d}{ds} \bigl(
\frac{1}{\kappa_\varepsilon }X_{sss}\cdot
X_{ss}\bigr)=-\frac{1}{\kappa_\varepsilon^3}\bigl(X_{sss}\cdot
X_{ss}\bigr)^2+\frac{1}{\kappa_\varepsilon }\bigl( X_{ssss}\cdot
X_{ss}+X_{sss}\cdot X_{sss}\bigr) \\&=
  -\frac{1}{\kappa_\varepsilon^3}\bigl(X_{sss}\cdot
  X_{ss}\bigr)^2+\frac{1}{2\kappa _\varepsilon}\frac{d^2}{ds^2}|X_{ss}|^2,
\end{align*}
thus arriving at
$$
\frac{d}{dt} \kappa_\varepsilon = \frac{d^2}{ds^2}\kappa_\varepsilon
+\frac{1}{\kappa_\varepsilon^3}\bigl(X_{sss}\cdot
X_{ss}\bigr)^2-\frac{1}{\kappa_\varepsilon
}|X_{sss}|^2+\frac{|X_{ss}|^4}{\kappa_\varepsilon}.
$$
Finally we compute
$$
\frac{d}{dt} \int \kappa_\varepsilon\, ds=\int\bigl(\frac{d}{dt}
\kappa_\varepsilon -|X_{ss}|^2\kappa_\varepsilon\bigr) \, ds
$$
and exploit the fact that
$$
\frac{|X_{ss}|^4}{\kappa_\varepsilon}\leq
|X_{ss}|^2\kappa_\varepsilon\text{ and
}\frac{1}{\kappa_\varepsilon^3}\bigl(X_{sss}\cdot 
X_{ss}\bigr)^2\leq \frac{|X_{sss}|^2}{\kappa_\varepsilon},
$$
as follows from $|X_{ss}|\leq\kappa_\varepsilon $, to obain that
$$
\frac{d}{dt} \int \kappa_\varepsilon ds\leq 0,
$$
which yields the claim in the limit as $\varepsilon\to 0$.
\end{proof}
\begin{rem}
In \cite{Alt91} the Frenet-Serret formul{\ae} for space curves ($d=3$)
are used and the decay of the total curvature is measured in terms of
the curvature and of the torsion $\tau$ to give the stronger estimate
$$
 \frac{d}{dt} \int |k|\, ds\leq -\int \tau^2|k|\, ds.
$$
This estimate is then used in the same paper to show that singularity
formation for space curves is a planar phenomenon. 
\end{rem}
\begin{rem}
During the flow, some curvature can occasionally (a finite number of
times) concentrate (blow up in a controlled way) at finitely many
singular points but, as a measure, it is always defined. In this sense
curvature never blows up during the evolution (and, in the limit as
$n\to\infty$, blows up only once and this happens at the extinction
time).
\end{rem}
\begin{rem}
When the evolving curve develops one or more of its finite
singularities, it is in the form of a cusp. In such a situation the
total (angular) curvature has a jump (in time) where it decreases
instantaneously by $\pi$ for each cusp that is formed. Notice also
that, thanks to Proposition \ref{cusps}, the (spatial) singularities
are asymptotically two dimensional (regular cusps) or one dimensional
(degenerate cusps). 
\end{rem}

\section{Existence of a Limiting Flow}

\begin{thm}\label{exThm}
Let $X_0\in \mathcal{IC}_e(\mathbb{R}^d)$ be given such that
$\kappa_{tot}(\Gamma_0)<\infty$. Then the sequence
$(X^n)_{n\in\mathbb{N}}$ converges uniformly to a solution
$X:[0,\infty)\to \mathbb{R}^d$ of the CSF \eqref{csf}. The solution 
$X$ has finite extinction time $t_e(X_0)$, which is bounded above by
the extinction time of a circle of radius $\frac{L(X_0)}{2\pi}$.

If the initial datum has a finite number of extremity points, the limit has
at most a finite number of singular times away from which it is
smooth (in space and time) and hence satisfies the equation in the
classical sense.
\end{thm}

\begin{proof}
We assume that $t_e(X_0)>0$ because in the degenerate situation that
$t_e(X_0)=0$ there is nothing to prove. When $d=2$ this is always the
case as soon as $X_0$ has a loop that contains a circle.\\[0.15cm]
(i) Propositions \ref{lengthDown} and \ref{tcDownApproxFlow} ensure that
$$
X^n\in \operatorname{W}^1_\infty \bigl(
[0,\infty),E_0\bigr)\cap\operatorname{L}^\infty \bigl(
[0,\infty),E_1\bigr), 
$$
with a bound on the norm which is independent of $n\in \mathbb{N}$. Here we set
$$
E_0=\mathcal{M}_\pi(I)\text{ and }E_1=\big\{ Y\in
\operatorname{W}^1_{\infty,\pi}(I)\,\big |\, Y_s\in
\operatorname{BV}_\pi(I)\big\}.
$$
Observing that
$\operatorname{W}^s_{p',\pi}(I)\overset{d}{\hookrightarrow}
\operatorname{C}_\pi(I)$ for $s>1/p'$ and $p'\in(1,\infty)$, we infer
the validity of the embedding
$E_0=\mathcal{M}_\pi(I)\hookrightarrow
\operatorname{W}^{-s}_{p,\pi}(I)=:F_0$. Notice that
$E_1\hookrightarrow\operatorname{W}^1_{p,\pi}(I)=:F_1$ for every $p>1$.
Then
\begin{equation}\label{emb}
\operatorname{W}^1_\infty \bigl(
[0,\infty),E_0\bigr)\cap\operatorname{L}^\infty \bigl(
[0,\infty),E_1\bigr)\hookrightarrow \operatorname{W}^1_\infty \bigl(
[0,\infty),F_0\bigr)\cap\operatorname{L}^\infty \bigl(
[0,\infty),F_1\bigr).
\end{equation}
Using the real interpolation functor we see that
$$
(F_0,F_1)_{\theta,p}=\operatorname{W}^{\theta-s(1-\theta)}_{p,\pi}(I),\: \theta\in(0,1),
$$
and, consequently that
$$
\operatorname{W}^1_\infty \bigl(
[0,\infty),E_0\bigr)\cap\operatorname{L}^\infty \bigl(
[0,\infty),E_1\bigr)\hookrightarrow \operatorname{C}^{1-\theta}\bigl(
[0,\infty), \operatorname{W}^{\theta-s(1-\theta )}_{p,\pi}(I)\bigr),
$$
using \eqref{emb}, the interpolation inequality for $(F_0,F_1)_{\theta,p}$, and
the fundamental theorem of calculus. By choosing
$1>\theta>\frac{1+s-\delta}{1+s}$, it can be achieved that
$\theta-s(1-\theta)>1-\delta$ for any small $\delta>0$. 
Thus, for any $\delta\in(0,1)$, there is $r>0$ (small) such that
$$
X^n\in\operatorname{C}^r\bigl([0,\infty),\operatorname{W}^{1-\delta}_{p,\pi}(I)\bigr)
$$
with a bound that is uniform in $n\in\mathbb{N}$. As the first
embedding below is compact,
$$
\operatorname{C}^r\bigl([0,T],\operatorname{W}^{1-\delta}_{p,\pi}(I)\bigr)
\hookrightarrow \operatorname{C}^{\tilde
  r}\bigl([0,T],\operatorname{W}^{1-\tilde\delta}_{p,\pi}(I)\bigr)
\hookrightarrow \operatorname{C}^\rho([0,T]\times I),\:
\tilde r<r, \tilde \delta >\delta,
$$
and the second is valid for $\rho>0$ small enough, it is possible to
extract a subsequence that converges in the second space for any $T>0$
to a limiting function $X$ that is an element of the first space and for
which convergence is uniform in space and time due the second
embedding. Notice that $X(T,\cdot)$ must have finite length by lower
semicontinuity of length. If
$L\bigl( X(T,\cdot)\bigr)>0$, then, given 
$0<\varepsilon<L\bigl( X(T,\cdot)\bigr)$, it is possible to  find
$m\in \mathbb{N}$ and a partition $(s^m_j)_{j=0}^{m-1}$ of parameter
space such that 
$$
 L\bigl( X(T,\cdot)\bigr)-\frac{\varepsilon}{2}\leq \sum_{j=0}^{m-1}
 \big |X(T,s_{(j+1)\operatorname{mod}m})-X(T,s_{j})\big |,
$$
since the length of $X(T,\cdot)$ is the supremum of the above expression over all partitions. Due to
the uniform convergence of $X^n(T,\cdot)$, $N\in \mathbb{N}$ can be found such that
\begin{equation}\label{lenComp}
L\bigl(X^n(T,\cdot)\bigr)\geq \sum_{j=0}^{m-1}
 \big |X^n(T,s_{(j+1)\operatorname{mod}m})-X^n(T,s_j)\big |\geq
 L\bigl(X(T)\bigr)-\varepsilon \text{ for }n\geq N.
\end{equation}
(ii) The approximants $X^n$ are smooth for $0<t<t^n_{1,sing}(X_0)$
with respect to the spatial variable due to the smoothing properties
of the heat equation and with respect of the time variable with the
exception of the discrete time of juncture. Here $t^n_{1,sing}$ is the
first (if any) singular time in their evolution, and
$$
t^\infty_{1,sing}=\liminf_{n\to\infty}t^n_{1,sing}(X_0)>0,
$$
since isolated extremal points cannot merge instantaneously by local
well-posedness of the CSF. As a consequence $\big |X^n_s(t,\cdot)\big
|$ is bounded away from zero independly of $n$ and the
reparametrizations at juncture times ($t=kh$, $k\in \mathbb{N}$)
preserve the regularity (as they does not change the curve\footnote{It
  is enough to exclude the possibility of a smoothly parametrized
  singular curve.}). The
bounding constant for spatial norms depends only on the initial datum
$X_0$ on any compact subinterval 
of $(0,t^\infty_{1,sing})$. As with time, the dependence is at most
Lipschitz due to the junctures. Thus, for any compact subinterval $J$
of $(0,t^\infty_{1,sing})$, we have that $X^n$ is bounded in
$\operatorname{C}^{1-}\bigl( J,C^{2+\varepsilon}_\pi(I)\bigr)$. It is
therefore possible to go to limit (along a further subsequence) in
$$
X^n_t(t)=\frac{1}{L^2\bigl(X^n(k_th)\bigr)}\partial_{ss}
e^{\frac{t-k_th}{L^2(X^n(k_th))}\partial_{ss}}R(X^n(k_th)),\: t\in
\bigl[k_th,(k_t+1)h\bigr),
$$
to obtain that $X_t=\frac{1}{L^2(X)}\partial_{ss} R(X)$ for
$t\in(0,t^\infty_{1,sing})$ as $k_th\to t$ as $h\to 0$ since $X_t$ can
only be the limit of the right-hand-side by closure. Notice that,
while the time derivative experiences jumps at the juncture time for
$X^n$, these gaps closes up in the limit as left and right limits
coincide in the limit due to the strong spatial convergence and time
continuity.
It follows that $X$ is a smooth solution of the CSF on an open
interval. It may develop one or more singularities at some 
time $t_{1,sing}>0$ but, even then, we can use Proposition \ref{tcDown} to
see that it remains a finite length curve of finite total curvature in
the limit. The limiting curve at the singular time is given by
$X(t_{1,sing})$ since $X$ exists for all times. It is possible that
$t_{1,sing}=t_e(X_0)$.\\[0.15cm] 
(iii) If $t_{1,sing}<t_e(X_0)$, it follows from (ii) that $X(t_{1,sing})$ is a curve of finite total
curvature and it clearly still has a finite number of extremity
points. We can therefore argue as above with $X(t_{1,sing})$ as a new
initial datum and obtain that $X$ is a smooth solution of the CSF for
$t\in(t_{1,sing},t_{2,sing})$ and some $t_{2,sing}>t_{1,sing}$. Again, if
$t_{2,sing}<t_e(X_0)$, the argument can be repeated. We end up with a
solution $X$ that is smooth up to a finite number of singular times where
it exhibits a finite number of singular points where a non-trivial
amount of total curvature concentrates (more precisely total curvature
has a density with a singular component concentrated at the singular
points). When $t_e(X_0)$ is reached, the curvature truly blows up as
the length vanishes (with total curvature settling on a limit).\\[0.15cm]
(iv) We already know that the length of $X$ depends
continuously on time and it is approached uniformly by $L(X^n)$ on compact 
time intervals away from singular times.  Using Proposition \ref{lenRecursion} it is seen that
$$
L_{k_{t+dt}}\leq e^{-4\pi^2h(\frac{1}{L^2_{k_{t+dt}-1}}+\dots
  +\frac{1}{L^2_{k_t}})}L_{k_t}
$$
where $t\in \bigl[ k_th,(k_t+1)h\bigr)$ for any $t$, and where $dt>0$.
Notice that $L_k$ also depends on $n\in \mathbb{N}$ but we omit the $n$
dependence for readability. Letting $n\to\infty$ ($h\to 0$), we obtain that 
$$
\frac{L(t+dt)-L(t)}{dt}\leq
\frac{e^{-4\pi^2\int_t^{t+dt} \frac{1}{L^2(\tau)}\, d\tau}-1}{dt}L(t),\: t>0.
$$
As this is valid regardless of the choice of $t>0$ and $dt>0$ as long
as $L(t+dt)>0$ we let $dt\to 0$ to see that
$$
\limsup_{dt\to 0}\frac{L(t+dt)-L(t)}{dt}\leq -\frac{4\pi^2}{L(t)}.
$$
As the limiting curve flow has at most countably many points of
discontinuity (finitely many in fact), the above limit coincides with
the actual derivative a.e. and, integrating, we see that
$$
L(t)\leq \sqrt{L_0^2-8\pi^2 t},
$$
for as long as the length remains positive. Thus the solution becomes
a point in finite time $t_e(X_0)<\infty$. This extinction time is at
most the time it takes a circle of radius $\frac{L_0}{2\pi}$ to shrink
to a point, which is precisely given by the right-hand expression in
the inequality as is explicitly computed in the next section.
\end{proof}

\begin{rem}\label{degenSol}
The construction of a solution works also for non-immersed initial
data in $\mathcal{C}_e(\mathbb{R}^d)$. This is a case in which it is
possible that $t_e(X_0)=0$. We conjecture that this happens only for
defomations of doubly (or multiply) covered segments (including
``curved'' ones). Such curves, in the doubly covered case, are
characterized by the symmetry given by
$$
X_0(\frac{1}{2}+s)=X_0(\frac{1}{2}-s),\: s\in[0,1),
$$
of their normal parametrization, i.e. curves that start in one (end)point,
reach another, and return along the same path to the starting
point. The case of a doubly covered segment is discussed in the next
section.
\end{rem}
\begin{rem}\label{uniqueness}
When $d=2$, using the same proof found in the paper \cite{Ang91} by
Angenent based on the maximum principle, one obtains unique local
existence in time for the solution to the CSF to any initial datum $X_0\in
\mathcal{IC}_e(\mathbb{R}^d)$ for which $\partial_{ss}R(X_0)$ is a
Radon measure or, more specially, has a singular part consisting
of Dirac measures only. Uniqueness in general dimensions holds in
as far it does for local in time classical solutions with initial data of the
regularity considered here.
\end{rem}
\begin{rem}
The solution construction procedure highlights the fact that the
singularities incurred during the flow are due to diffusion and not to
the nonlinear nature of the equation. The reparametrization operator
$R$ does not engender any geometric singularity and the decreasing
length only modulates the diffusive strength causing extinction to
occur in finite rather than infinite time. Diffusion is also
responsible for driving each component of 
$X^n(t)$ to a constant, i.e. for eventual convergence to a point
independently of the initial datum.
\end{rem}
\begin{prop}
At a singularity time $t_{sing}\leq t_e(X_0)$, $X(t_{sing})(I)$ either is a point (if
$t_{sing}=t_e(X_0)$) or it possesses one or more cuspidal points. A cusp is
formed when the evolving curve looses a loop which shrinks to the
corresponding cuspidal point. 
\end{prop}
\begin{rem}
When $d=2$, a smooth embedded curve does consist of a single loop and
$X^n$  does not experience any singularity during its entire evolution
and, in the limit, until its extinction time.
\end{rem}
\begin{proof}
For any fixed $n$, a (cusp) singularity can only form by the merger of
extrema, i.e. from the disappearance of a loop. Thus no singularity can
form since the initial curve consists of a single loop and its disappearance
would signify collapsing to a point in finite time. Taking a time
$0<t<t_e(X^0)$, any finite regularity space norm of $X^n$ and its
Lipschitz time norm can be bounded on $[0,t]\times I$ since $|Y^n_s|$
is never zero and thus bounded away from zero on $[0,t]\times I$ thanks to
\eqref{lenComp} and $L \bigl( X(t,\cdot)\bigr)>0$. Thus, on the same interval, 
$X=\lim_{n\to\infty} X^n$ is a classical solution (as argued in the
proof of Theorem \ref{exThm}) and does not develop
singularities. As $t<t_e(X^0)$ is arbitrary, the claim follows.
\end{proof}

\section{Special Solutions}
Next we study the constructed solution to the CSF in three special
cases: circles, a doubly covered segment, and a class of figure
infinity shapes (also known as figure eight curves in the literature).
\subsection{Circles}\label{circlesSection}
We already know that the only self-similarly evolving curves that
remain normally parametrized for all times are circles (possibly
multiply covered) in a plane. It is enough to consider the simply
covered case where we can assume without loss of generality that
$$
X_0=\bigl(\cos(2\pi s),\sin(2\pi s),0,\dots,0\bigr),\: s\in [0,1],
$$
so that the problem can be formulated in a plane and the additional
coordinates neglected. Clearly $R(X_0)=X_0$ so that
$$
X^n(h)=e^{\frac{h}{L_0^2}\partial_{ss}}X_0=
e^{-\frac{h}{4\pi^2}4\pi^2}X_0=e^{-h}X_0\text{ for }L_0=2\pi, 
$$
and similarly
$$
X^n(2h)=^{\frac{h}{L_1^2}\partial_{ss}}X^n(h)=e^{-\frac{h}{e^{-2h}}}X^n(h)\text{
  for }L_1=L(X^n(h)=2\pi e^{-h},
$$
so that the recursion
$$
X^n \bigl( (k+1)h\bigr)=e^{-\frac{h}{L_k^2}\partial_{ss}}X^n(kh)\text{ for
}L_k=L\bigl(X^n(kh)\bigr),\: k\geq 2,
$$
is obtained. It is, however, clearly enough to follow the length as
the shape does not change, in which case we arrive at
$$\begin{cases}
  L_{k+1}=e^{-\frac{4\pi^2 h}{L_k^2}}L_k,&k\geq 0\\
  L_0=2\pi.&
\end{cases}
$$
Then
$$
L_k=e^{-\frac{t}{n}\frac{4\pi^2}{L_{k-1}^2}-\cdots-\frac{t}{n}\frac{4\pi^2}{L_0^2}},
$$
which is a discretization of
$$
L(t)=e^{-\int_0^t \frac{4\pi^2}{L^2(\tau)}\, d\tau}L_0,
$$
and thus of the ODE
$$
\dot L(t)=-\frac{4\pi^2}{L^2(t)}e^{-\int_0^t
\frac{4\pi^2}{L^2(\tau)}\, d\tau}L_0=-\frac{4\pi^2}{L(t)},\: L(0)=L_0,
$$
which we know to be the exact evolution of the length under the CSF in
this case, as follows from \eqref{noReparam}, for instance.
\section{Doubly Covered Segments}
Again we can assume to be working in the plane and start with a smooth
parametrization of the doubly covered segment
$[-\frac{L_0}{4},\frac{L_0}{4}]$, which has length $L_0$, given by
$$
X_0(r)=\bigl( \frac{L_0}{4}\cos(2\pi r),0\bigr),\: r\in [0,1). 
$$
We parametrize in this way to underscore the connection to the first
non-trivial mode of $-\partial_{ss}$ and to suggest that its evolution
may be thought of as a degenerate self-similar flow. It holds that
$$
R(X_0)(s)=
\begin{cases}
  \bigl(-\frac{L_0}{4}+L_0s,0\bigr),&s\in[0,\frac{1}{2}),\\
  \bigl(\frac{3}{4}L_0-L_0s,0\bigr),&s\in[\frac{1}{2},1).
\end{cases}
$$
Notice that this curve has infinitely many extremity points and does
therefore not satisfy the assumption made earlier, if considering it a
curve in the plane or higher dimensional space. Confined to the one
dimensional line it spans, however, it has a finite number of
extremity points. The construction of the solution can be carried
out. As computed previously, we have that
$$
\partial_{ss}R(X_0)=\bigl( 2L_0(\delta_0-\delta_{\frac{1}{2}}),0\bigr).
$$
In order to obtain $e^{\frac{h}{L_0^2}\partial_{ss}}R(X_0)$ we solve
$$
\begin{cases}
  Y_t=\frac{1}{L_0^2}Y_{ss},&\\
  Y(0)=2L_0 \partial _{ss}^{-1}(\delta_0-\delta_{\frac{1}{2}}),&
\end{cases}
$$
where the inverse $\partial_{ss}^{-1}$ is taken in the space of mean
zero functions. It follows that
$$
\widehat{Y(0)}_\ell=\frac{2L_0}{4\pi^2\ell^2}\Bigl( e^{-2\pi i \ell
  \frac{1}{2}}-e^{-2\pi i \ell 0}\Bigr) = \frac{L_0}{2\pi^2\ell^2} \bigl(
(-1)^\ell-1\bigr)=\begin{cases} 0,&\ell\text{ even,}\\
  -\frac{L_0}{\pi^2\ell^2},&\ell\text{ odd}, \end{cases}
$$
so that
$$
Y^n(h)=-\sum_{\ell\text{ odd}}\frac{L_0}{\pi^2
 \ell^2}e^{-\frac{4\pi^2}{L_0^2}l^2 h}e^{2\pi i ls},
$$
and then
$$
L_1=L_0\sum_{\ell\text{ odd}}\frac{2}{\pi^2
  \ell^2}e^{-\frac{4\pi^2}{L_0^2}l^2h}\leq L_0
\frac{\pi^2}{8}\frac{2}{\pi^2}=\frac{1}{4}L_0,
$$
using the fact that $\sum_{\ell\text{
    odd}}\frac{1}{\ell^2}=\frac{\pi^2}{8}$. This remains true in any
step of the approximate evolution and it therefore holds that
$$
L_n(t)=L_n\leq \frac{1}{4^n}L_0\to 0\text{ as }n\to\infty,
$$
for $h=\frac{t}{n}$. Thus the limiting solution of the curve
shortening flow instantaneously collapses to a point, i.e.
$$
X(t)=\begin{cases} X_0,&t=0,\\ \{ (0,0)\},& t>0.\end{cases}
$$
\begin{rem}
A doubly covered segment can be viewed as the limit of an ellipse with
vanishing aspect ratio such as the one satisfying
$$
(\frac{2x}{L_0})^2+(\frac{2y}{\varepsilon})^2=1.
$$
As the area of this simple curve is reduced at the constant rate
$2\pi$ by the CSF, it limits to a circular point at $t_e=\frac{l
  \varepsilon}{8}$. The doubly covered segment inside shrinks
faster than the ellipse does, regardless of $\varepsilon>0$. It must therefore
vanish instantaneously in line with what we proved above.
\end{rem}
\begin{rem}
A doubly covered segment can also be thought of as a degenerate
loop. We conjecture that any smooth initial datum with added
``hairs'', i.e. with added degenerate loops, instantaneously absorbes
all of them and evolves just as if it did not have any
hair. Numerical experiments (not shown) and the behavior of the doubly
covered segment support this conjecture.
\end{rem}
\section{The Figure Infinity Shape}
\begin{figure} 
\begin{center}
\includegraphics[scale=.5]{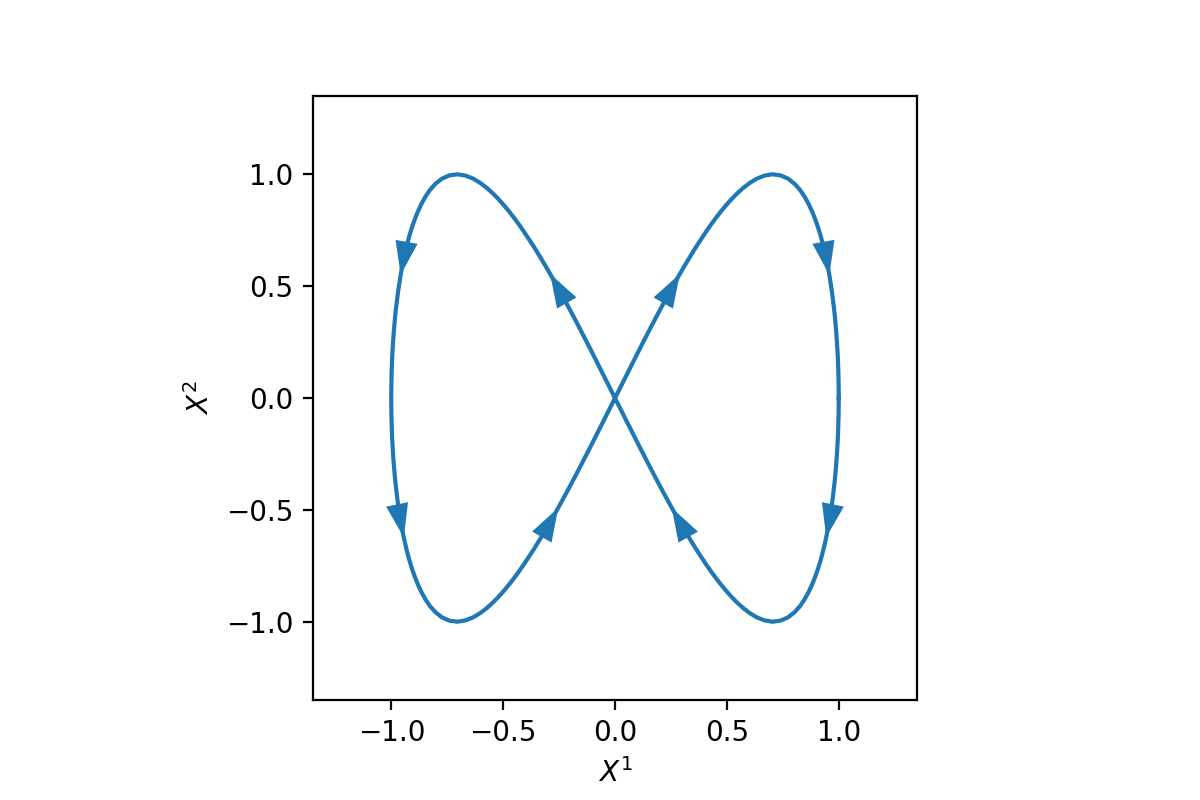}
\caption{A symmetric figure infinity (or figure eight).}
\label{fig:infinityLoop}
\end{center}
\end{figure}

Consider any closed curve with the symmetries of the planar shape $\infty$
with the origin in the self-intersection point and that is specularly
symmetric with respect to both axes. See Figure
\ref{fig:infinityLoop}. Such a shape will have a
parametrization $X_0$ that we can assume to be normal and to start at
the extreme point to the right. We shall also assume that the curve
self-intersects at the origin transversally (as opposed to the
corresponding ``almost embedded singular curve'' which merely self-touches
in the origin). An example of a curve with the desired properties but
not normally parametrized is given by
$$
 X_0= \bigl( \cos(r -\frac{\pi}{2}),\sin(\pi-2r)\bigr),\: r\in [0,2\pi).
$$
Notice that the curve starts at the origin and loops back to it
visiting the first and last quadrant before entering the second
quadrant and looping back to the origin through the third. 
The imposed symmetries imply that
\begin{equation}\label{x-sym}
X_0^1(\frac{k}{4}+s)=(-1)^{k}X_0^1(\frac{k}{4}-s)\text{ for
}s\in[0,1)\text{ and }k=1,\dots,4,
\end{equation}
and that
\begin{equation}\label{y-sym}
X_0^2(\frac{k}{8}+s)=(-1)^{k+1}X_0^2(\frac{k}{8}-s)\text{ for
}s\in[0,1)\text{ and }k=1,\dots,8.
\end{equation}
These symmetries are preserved for all times along any approximant
$X^n(t)$ by the unique solvability of the heat equation and by the
fact that they are clearly preserved when reparametrizing by an
application of $R$ in the sense that, while formulated in terms of the
normal parametrization, they are geometric in nature, i.e. symmetries of the
image set. As $n$ is arbitrary the limiting CSF solution $X$
will also enjoy these symmetries due to uniform convergence. This
means that the signed enclosed
area always vanishes (and hence is independent of the direction of the
parametrization) since it does for $X_0$ as a consequence of the
shapes' symmetries. This already excludes that this curve converge to
a circular point, which has positive enclosed area that only vanishes
at the extinction time. Next observe that the symmetries of the component
functions $X^1$ and $X^2$ also imply that, for every $t>0$,
\begin{gather*}
\int_0^1X^1(t,s)\sin(2\pi s)\, ds=0,\: \int_0^1 X^1(t,s)\cos(2\pi
s)\, ds\neq0\text{ for every }t\in [0,t_e),\\
\int_0^1X^2(t,s)\varphi(s)\, ds=0\text{ for every }t\in[0,t_e)\text{
and } \varphi=\cos(2\pi\cdot), \sin(2\pi\cdot),\\
\int_0^1X^2(t,s)\sin(4\pi s)\, ds\neq 0\text{ for every }t\in [0,t_e),\\
\end{gather*}
Notice how the first component will have a non-zero first cosine mode,
while it will vanish for the second component. For the latter, the
second sine mode will not vanish. This means that, for all
approximants, the first component will decay at a exponential rate of
at least $-\frac{4\pi^2}{L^2}$ where $L$ is the appropriate length for the
corresponding time interval. The second component, however, will decay
with an exponential rate of $-\frac{16\pi^2}{L^2}$ at least. It follows that,
independently of the approximant, the aspect ratio
$\frac{\operatorname{diam}_x\bigl(X^n(t)\bigr)}{\operatorname{diam}_y\bigl(X^n(t)\bigr)}$
of the evolving shape will converge to zero as the extinction time is approached.
\begin{thm}
Let $X_0$ be an infinity shape with the above symmetries and a
transversal self-intersection. Then it will not limit to a circular
point, but rather to a zero aspect ratio doubly covered segment-like
point.
\end{thm}
\begin{proof}
First observe that the symmetries assumed on $X_0$ translate into
symmetries of the normal and curvature vectors and thus any solution
of the CSF (as well as of the approximating flow) will preserve
them. A first consequence of these symmetries is that the evolving
curve will bound a region of vanishing signed area centered in the
origin. Indeed, it follows from
$$
Y(\frac{1}{2}+s)=-Y(\frac{1}{2}-s),\: s\in [0,\frac{1}{2}),
$$
that
\begin{align*}
\int_\Gamma x^1dx^2&=\int ^{\frac{1}{2}}_0Y^1(s)\frac{d}{ds}Y^2(s)\,
ds+\int _{\frac{1}{2}}^1 Y^1(s)\frac{d}{ds}Y^2(s)\, ds\\
&=\int ^{\frac{1}{2}}_0Y^1(s)\frac{d}{ds}Y^2(s)\, ds+
\int ^{\frac{1}{2}}_0Y^1(\frac{1}{2}+s)\frac{d}{ds}Y^2(\frac{1}{2}+s)\, ds\\
&=\int ^{\frac{1}{2}}_0Y^1(s)\frac{d}{ds}Y^2(s)\, ds-
\int ^{\frac{1}{2}}_0Y^1(\frac{1}{2}-s)\frac{d}{ds}Y^2(\frac{1}{2}-s)\, ds\\
&=\int ^{\frac{1}{2}}_0Y^1(s)\frac{d}{ds}Y^2(s)\, ds-
\int ^{\frac{1}{2}}_0Y^1(s)\frac{d}{ds}Y^2(s)\, ds=0,
\end{align*}
where $\Gamma=Y \bigl( [0,1]\bigr)$ is any curve with the discussed
symmetries. Thus the signed area indeed vanishes along the
(approximate) CSF flow and this alone excludes the possibility that
the curve shrinks to a round point. If that were the case, the area
would have to be positive for a non-empty interval of time before
extinction. Now notice that
\begin{align*}
\cos\bigl( 2\pi(\frac{k}{4}+s)\bigr)&=(-1)^k\cos\bigl(
2\pi(\frac{k}{4}-s)\bigr)\\
\sin\bigl( 2\pi(\frac{k}{4}+s)\bigr)&=(-1)^{k+1}\sin\bigl(
2\pi(\frac{k}{4}-s)\bigr).
\end{align*}
Using \eqref{x-sym}, we see that
\begin{align*}
\int_0^1 X^1_0(s)\sin(2\pi s)\,
  ds&=\sum_{k=0}^3\int_{\frac{k}{4}}^{\frac{k+1}{4}}X^1_0(s)\sin(2\pi s)\,
  ds=\sum_{k=0}^3\int_0^{\frac{1}{4}}X^1_0(\frac{k}{4}+s)\sin\bigl(2\pi
  (\frac{k}{4}+s)\bigr)\,
  ds\\&=(-1)^{2k+1}\sum_{k=0}^3\int_0^{\frac{1}{4}}X^1_0(\frac{k}{4}-s)\sin\bigl(2\pi 
  (\frac{k}{4}-s)\bigr)\, ds\\&=\sum_{k=0}^3\int^0_{\frac{1}{4}}X^1_0(\frac{k}{4}+s)\sin\bigl(2\pi
  (\frac{k}{4}+s)\bigr)\, ds\\&=-\int_0^1 X^1_0(s)\sin(2\pi s)\, ds
\end{align*}
and thus
$$
\int_0^1 X^1_0(s)\sin(2\pi s)\, ds=0.
$$
Next observe that
$$
X^1_0(s)> 0\text{ on }(0,\frac{1}{2})\text{ and }X^1_0(s)<
0\text{ on }(\frac{1}{2},1),
$$
and that
$$
\cos(2\pi s )>0\text{ on }(0, \frac{1}{2})\text{ and }\cos(2\pi s
)< 0\text{ on }(\frac{1}{2},1).
$$
It follows that
$$
\int_0^1 X^1(t,s)\cos(2\pi s)\, ds=\int_0^{\frac{1}{2}}
X^1(t,s)\cos(2\pi s)\, ds+ \int^1_{\frac{1}{2}} X^1(t,s)\cos(2\pi s)\,
ds>0.
$$
As for $X^2_0$ we have that
\begin{align*}
\int_0^1 X^2_0(s)\sin(2\pi s)\,
  ds&=\sum_{k=0}^3\int_{\frac{k}{4}}^{\frac{k+1}{4}}X^2_0(s)\sin(2\pi s)\,
  ds=\sum_{k=0}^3\int_0^{\frac{1}{4}}X^2_0(\frac{k}{4}+s)\sin\bigl(2\pi
  (\frac{k}{4}+s)\bigr)\,
  ds\\&=-\sum_{k=0}^3\int_0^{\frac{1}{4}}X^2_0(\frac{k}{4}-s)\sin\bigl(2\pi 
  (\frac{k}{4}-s)\bigr)\, ds\\&=\sum_{k=0}^3\int^0_{\frac{1}{4}}X^2_0(\frac{k}{4}+s)\sin\bigl(2\pi
  (\frac{k}{4}+s)\bigr)\, ds\\&=-\int_0^1 X^2_0(s)\sin(2\pi s)\, ds,
\end{align*}
which yields
\begin{equation*}
\int_0^1X^2_0(s)\sin(2\pi s)\, ds=0.
\end{equation*}
Next, similar calculations show that
\begin{equation*}
  \int_0^{\frac{1}{2}} X^2_0(s)\cos(2\pi s)\, ds=0\text{ and }
  \int^1_{\frac{1}{2}} X^2_0(s)\cos(2\pi s)\, ds=0
\end{equation*}
so that, indeed,
$$
\int_0^1 X^2_0(s)\cos(2\pi s)\, ds=0.
$$
Next we see that
\begin{align*}
\int_0^1 X^2_0(s)\cos(4\pi s)\, ds&=\sum_{k=0}^7\int_0^{\frac{1}{8}}X^2_0(
\frac{k}{8}+s) \cos\bigl(4\pi (\frac{k}{8}+s)\bigr)\, ds\\&=\sum_{k=0}^7(-1)^{2k+1}\int_0^{\frac{1}{8}}X^2_0(
\frac{k}{8}-s) \cos\bigl(4\pi (\frac{k}{8}-s)\bigr)\, ds\\&=\sum_{k=0}^7\int^0_{\frac{1}{8}}X^2_0(
\frac{k}{8}+s) \cos\bigl(4\pi (\frac{k}{8}+s)\bigr)\, ds=-\int_0^1 X^2_0(s)\cos(4\pi s)\, ds.
\end{align*}
As for the last integral, notice that
$$
X^2_0(s)\sin(4\pi s)>0\text{ on }[0,1)\setminus\{ 0, \frac{1}{4}, \frac{1}{2}, \frac{3}{4}\} 
$$
and, hence, that
$$
\int _0^1 X^2_0(s)\sin(4\pi s)\, ds>0,
$$
as announced. As the linear heat equation and any reparametrization
maintains the symmetries \eqref{x-sym} and \eqref{y-sym}, it can be
inferred that the same integral (spectral) properties are enjoyed by
the approximants $X^n(t)$ and for all times $t>0$, regardless of the
value of $n\in \mathbb{N}$. We now estimate the components. Recall
that
$$
Y(t)=e^{-\frac{t-kh}{L^2(X^n_k)}\partial_{ss}}R(X^n_k),\: t\in
\bigl[kh,(k+1)h\bigr).
$$
It follows that
\begin{align*}
\| Y^j(t)\|_\infty&\leq \| Y^j(t)\|_{1,2}\leq
e^{-\frac{4\pi^2(t-h)}{L^2(X^n_k)}}\| R(X^n_k)\|_{1,2}\leq
e^{-\frac{4\pi^2(t-h)}{L^2(X^n_k)}}\| R(X^n_k)\|_{1,\infty}\\
&\leq 2 e^{-\frac{4\pi^2(t-h)}{L^2(X^n_k)}}L(X^n_k) ,\: t\in
\bigl[kh,(k+1)h\bigr),\: j=1,2.
\end{align*}
Here we used the facts that
$$
\| R(X^n_k)_s\|_\infty=L(X^n_k)\text{ and that }\| R(X^n_k)\|
_\infty\leq \frac{L(X^n_k)}{2},
$$
where the latter follows from the symmetries of the evolving curves
about the $x$- and $y$-axis as well as the fact that the diameter of a
closed curve is at most half its length. Now the decay enjoyed by
$Y^1(t)$ is precisely that suggested by the inequality since the
$$
\langle Y^1, \cos(2\pi\cdot)\rangle\neq 0,
$$
while that of $Y^2$ has a factor of $16\pi^2$ instead of $4\pi^2$
since
$$
\langle Y^2, \varphi(2\pi\cdot)\rangle= 0,\: \varphi=\cos,\sin.
$$
This shows that the aspect ratio of the evolving curve tends to zero
as the curve shrinks to a point while maintaining its symmetries. The
rescaled limit does therefore have to be the doubly covered segment of
length 1 that enjoys $x$- and $y$-axis symmetry as desired.
\end{proof}
\begin{rem}
It was already conjectured in \cite{AL86} that the limiting (normalized)
curvature of a symmetric figure eight curve be that of a doubly
covered segment. In this case the (normalized) curvature does not
converge to its asymptotic limit in the $L^1$-topology. By \cite{AL86}
all (closed) curves for which it does must converge to a (smooth)
homotetic solution. We now know that a doubly covered segment
instantaneously shrinks to a point and can therefore also be thought
as a very fast self-similar shrinker.
\end{rem}
\begin{rem}
In \cite{Gray89}, Grayson studied the evolution of planar figure eight
shapes (in his definition, curves with exactly a double point and vanishing signed area)
and showed that the isoperimetric ratio $\frac{L^2}{A}$ tends
to infinity as the extinction time is approached if and only if the
two loops enclose regions of equal area.
\end{rem}
\begin{rem}
For an initial curve to converge to a doubly covered segment, it is
enough to require lesser symmetry for the initial datum $Y_0$ as, for
instance, central symmetry
$$
Y_0(\frac{1}{2}+s)=-Y_0(\frac{1}{2}-s),\: s\in [0,1),
$$
or $y$-axis symmetry
$$
Y_0(\frac{1}{2}+s)=\begin{bmatrix}
  -1&0\\0&1\end{bmatrix}Y_0(\frac{1}{2}-s),
$$
which ensure that the total signed area vanishes all along the
evolution as all the approximants and hence the limit preserve these
symmetries. In this case, however, the location and direction of the limiting 
segment is not immediately available. See Figure
\ref{fig:infinityLoops} for such situations. The additional
symmetries enjoyed by $X_0$ make sure that the limiting segment is on
the $x$-axis in the chosen coordinates.
\end{rem}
\begin{rem}
The discussion of the long time behavior of the infinity loop and the
construction of the solution show that the curve shortening flow
exhibits a long time behavior that could be described as qualitatively
of heat type. By this we mean that a generic intial datum that has a persistent
non-zero projection onto the first non-constant modes will shrink to a
circular point as the first mode will dominate the decay. Clearly the
effect of diffusion is amplified and sped up by the length dependent
diffusivity which also causes finite exinction time. For a curve to
shrink to a non-circular point, it will need to have qualitative
properties (symmetries) which prevent the first non-constant modes to
dominate the evolution. The infinity loop considered above is such an
example. We will show in the section dedicated to numerical
experiments that small perturbations of these symmetries lead back to
convergence to a round point.
\end{rem}

\section{Symmetries and Long Time Behavior}
\subsection{Symmetries}
\begin{figure}
\begin{center}
  \includegraphics[scale=0.5]{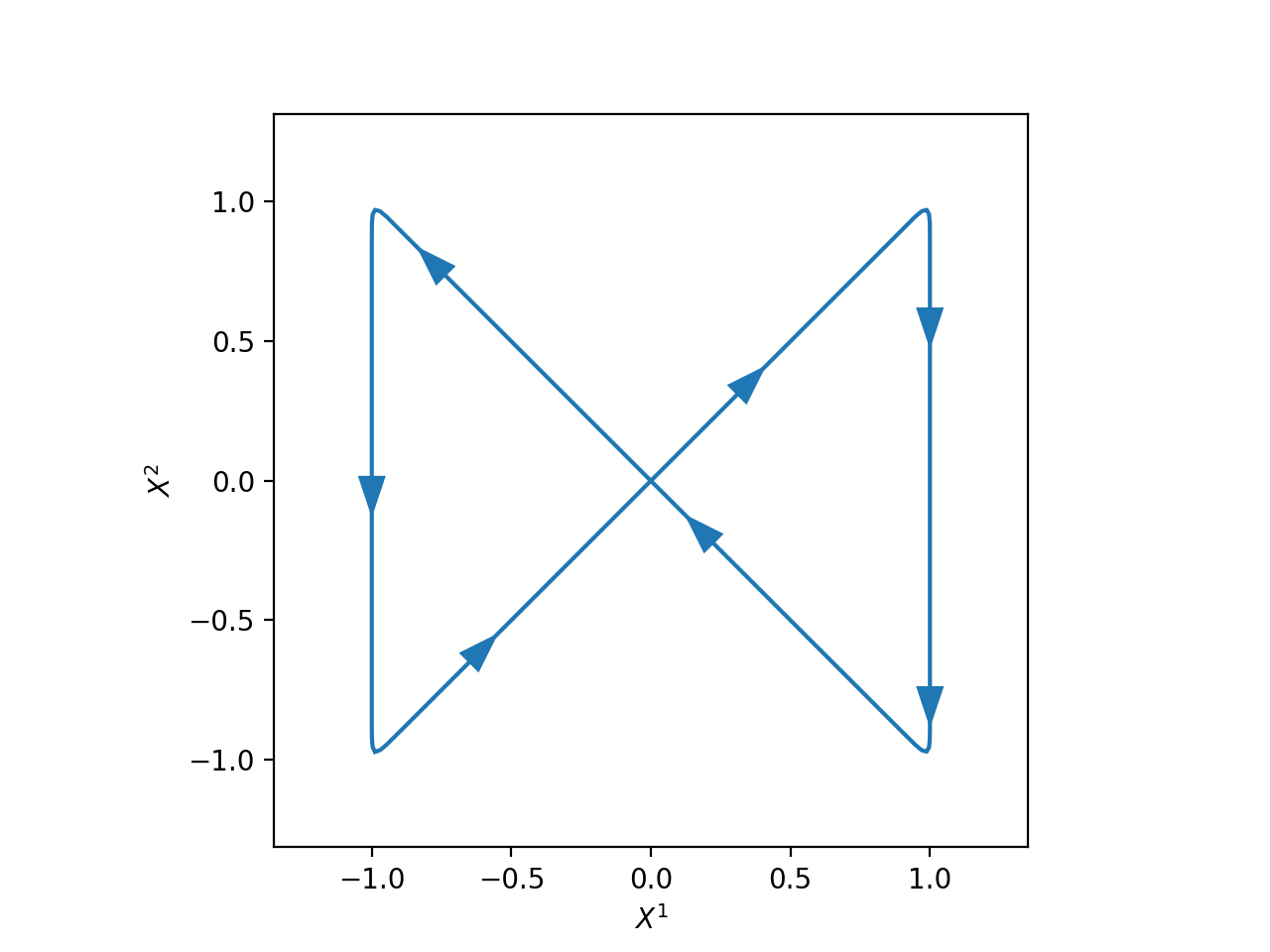}
  \includegraphics[scale=0.5]{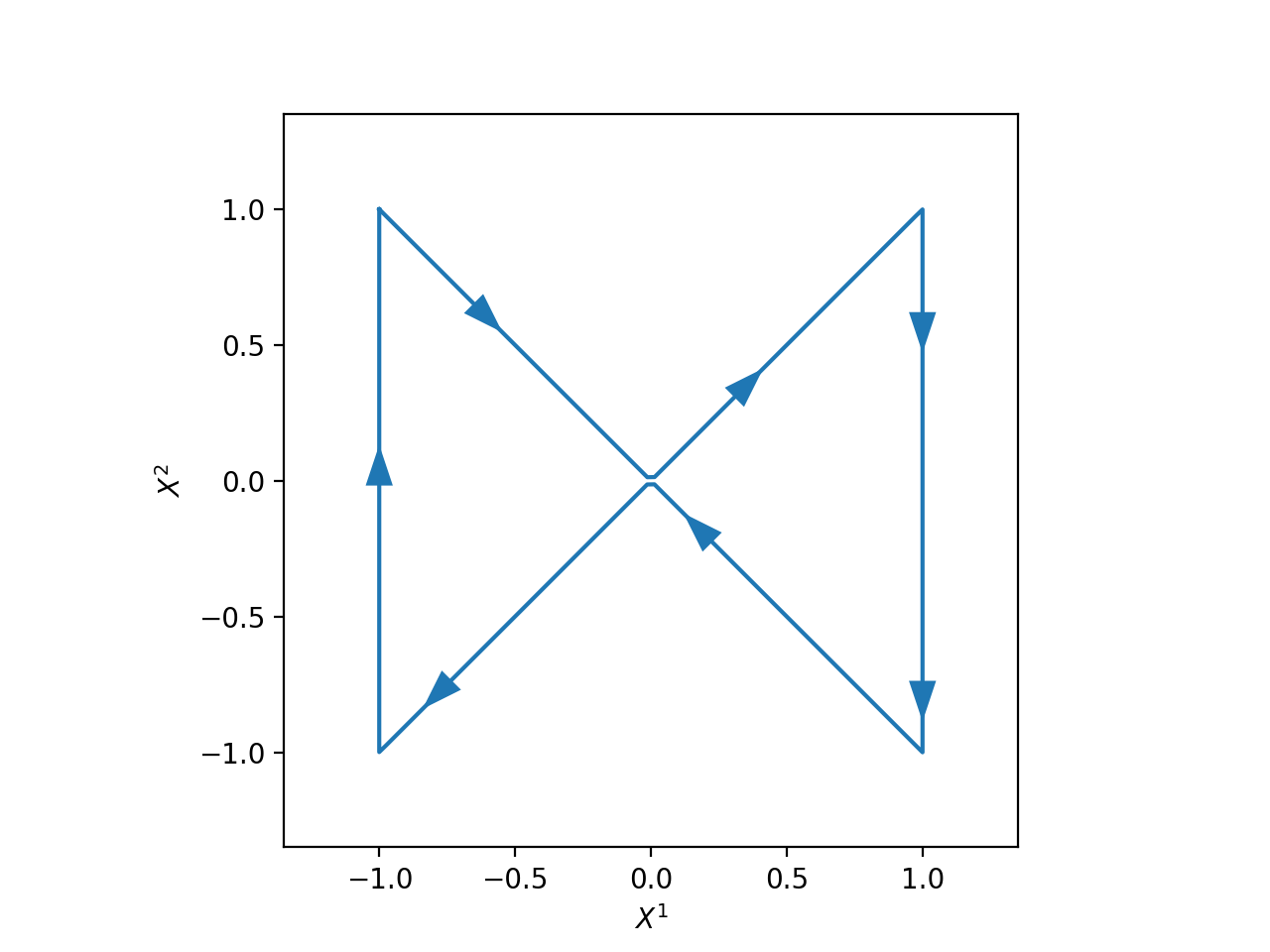}
\caption{Two distinct immersed curves sharing the same trace set. They
do carry different topologies that are both distinct from the topology
induced by the ambient plane.}
\label{fig:infinityLoopsCurves}
\end{center}
\end{figure}

\begin{figure}
\begin{center}
  \includegraphics[scale=0.5]{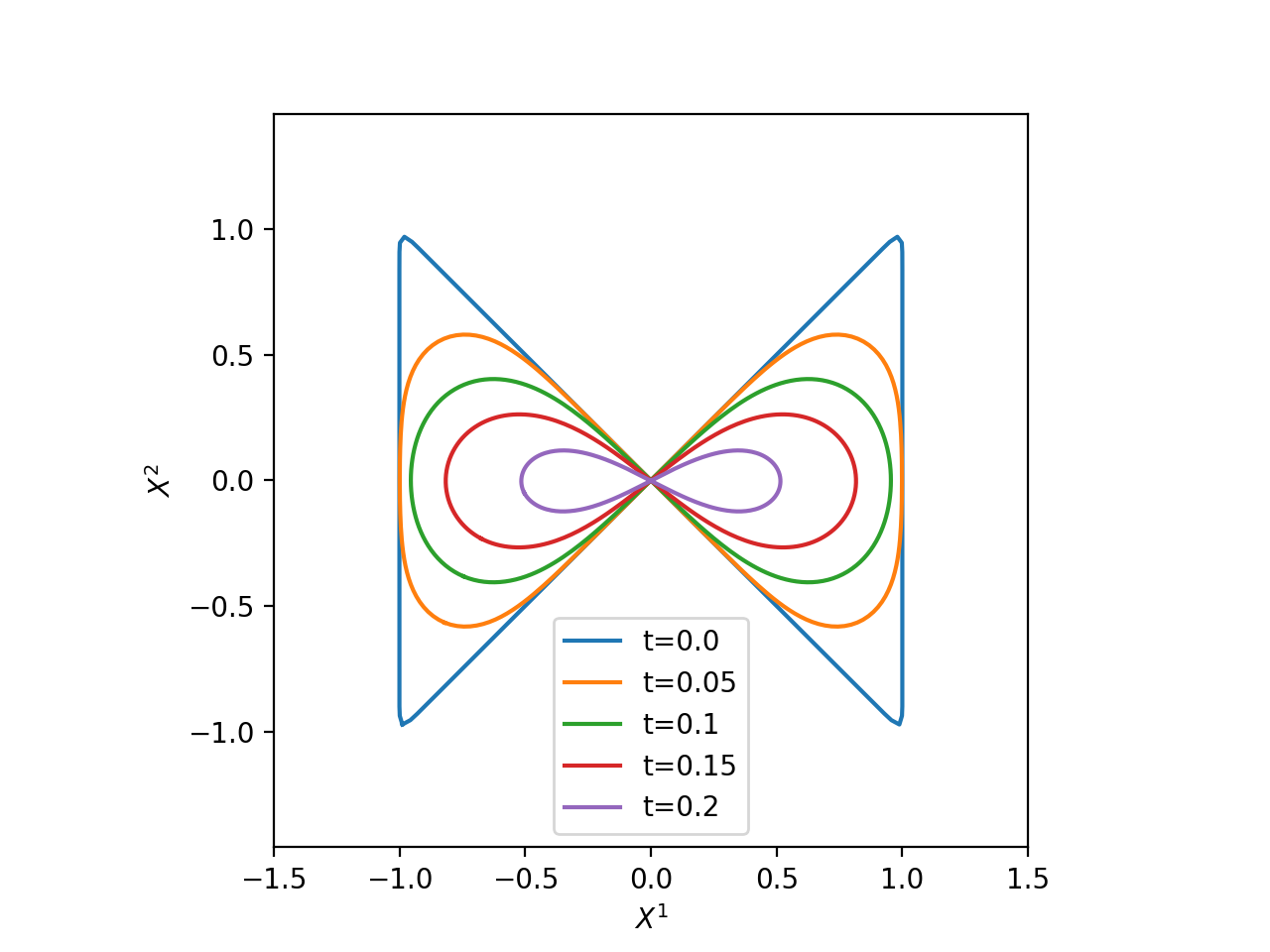}
  \includegraphics[scale=0.5]{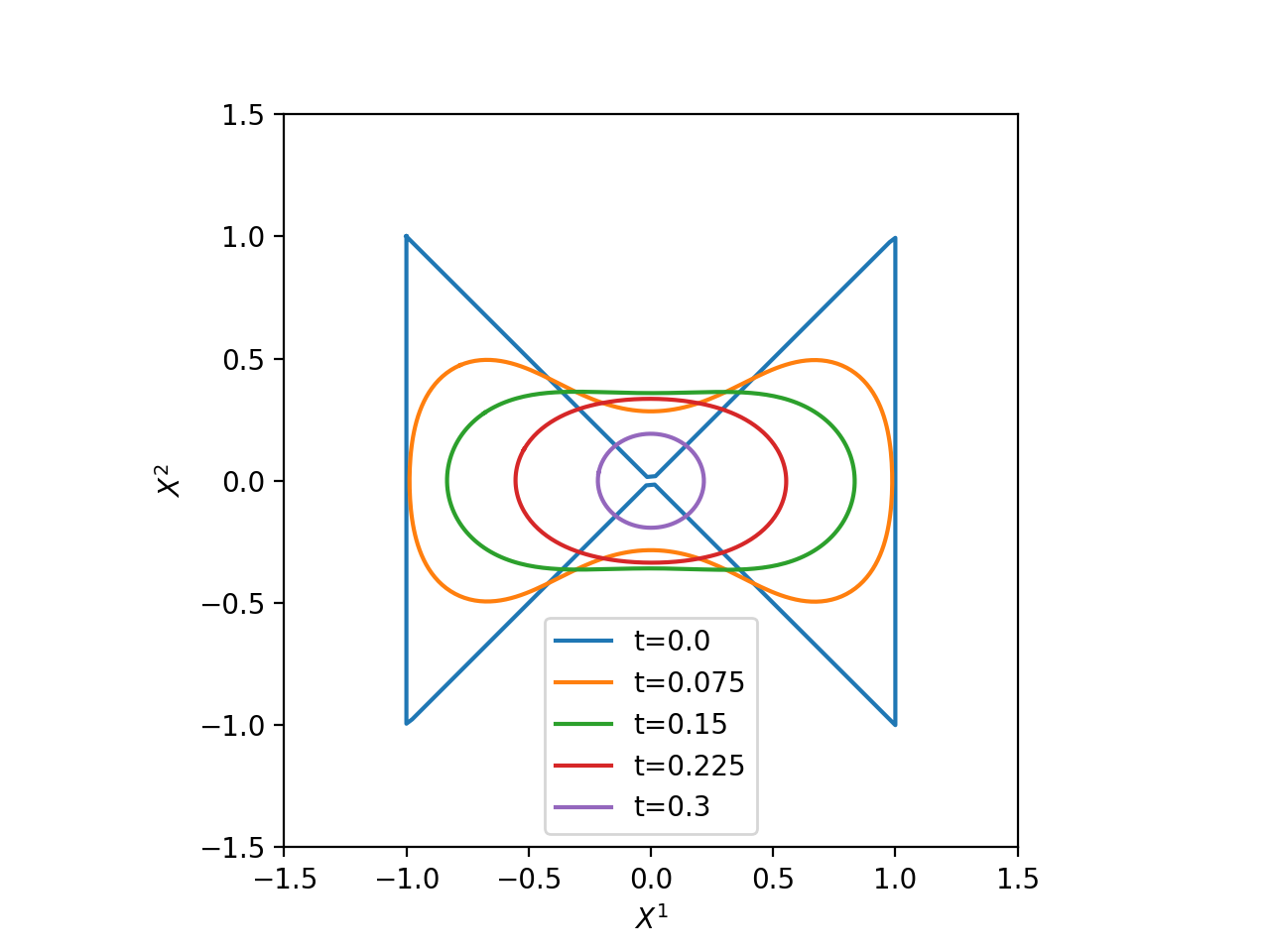}
\caption{Snapshots of the evolution of two topologically distinct
  curves sharing the same trace set.}
\label{fig:loopsEvo}
\end{center}
\end{figure}

Any parametrization $X$ of an immersed curve will have velocity $X_r$
and acceleration $X_{rr}$. These depend on the way the curve is
traversed and do not immediately provide intrinsic geometric information about the
immersed curve. The velocity can, however, be split into the product
of speed and direction
$$
 X_r=|X_r| \frac{X_r}{|X_r|},
$$
leading to a decomposition of the acceleration
$$
 X_{rr}= \frac{d}{dr}|X_r| \frac{X_r}{|X_r|}+|X_r|
 \bigl(\frac{X_r}{|X_r|}\bigr)_r= |X_r|_r \frac{X_r}{|X_r|}+ |X_r|^2 k(X)
$$
into a tangential (parametrization dependent) and a normal
component. The curvature $k(X)$ is the acceleration experienced by
traversing the curve at a constant unit speed. As such, it is an intrinsic property
of an immersed curve and not a manifestation of the specific way the
curve is run through. It is, however, not an intrinsic property of the
trace set $X(I)\subset \mathbb{R}^d$, which by itself cannot be
considered an immersed curve.
Parametrized curves are not simple subsets of $\mathbb{R}^d$. They are
topological spaces with respect to the topology induced by their
parametrization equivalence class. It follows that distinct curves $X$
and $Y$ with the same trace set, i.e. for which $X(I)=Y(I)$, need not
be the same topological space and, hence, manifold. While embedded
curves also carry the topology induced by their parametrization
classes, this topology does coincide with the ``natural''
topology of their trace set, which is the one induced by the ambient
space $\mathbb{R}^d$. Immersed curves, by
contrast, have topologies that do not coincide with the induced
topologies of their trace sets. To see this, take an infinity like
shape: any open neighborhood of the crossing point in the topology
induced by $\mathbb{R}^2$, will contain two crossing curve
segments. That is not the case for any of the two possible curve
topologies corresponding to the two distinct ways in which the trace
set can be run across as depicted in Figure
\ref{fig:infinityLoopsCurves}. These topologies are what distinguishes
their CSF evolutions. See Figure \ref{fig:loopsEvo}. If the CSF is
interpreted in a sense that only uses the initial trace set as a
topological subspace of the ambient space, then uniqueness may be lost
for any concept of weak solution.

This shows that symmetries enjoyed by the trace set of an initial
curve in the ambient space do not suffice to determine the evolution
of the curve by the CSF. One indeed needs
to consider the symmetries of its normalized parametrization since curves
with the same trace set can enjoy different symmetries as immersed
manifolds. Again we are confronted with the fact that embedding
symmetries coincide with ambient space symmetries, while immersion
symmetries cannot be inferred from ambient symmetries. In Figure
\ref{fig:infinityLoopsCurves}, the curve $X_l$ on the left satisfies
$$
X_l(\frac{1}{2}+s)=\underset{R_1}{\underbrace{\begin{bmatrix}
  -1&0\\0&-1\end{bmatrix}}}X_l(\frac{1}{2}-s),\: s\in[0,1),
$$
while the one on the right
$$
X_r(\frac{1}{2}+s)=\underset{R_c}{\underbrace{\begin{bmatrix}
  -1&0\\0&1\end{bmatrix}}}X_r(\frac{1}{2}-s),\: s\in[0,1).
$$
The common trace set $C=X_l(I)=X_r(I)$, however, satisfies
$$
R_1(C)=C,\: R_c(C)=C,\text{ and }-R_1(C)=C.
$$
We also observe that $\kappa_{tot}(X_l)\neq \kappa_{tot}(X_r)$.
\subsection{Long Time Behavior}
It was shown by Altschuler \cite{Alt91} that space curves, i.e. curves
in $\mathbb{R}^3$, that develop a singularity, are asymptotic to
planar Abresch-Langer homotetic solutions if the singularity is of
type I (and thus not a cusp). Type II singular behavior is also
studied and shown to lead to Gream Reapers by blow up around points
where the curvature explodes as the singular time is approached. This
shows that space curves of the CSF flow constructed in this paper
shrink to a point in a self-similar fashion (asymptotic to an
Abresch-Langer curve, which include multiply covered circles) or to a
$2m$-covered segment for $m\geq 1$. The latter can be considered as a 
degenerate Abresch-Langer curve with ``rotation index'' $m$  and such that
the curvature closes up in $m$ periods, and, hence, in some sense, a
degenerate $m$-covered circle.
\section{A Numerical Scheme}
The inspiration for the numerical method to be described shortly is to
be found in the construction of the solution presented earlier and in
the smoothing properties of the curve shortening flow. We think of a
discrete curve as a $n$-tuple of points $\mathbb{X}=(X_0,X_1,\dots,
X_{n-1})$. As we only consider closed curves, we can visualize
$\mathbb{X}$ as the polygon $P(\mathbb{X})$ obtained by connecting
$X_i$ to $X_{(i+1)\operatorname{mod}n}$ for $i=0,\dots,n-1$. The
length of the curve $\mathbb{X}$ is taken to be the length of the
polygon $P(\mathbb{X})$ so that
$$
L(\mathbb{X})=L \bigl( P(\mathbb{X})\bigr) = \sum_{i=0}^{n-1}\big |
X_{(i+1)\operatorname{mod}n}-X_i\big |.
$$
Given a discrete closed curve $\mathbb{X}$, its natural arclength
parametrization is given by
$$
X_k=X(l_k), k=0,\dots, n-1,
$$
where arclength is itself given by
$$
l_k=\sum_{i=0}^{k-1}\big | X_{(i+1)\operatorname{mod}n}-X_i\big  |,\: k=1,\dots,n,
$$
with the understanding that $l_0=0$. It can naturally be extended to a
continuous parametrization of $P(\mathbb{X})$ by setting
$$
X(t)=\frac{t-l_k}{l_{k+1}-l_k}X_k+\frac{l_{k+1}-t}{l_{k+1}-l_k}X_{(k+1)
  \operatorname{mod}n},
$$
for $t\in (l_k,l_{k+1})$ and $ k=0,\dots,n-1$. In order to simplify
the computation of curvature, needed for the evolution, we replace
the non-uniform discrete parametrization of $P(\mathbb{X})$ by a uniform
one comprised of $N\geq n$ points $\mathbb{Y}=(Y_0,\dots,Y_{N-1})$ along
$P(\mathbb{X})$ such that
$$
d_{P(\mathbb{X})}(Y_k,Y_{k+1})=\frac{L(\mathbb{X})}{N},\: k=0,\dots, N-1,
$$
where $d_{P(\mathbb{X})}$ is the distance measured along $P(\mathbb{X})$.
\begin{rem}
If $\mathbb{X}$ is a discretization of a smooth curve of some
order (of accuracy) $p>1$, then so is $\mathbb{Y}$. Thus the
reparametrization step does not affect accuracy.
\end{rem}
\begin{proof}
If $\mathbb{X}^h=\big\{ x^h_i\,\big |\, i=n,\dots, n(h)-1\big\}$ is an
approximating set sequence for a curve $\Gamma\subset \mathbb{R}^d$
such that
$$
d(x^h_i,\Gamma)\leq Ch^p,\: i=0,\dots,n(h)-1,\: h>0,
$$
then, for any point $x$ along the polygon $P(\mathbb{X}^h)$ it holds
that $x=(1-\tau)x^h_i+\tau x^h_{(i+1)\operatorname{mod}n(h)}$ for some
$i\in\{0,\dots,n(h)-1\}$ and $\tau\in(0,1)$. Thus we have that
$$
d(x,\Gamma)\leq (1-\tau)d(x^h_i,\Gamma)+\tau d(x^h_{(i+1)\operatorname{mod}n(h)},
\Gamma)\leq Ch^p,
$$
as claimed.
\end{proof}
We denote the constructed polygon $\mathbb{Y}\in \mathbb{R}^N$ by
$R_N(\mathbb{X})$ as it corresponds to the reparametrization operator
$R$ used earlier in the continuous setting. It is then natural to
define the curvature $k(\mathbb{X})$ of $\mathbb{X}$ by
$$
k(\mathbb{X})=\frac{1}{L^2(\mathbb{X})}\mathcal{F}_N^{-1}D_N^2
\mathcal{F}_NR_N(\mathbb{X}),
$$
where $\mathcal{F}_N$ denotes the discrete Fourier transform and $D_N$
denotes the diagonal matrix with diagonal entries given by $-2\pi i k$
for $k=-\frac{N}{2},\dots,\frac{N}{2}-1$. In this case $D_N^2$ is the
symbol of the spectral discrete Laplacian. With these notations
the numerical scheme reads
\begin{equation}\label{dcsf}
  \begin{cases}
    \mathbb{X}^{n+1}=\mathcal{F}_N^{-1}
    e^{-\frac{h}{L^2(\mathbb{X}^n)}D_N^2}\mathcal{F}_NR_N(\mathbb{X}^n),&n\geq
    0,\\
   \mathbb{X}^0,&\text{ given,}
 \end{cases}
\end{equation}
with a discrete time step $h>0$.
\section{Numerical Exeriments}
\subsection{Infinity Like Shapes}
\begin{figure}
\begin{center}
  \includegraphics[scale=0.4]{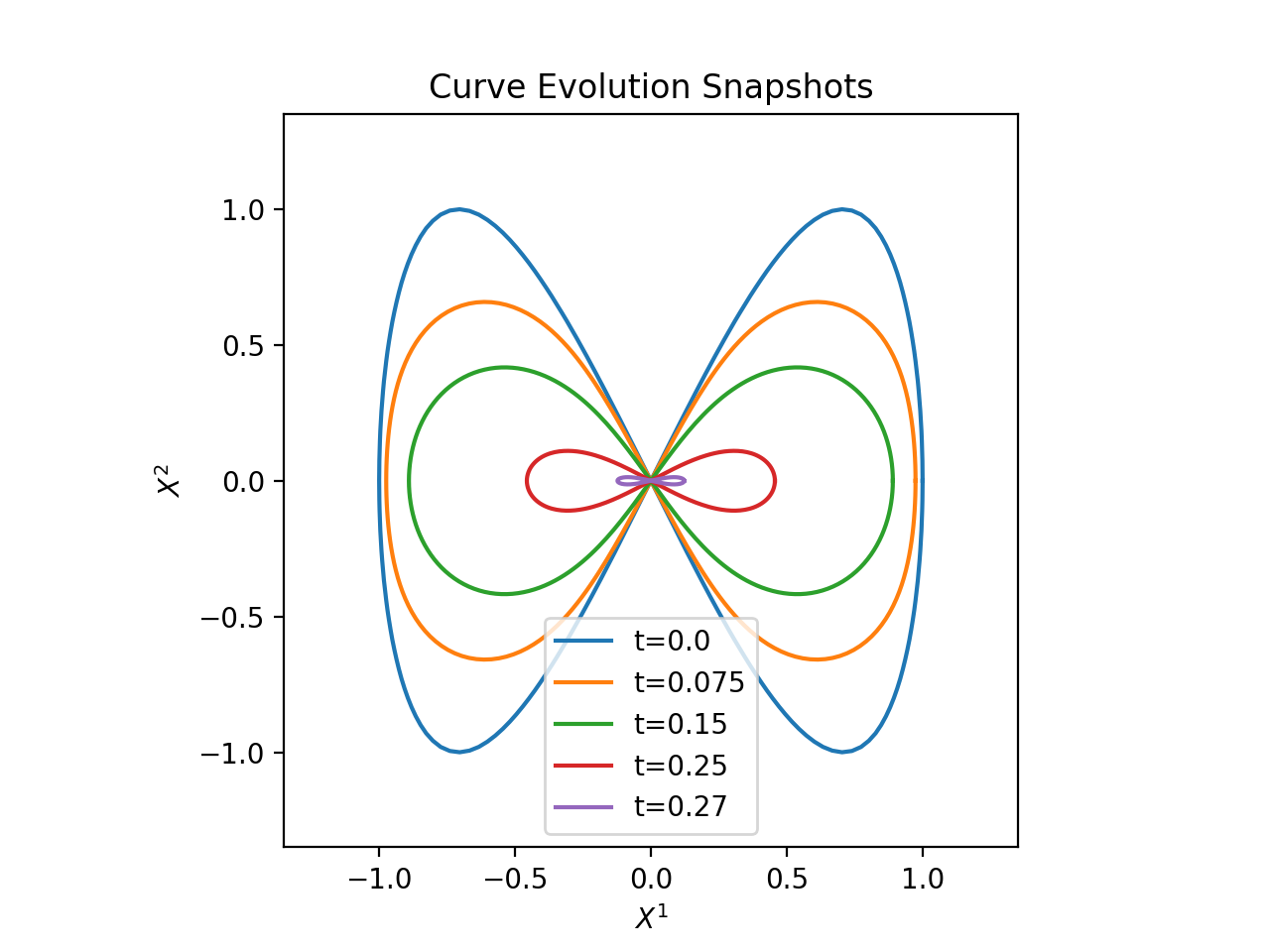}
  \includegraphics[scale=0.4]{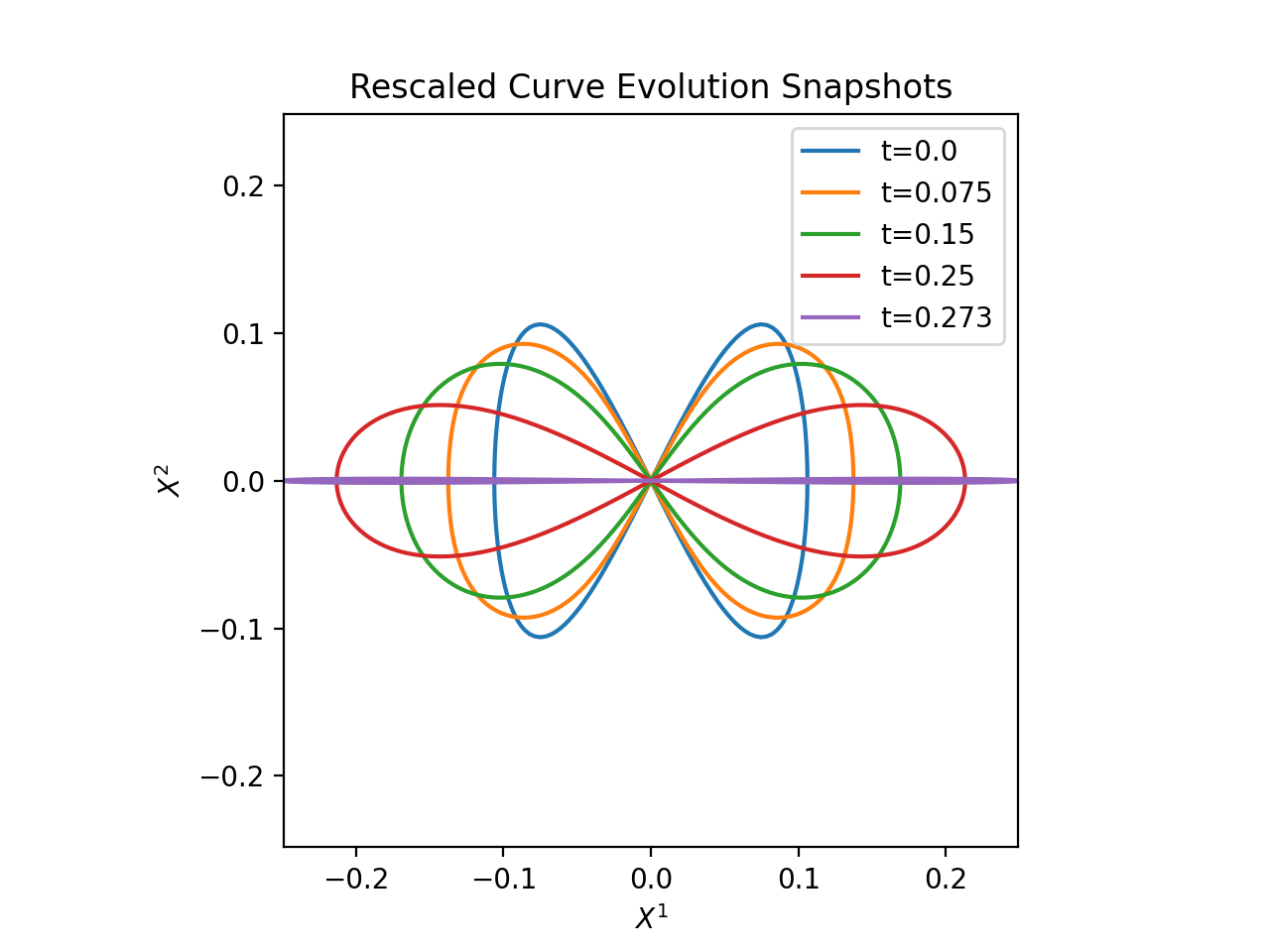}\\
  \includegraphics[scale=0.4]{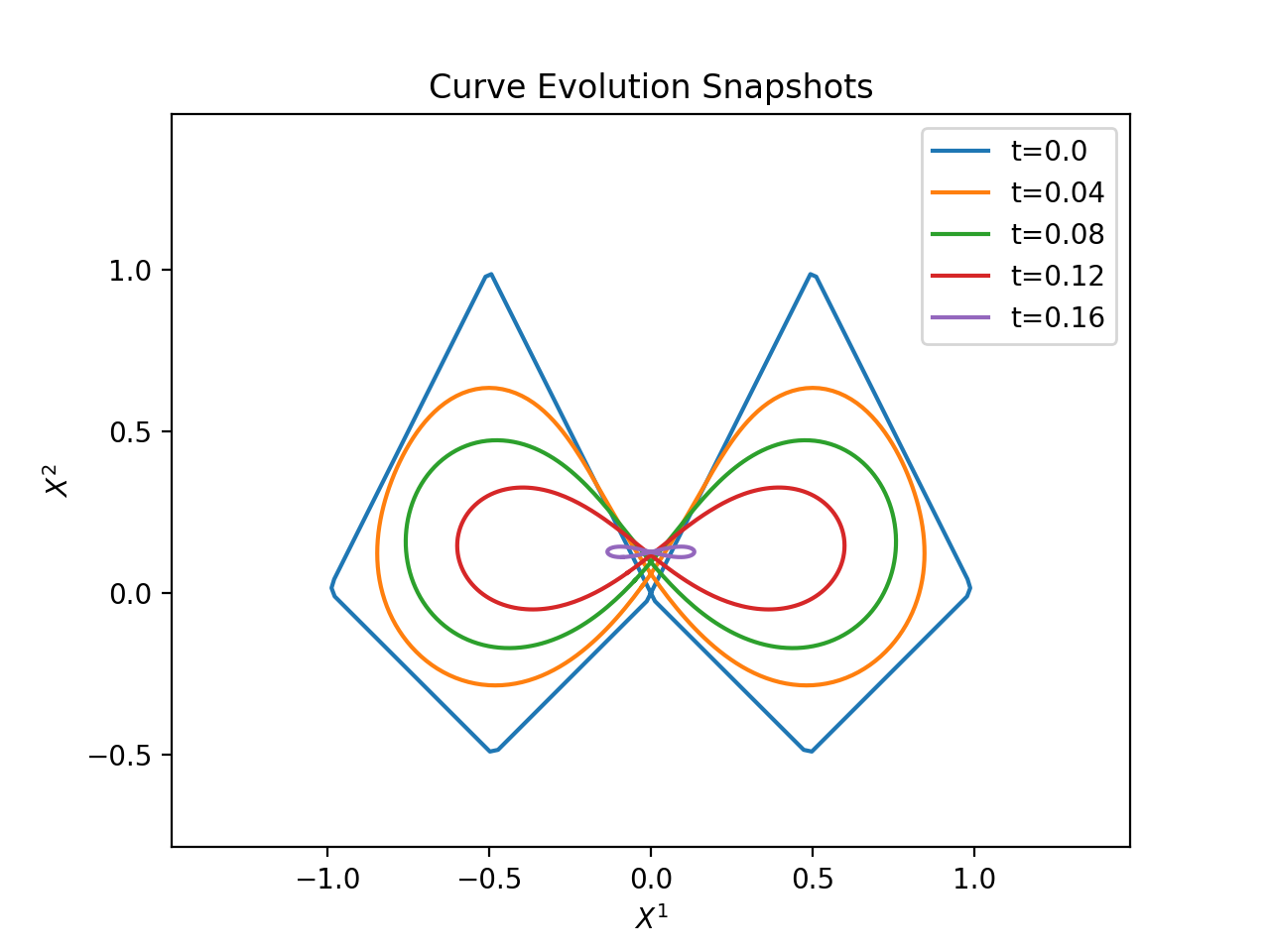}
  \includegraphics[scale=0.4]{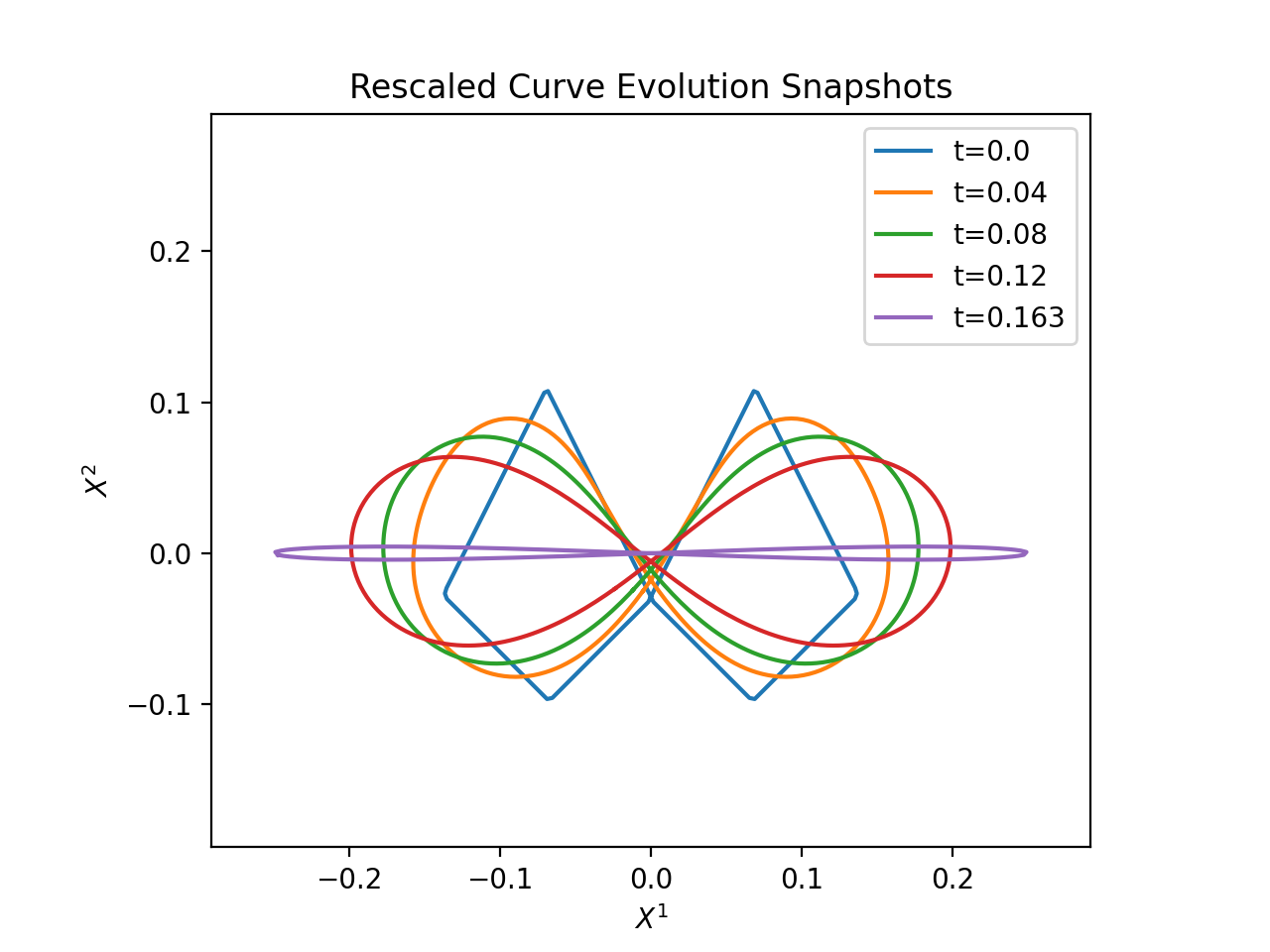}\\
  \includegraphics[scale=0.4]{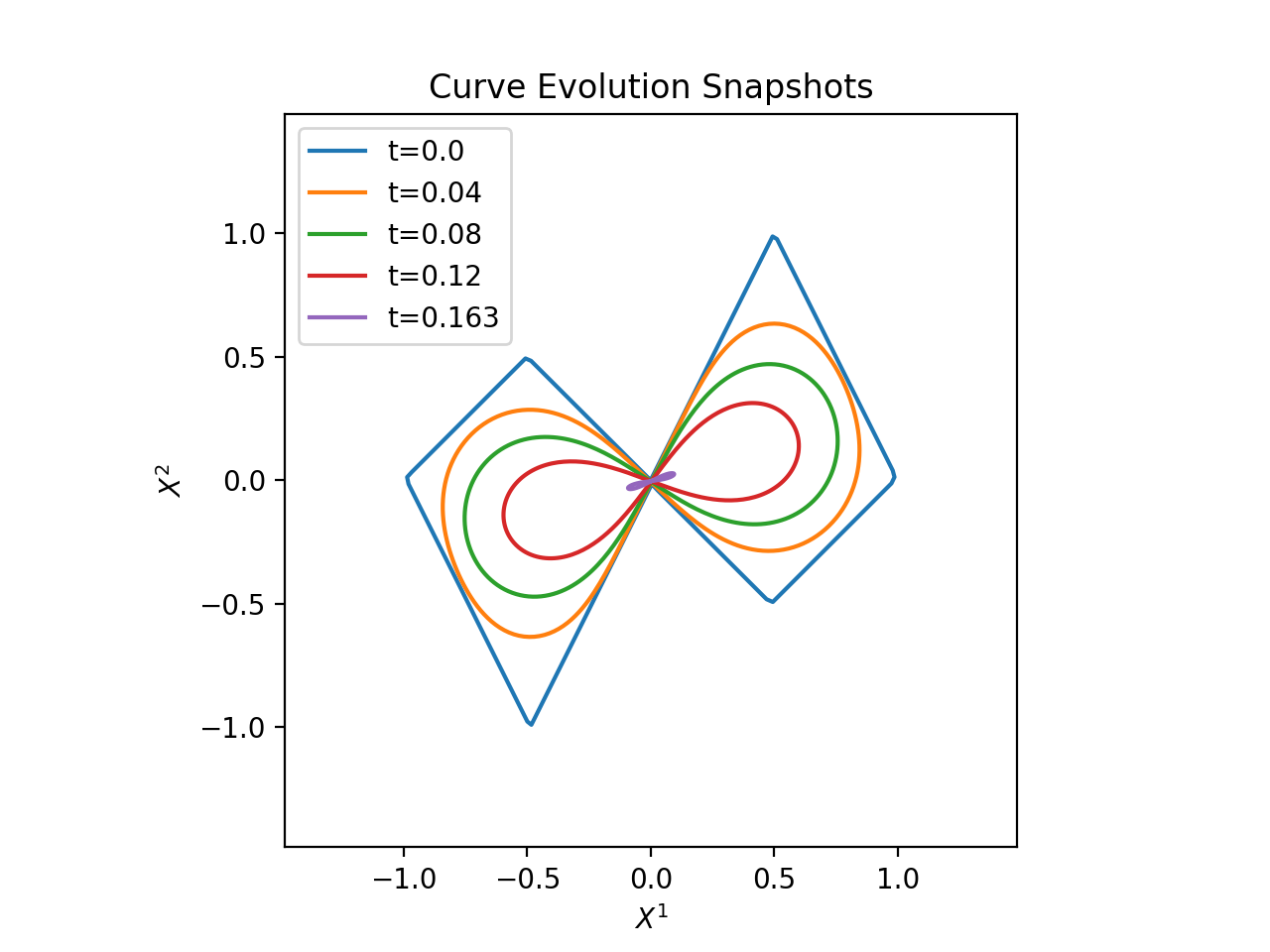}
  \includegraphics[scale=0.4]{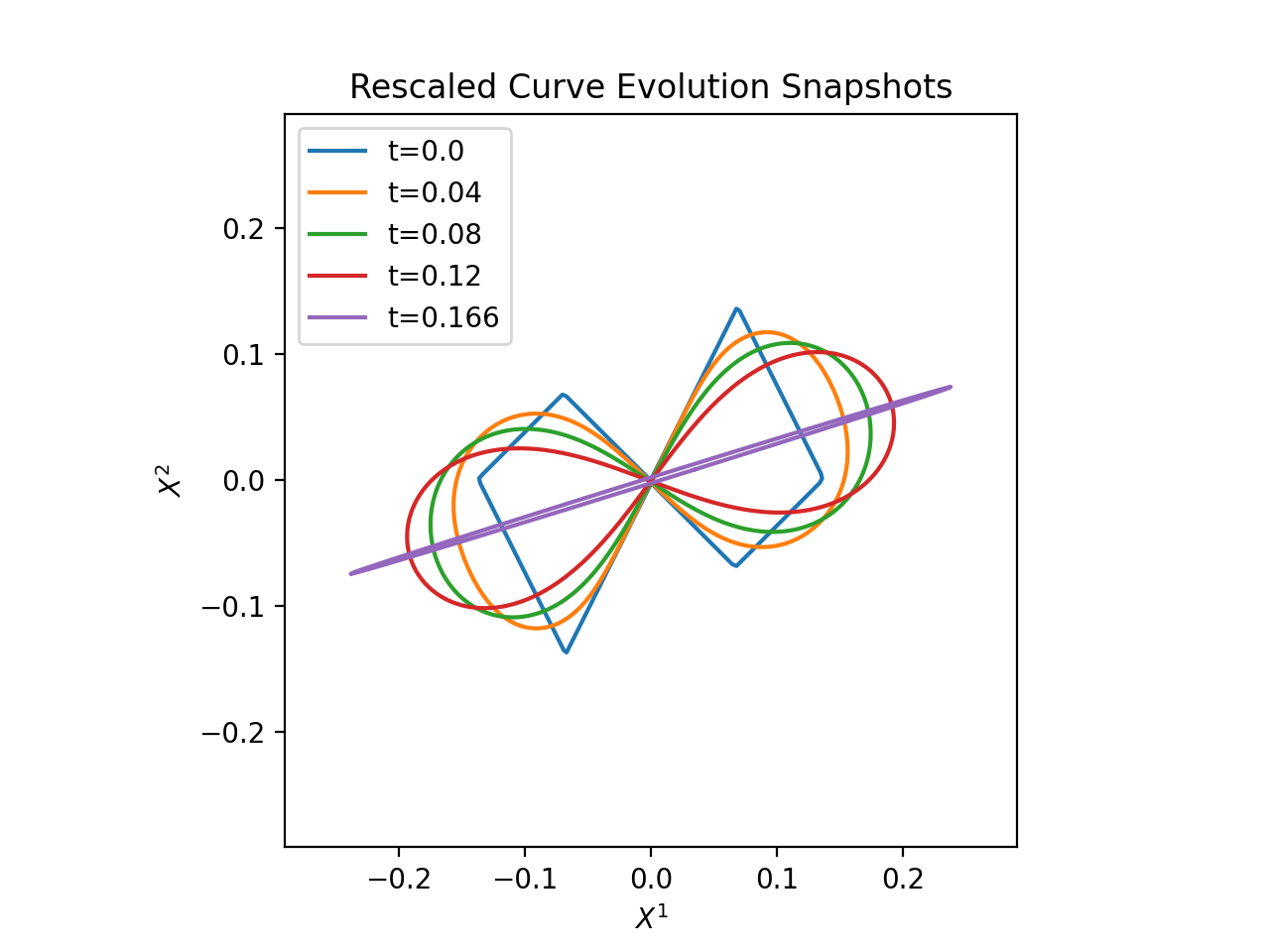}
\caption{On the left, the evolution of different initial infinity-like
  curve (in blue) is shown, while, on the right, the same curves are
  rescaled to have unit length to clearly show the asymptotic shape in
the singularity. In the second example the rescaled curves are
recentered around their center of mass since, in this case, the curve
evolution moves the center of the shape due to a lack of symmetry with
respect to the $x$-axis.}
\label{fig:infinityLoops}
\end{center}
\end{figure}

\begin{figure}
\begin{center}
  \includegraphics[scale=0.5]{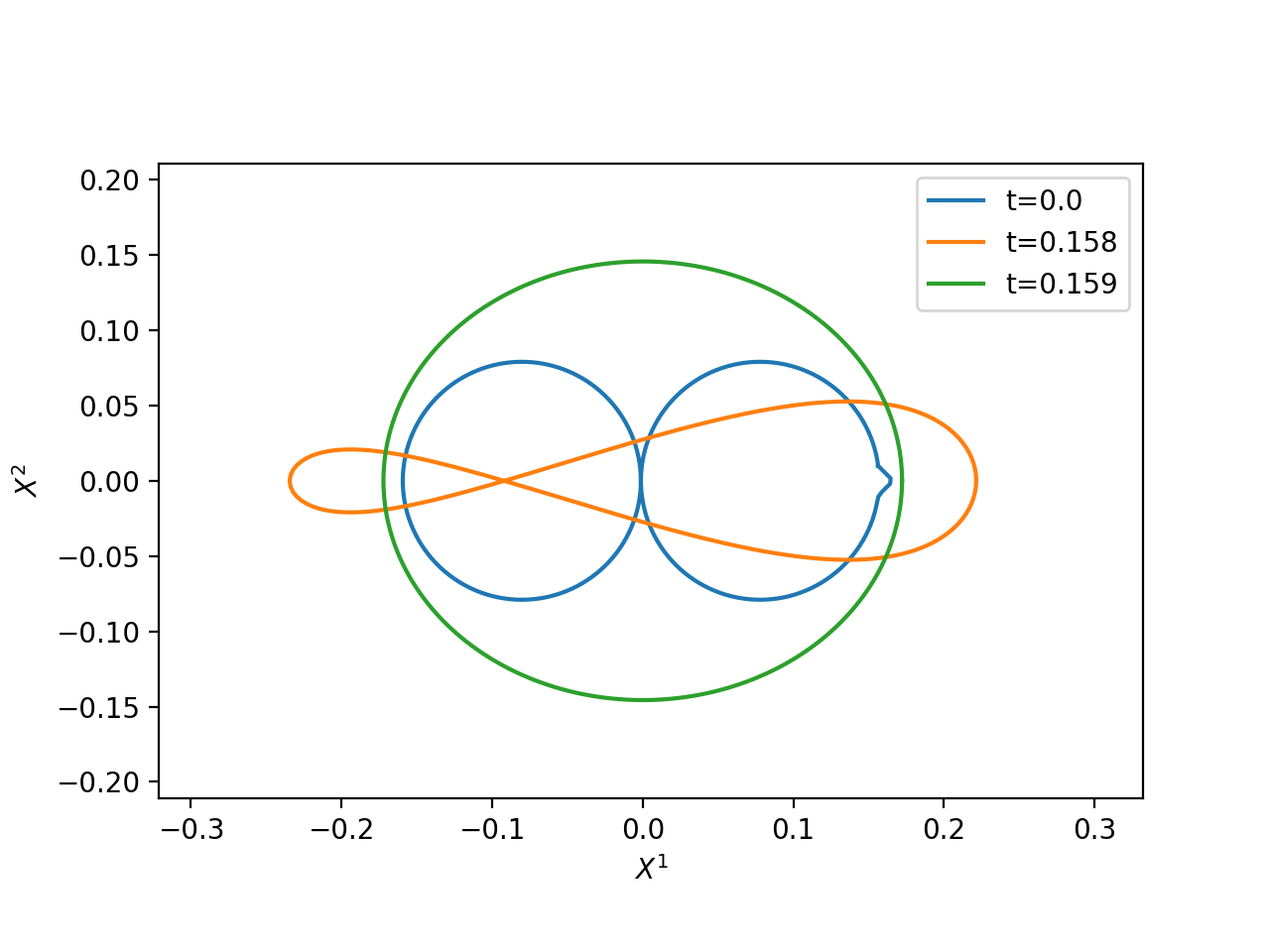}
\caption{Evolution of a slightly perturbed infinity shape. Curves are
  rescaled to have unit length.}
\label{fig:infPert}
\end{center}
\end{figure}

We consider three infinity-like shapes: one with both $x-$axis and
$y$-axis symmetries, one with only $y$-axis symmetry and one with
central symmetry. These evolutions, depicted in Figure
\ref{fig:infinityLoops}, exhibit the expected asymptotic 
behavior, i.e. convergence to a point not asymptotically circular, but
rather in the shape of a doubly covered segment. Shown are
snapshots of the actual curve evolution and a version rescaled to have
length one and recentered in the center of mass of the curve.
If we modify the infinity shape by adding a small bump to the loop
on the right. This break of symmetry is enough to change the long time
behavior of the solution as shown in Figure \ref{fig:infPert}. The curve
now shrinks to a round point.

\subsection{A Convoluted Curve}
\begin{figure}
\begin{center}
  \includegraphics[scale=0.33]{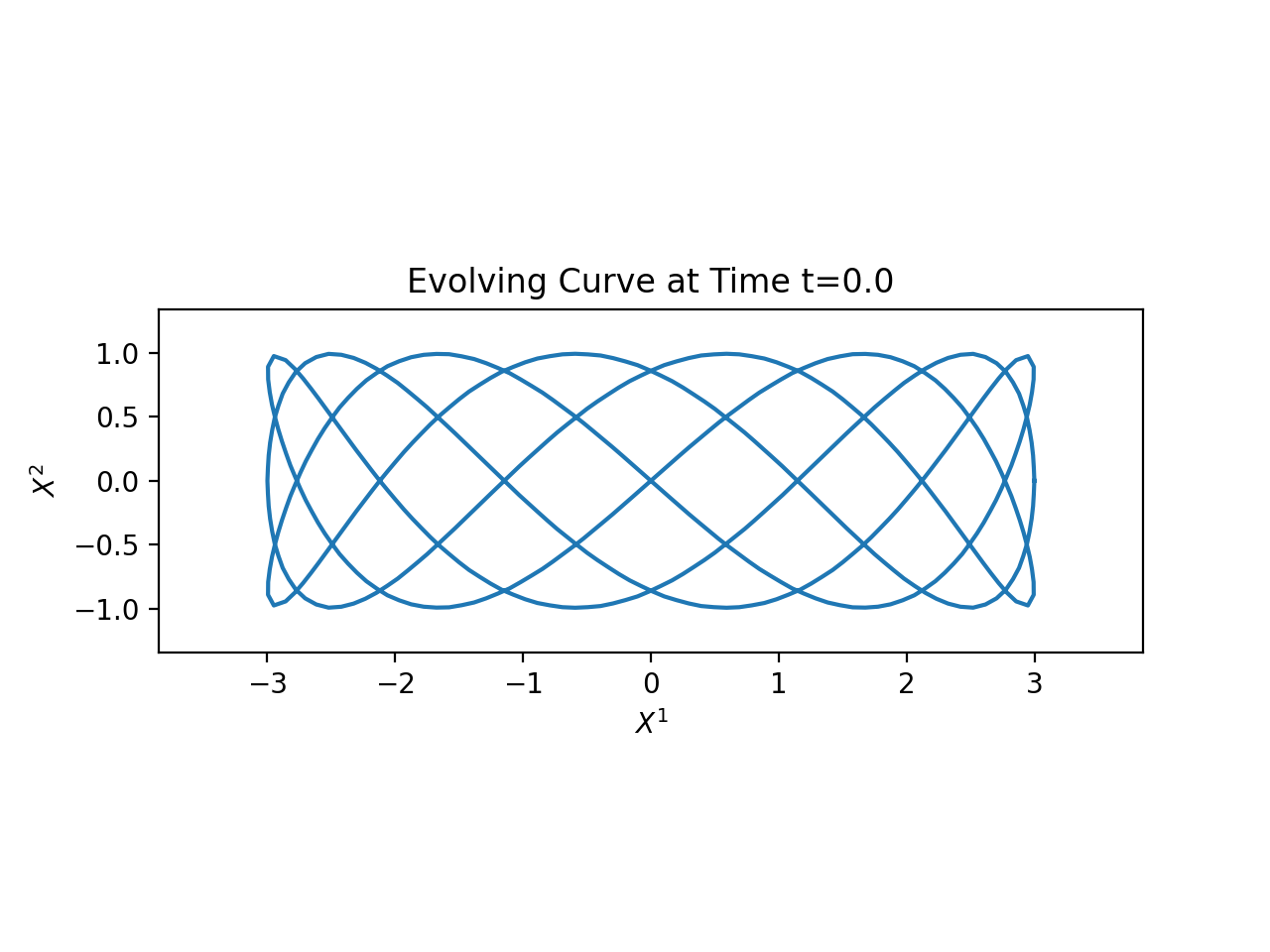}
  \includegraphics[scale=0.33]{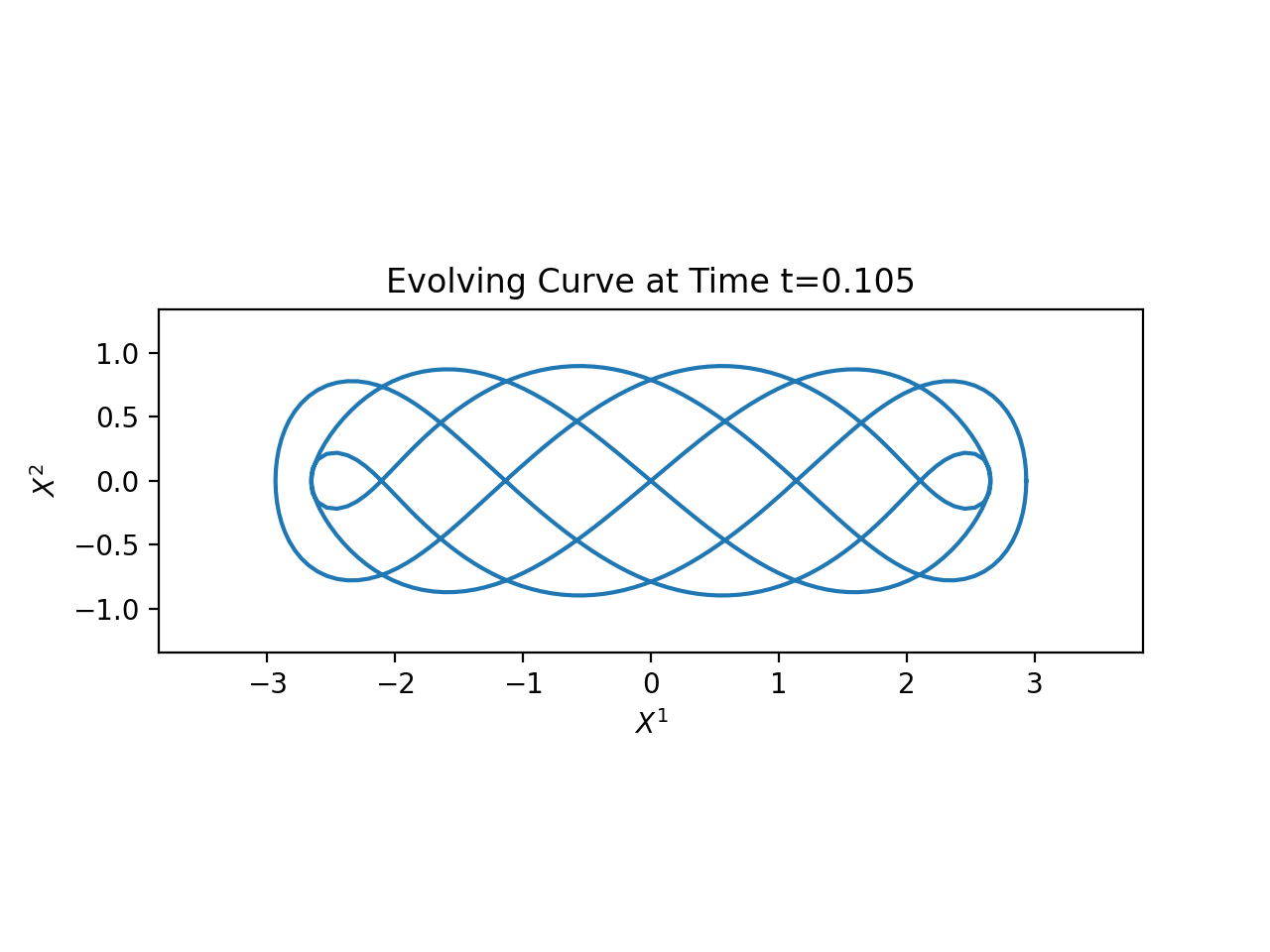}
  \includegraphics[scale=0.33]{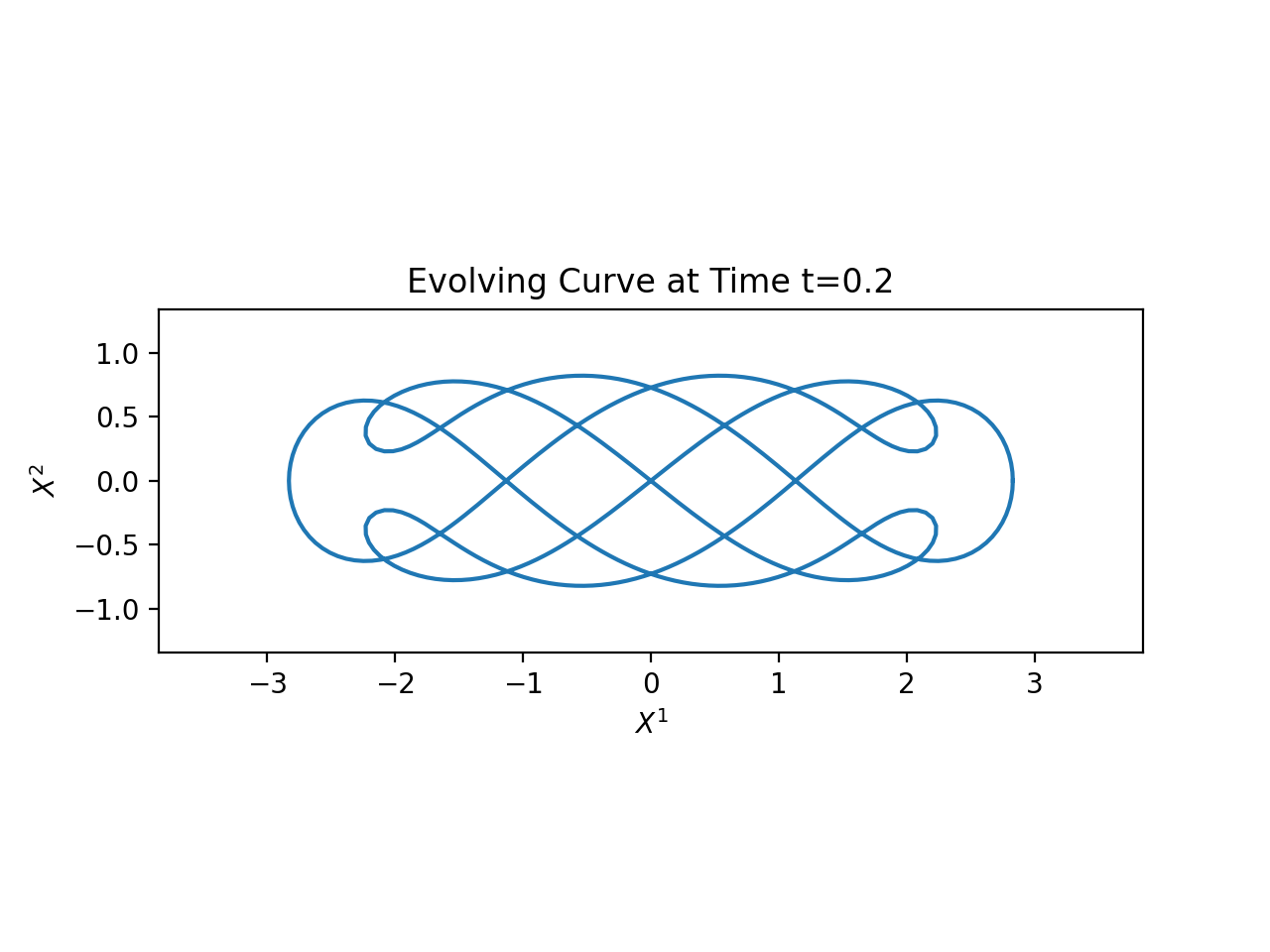}\\
  \includegraphics[scale=0.33]{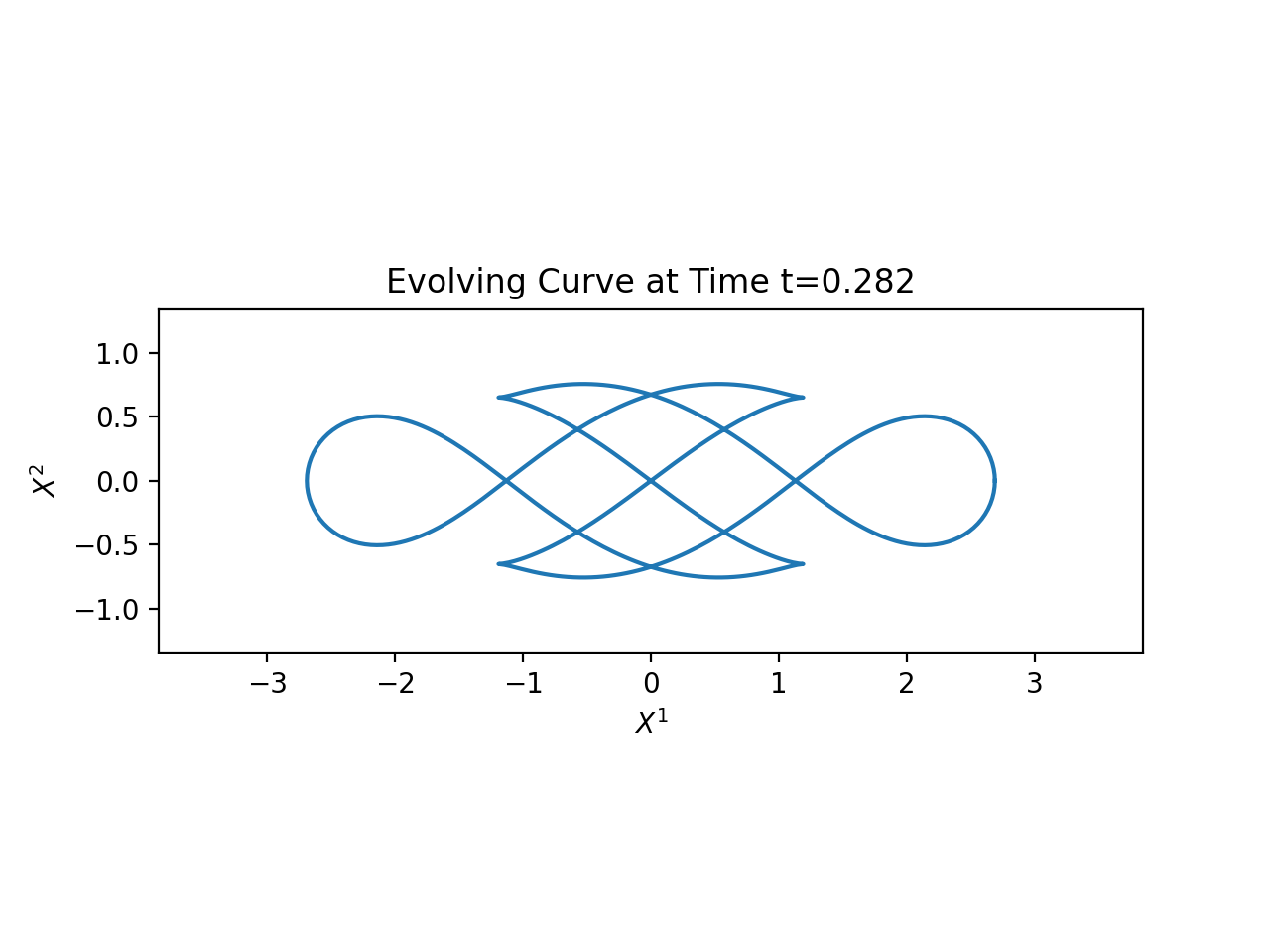}
  \includegraphics[scale=0.33]{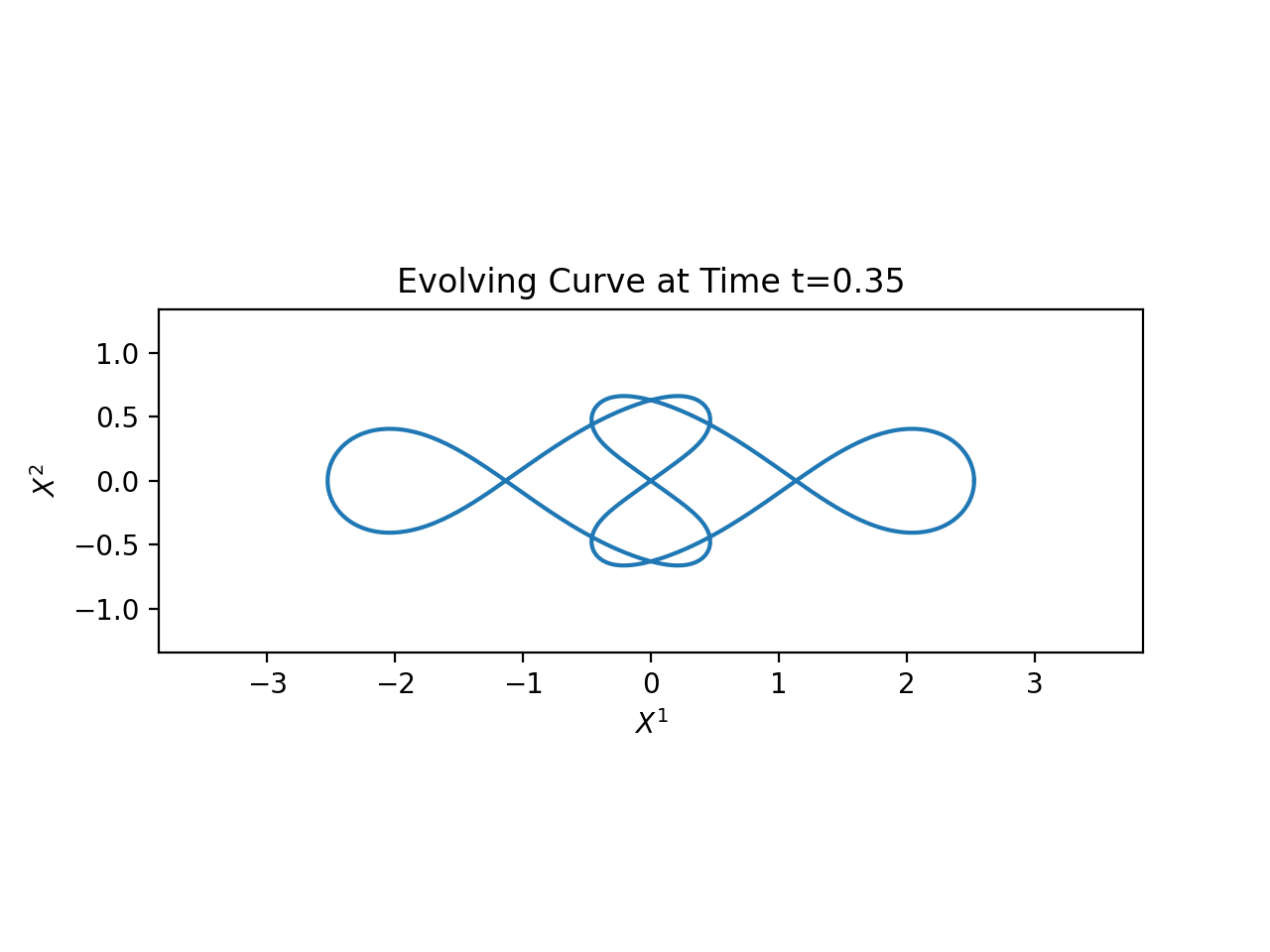}
  \includegraphics[scale=0.33]{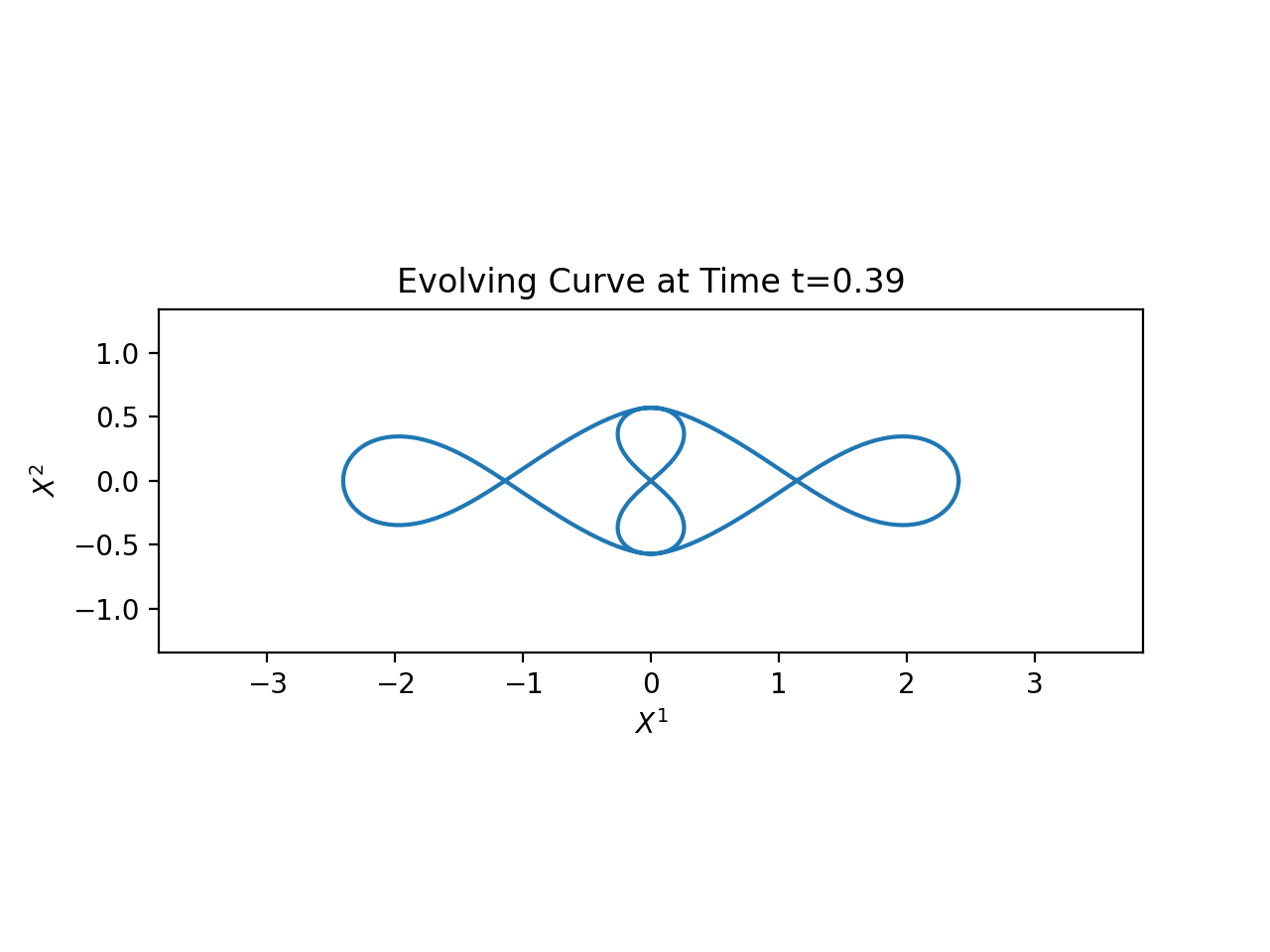}\\
  \includegraphics[scale=0.33]{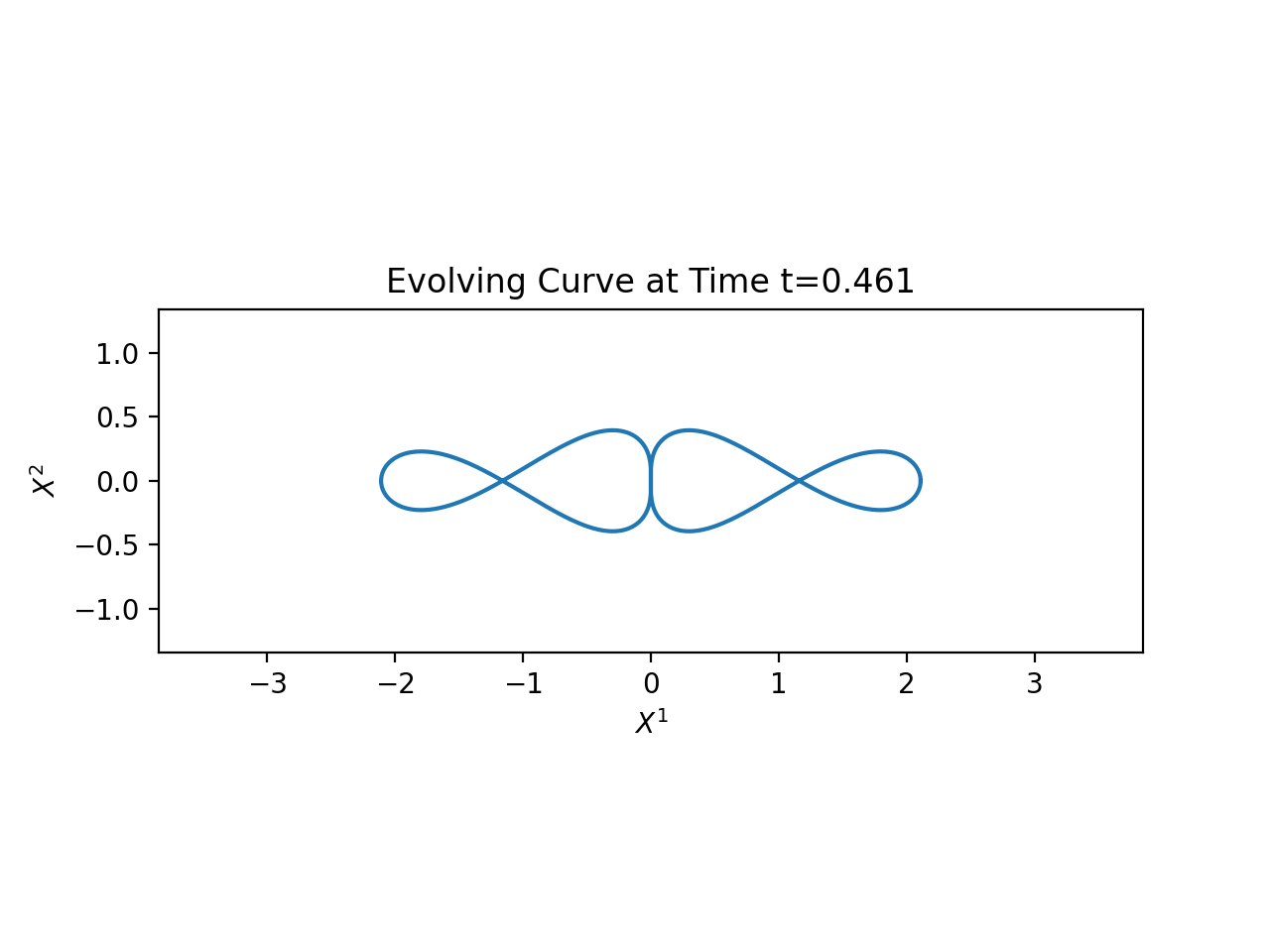}
  \includegraphics[scale=0.33]{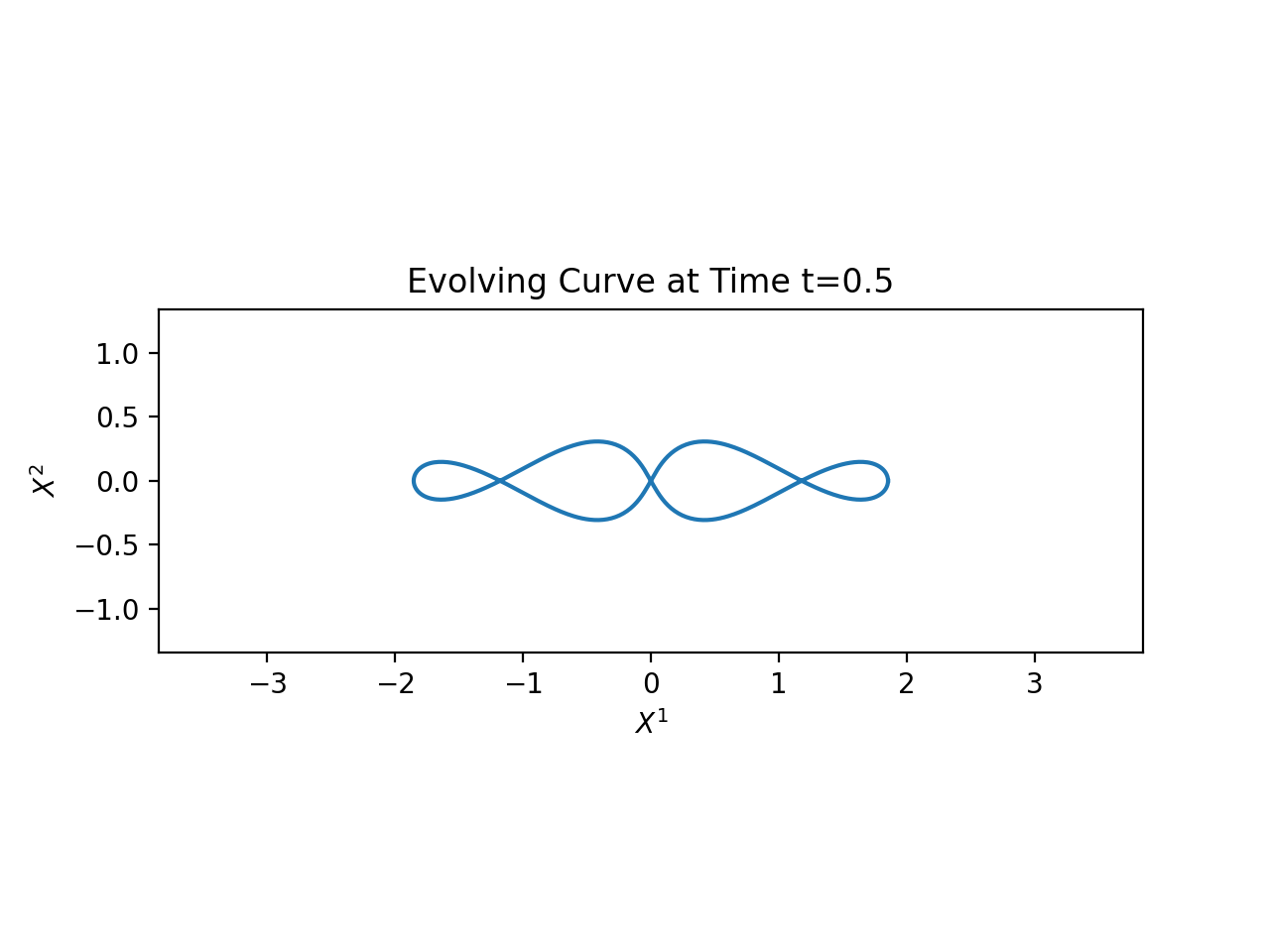}
  \includegraphics[scale=0.33]{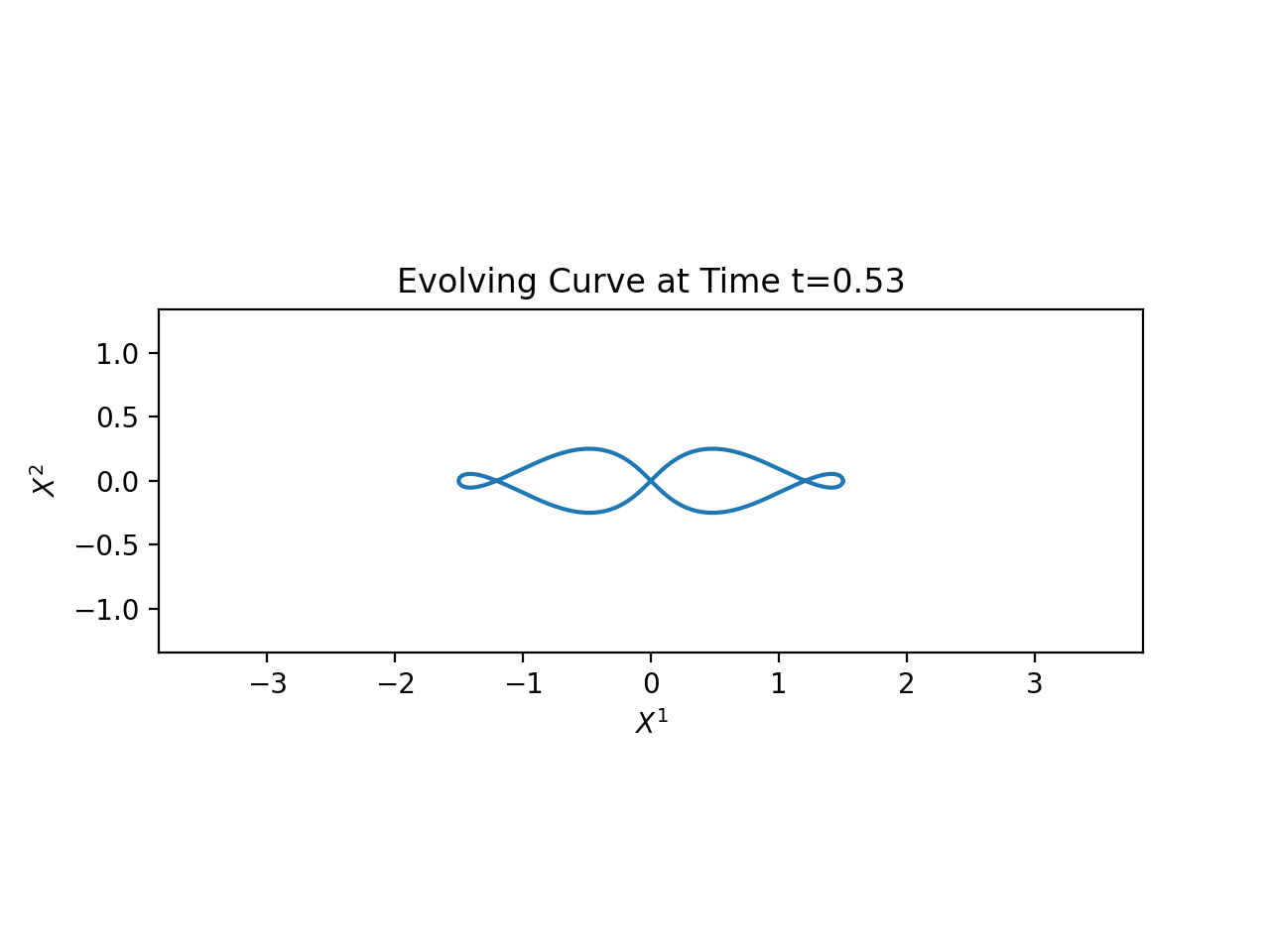}
  \includegraphics[scale=0.33]{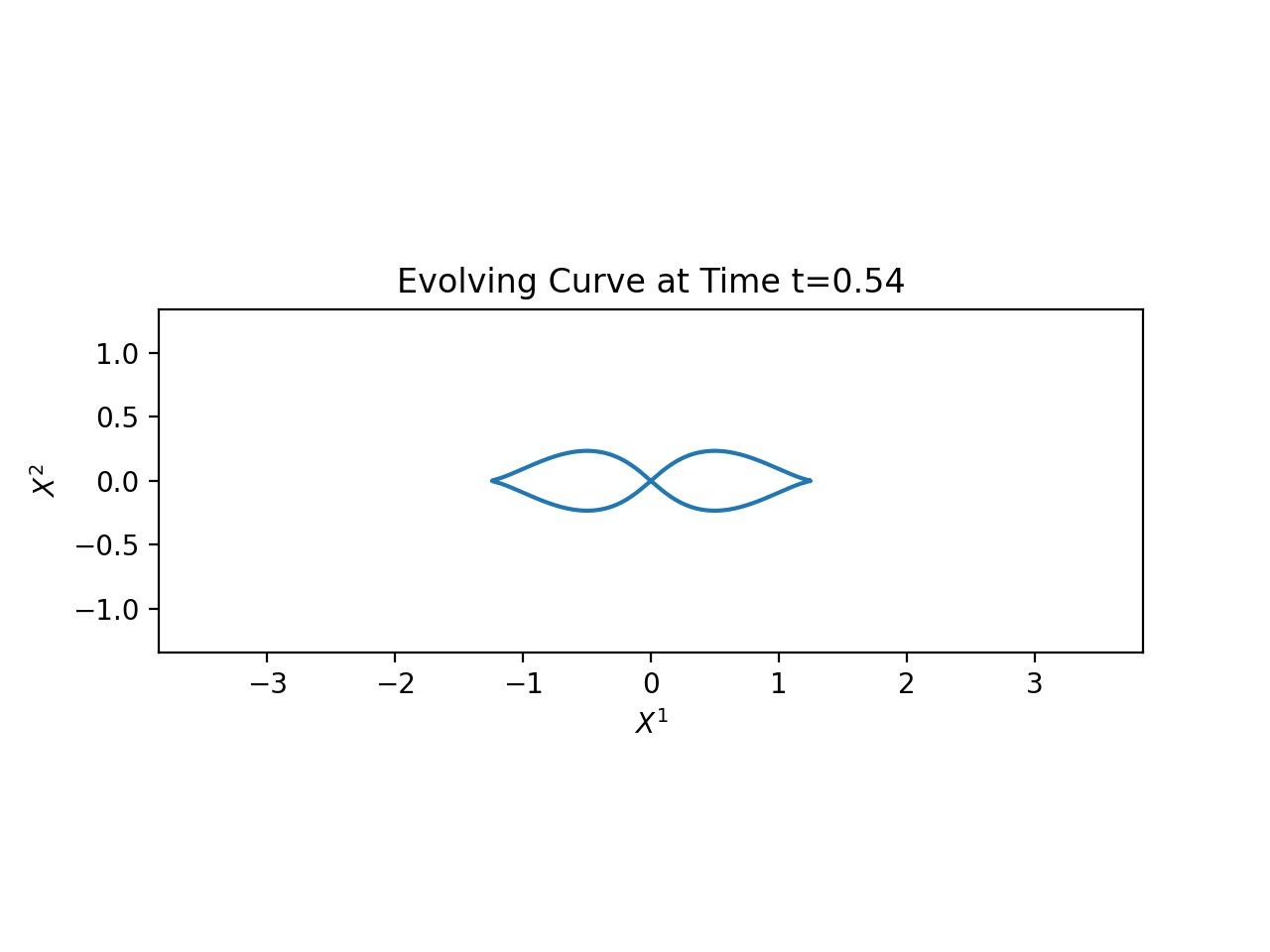}
  \includegraphics[scale=0.33]{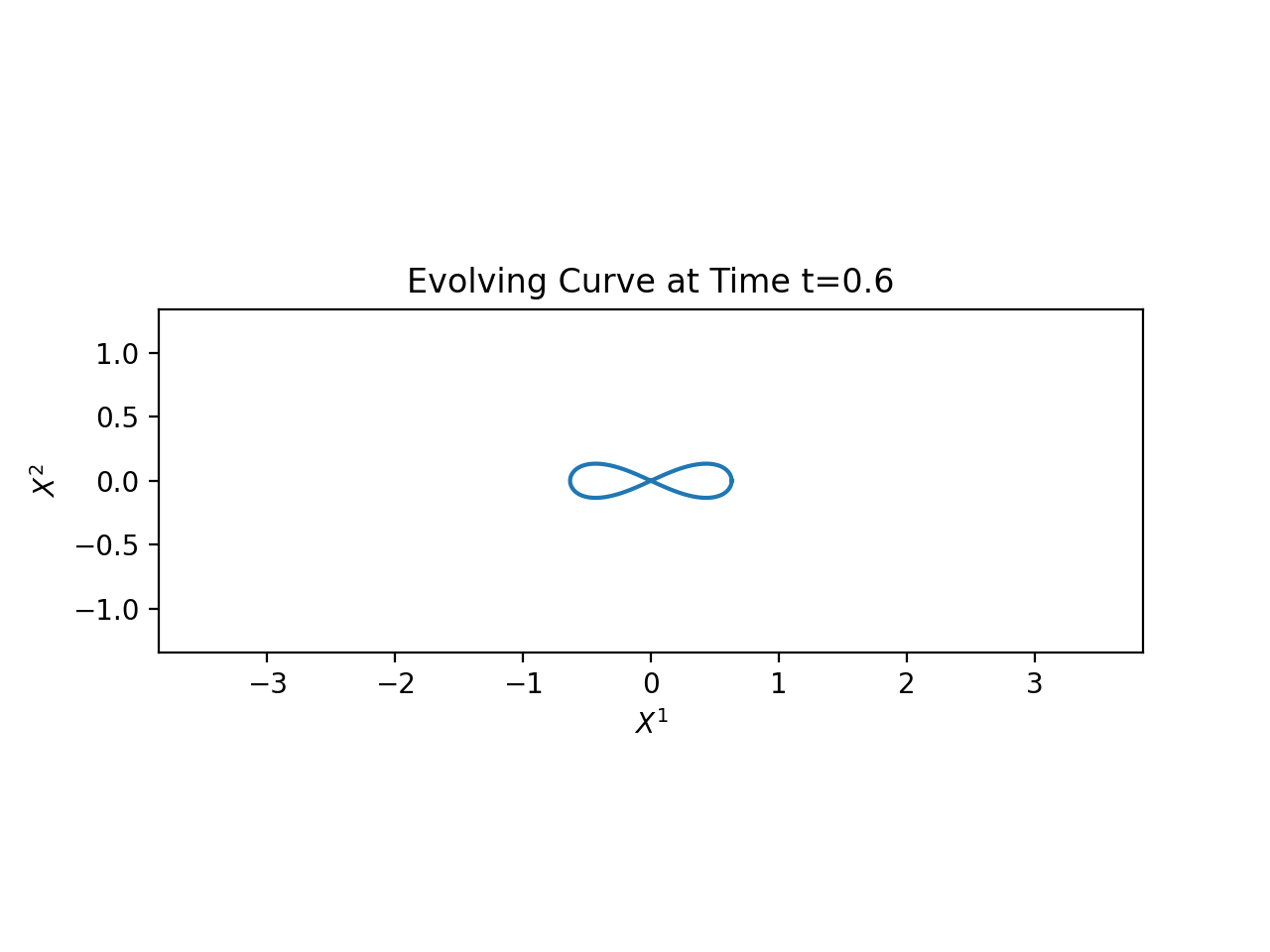}
  \includegraphics[scale=0.33]{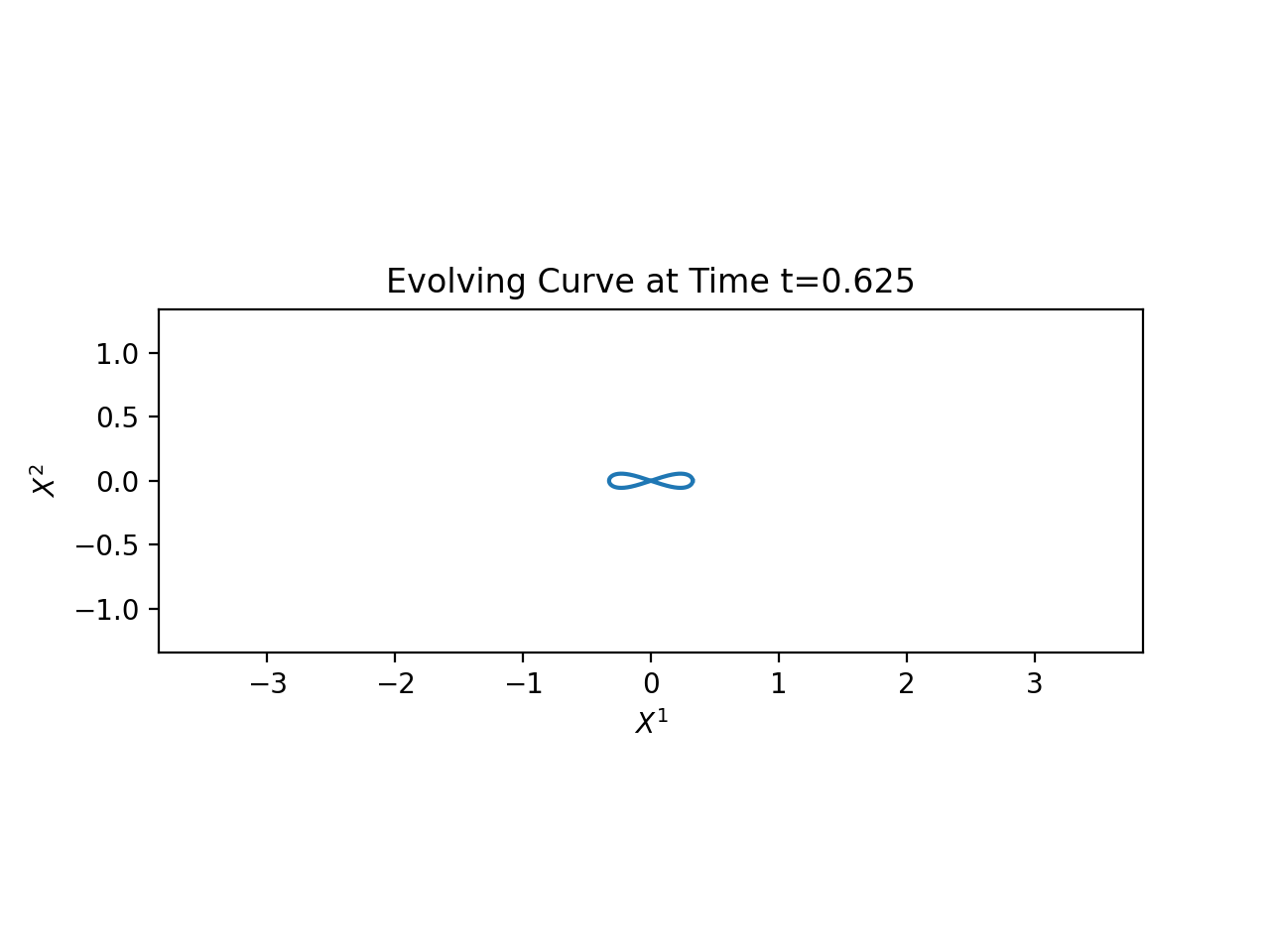}
\caption{Evolution of a convoluted curve with initial parametrization
  given by $X_0(r)=\bigl( 3\cos(6\pi r),\sin(16\pi r)\bigr)$,
  $r\in[0,1)$. There are two singularity times during this evolution,
  the first sheds four loops simultaneously ($t\simeq 0.282 $) and the
  second two ($t\simeq 0.54$).}
\label{fig:convolutedCurve}
\end{center}
\end{figure}

\begin{figure}
\begin{center}
  \includegraphics[scale=0.5]{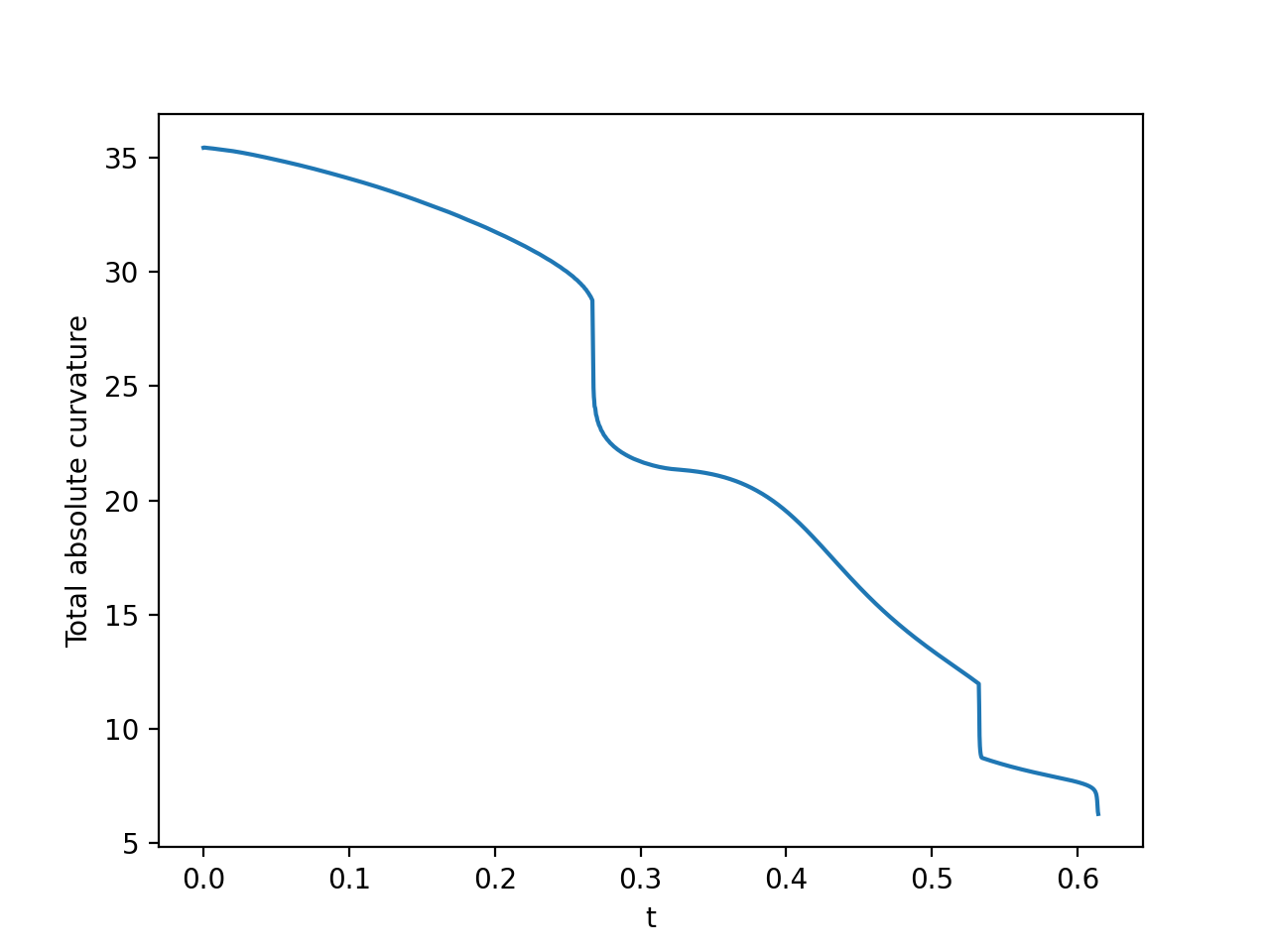}
  \includegraphics[scale=0.5]{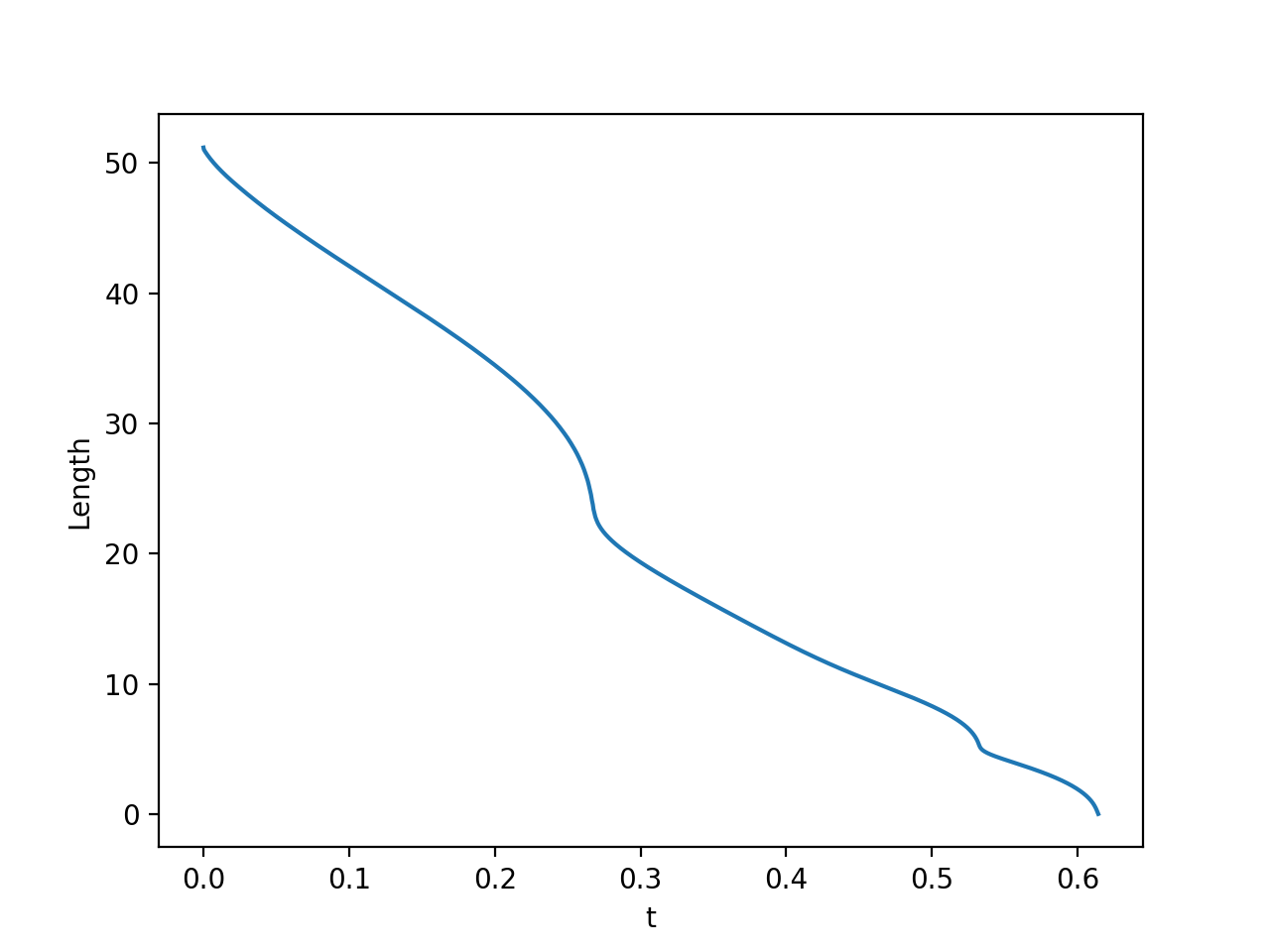}
\caption{Length and total absolute curvature along the orbit of the
  convoluted curve of Figure \ref{fig:convolutedCurve}}
\label{fig:l+tac}
\end{center}
\end{figure}

In Figure \ref{fig:convolutedCurve} we show the evolution of curve
with many self-intersections originating in the immersed curve
determined by the parametrization given by
$$
X_0(s)=\bigl( 3\cos(3s),\sin(8s)\bigr),\: s\in I. 
$$
This curve evolves by reducing the number of extremity points as
expected and, in doing so, it also reduces the number of
intersections. This does not always happen via the formation of
singularities. Singularities are only observed at time $t\simeq 0.282$, when
four loops are shedded simultaneously, and at time $t\simeq 0.54$, when two
additional loops disappear. At the extinction time, the curve is
asymptotic two a doubly covered segment after having simplified to an
infinity loop. In Figure \ref{fig:l+tac} we also plot the corresponding evolution of
length and total absolute curvature. The singular times are clearly
visibile in both graphs. The behavior of the numerical total absolute curvature
shortly before the singularity times is a rough approximation since
the loops disappear in the singularities and a fixed number of
discrete parameter points are used in the computation that do not
suffice to resolve the total curvature of a small loop.

\subsection{Linear Diffusion}
\begin{figure}
\begin{center}
  \includegraphics[scale=0.5]{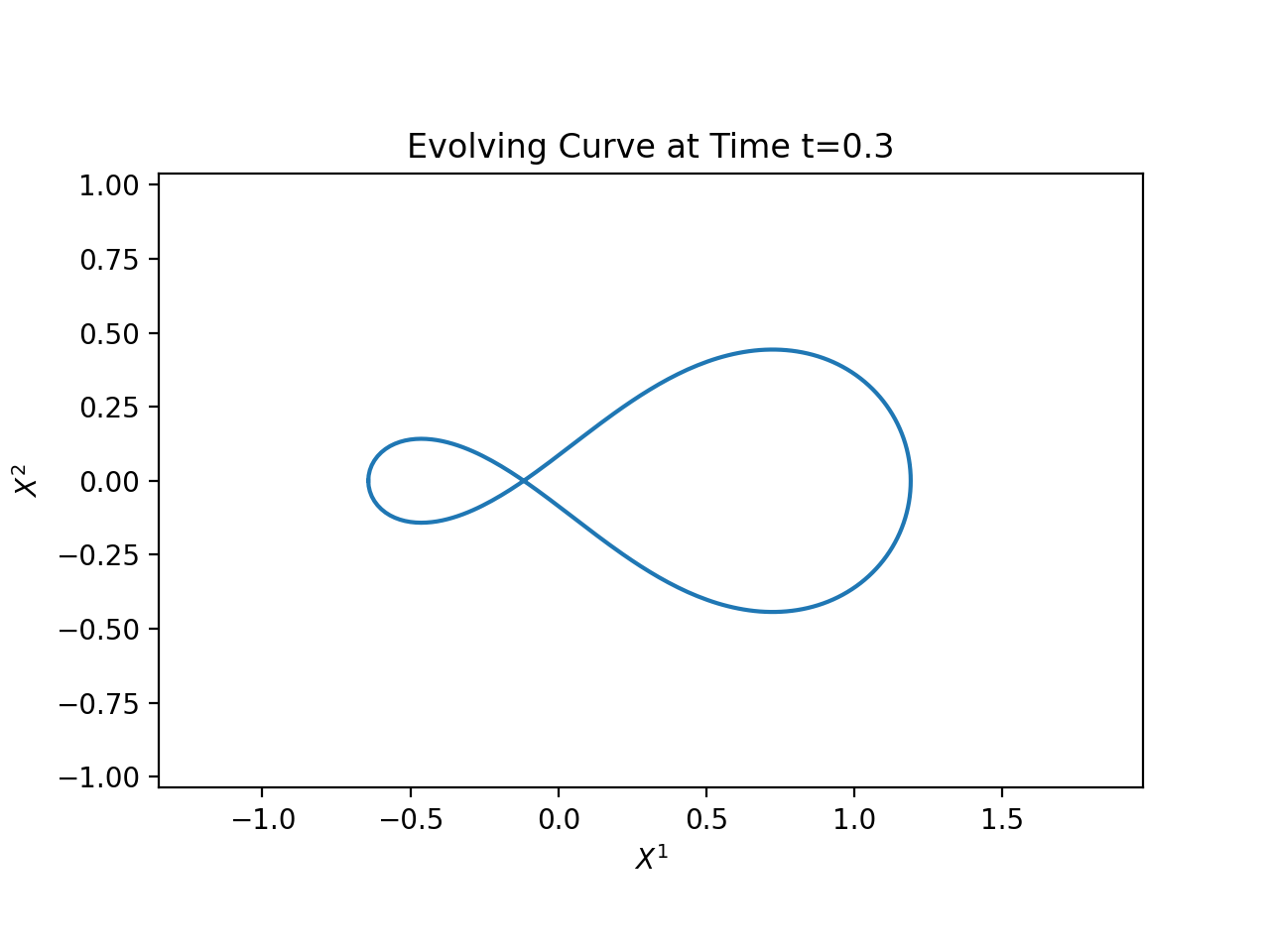}
  \includegraphics[scale=0.5]{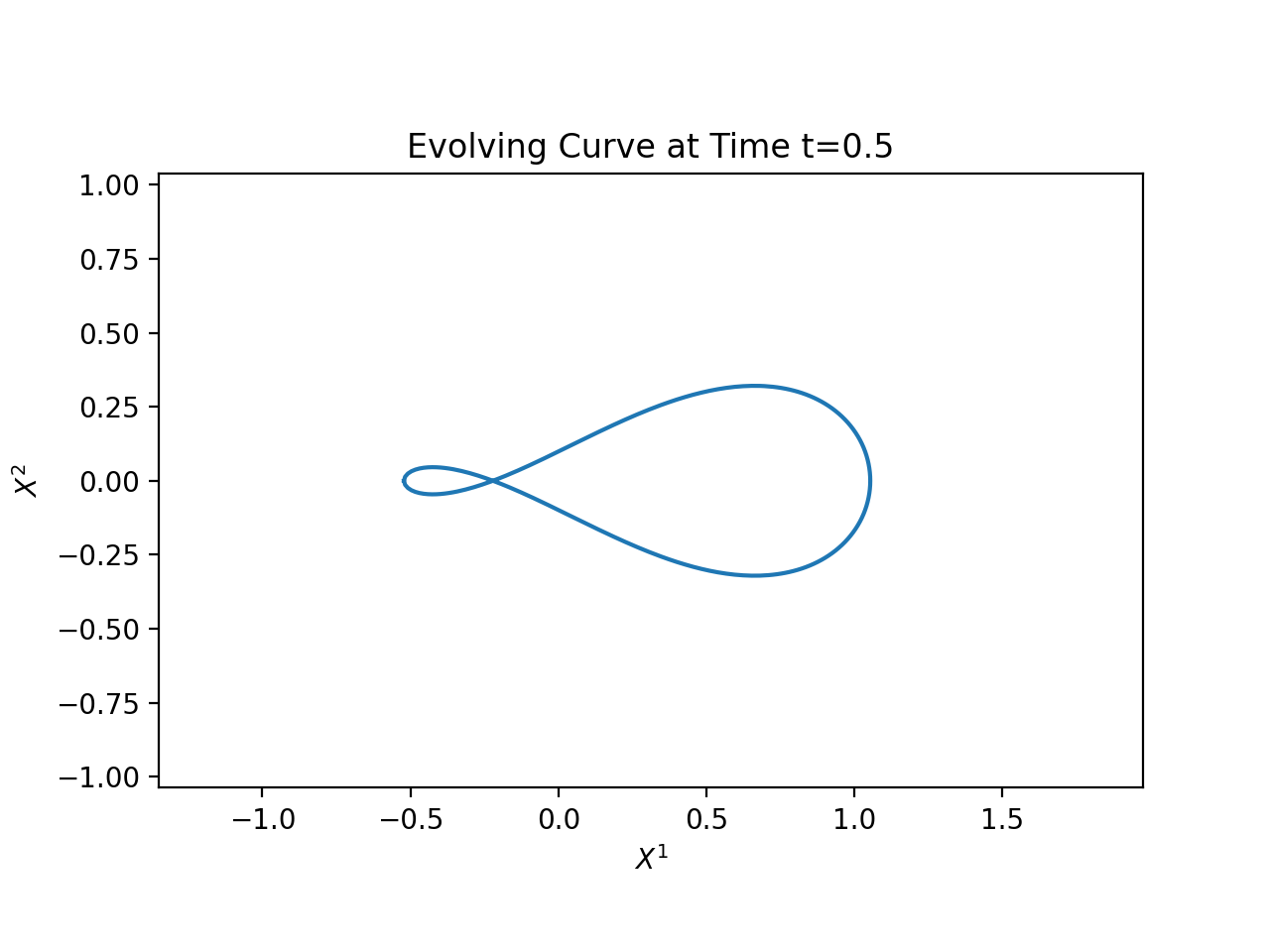}\\
  \includegraphics[scale=0.5]{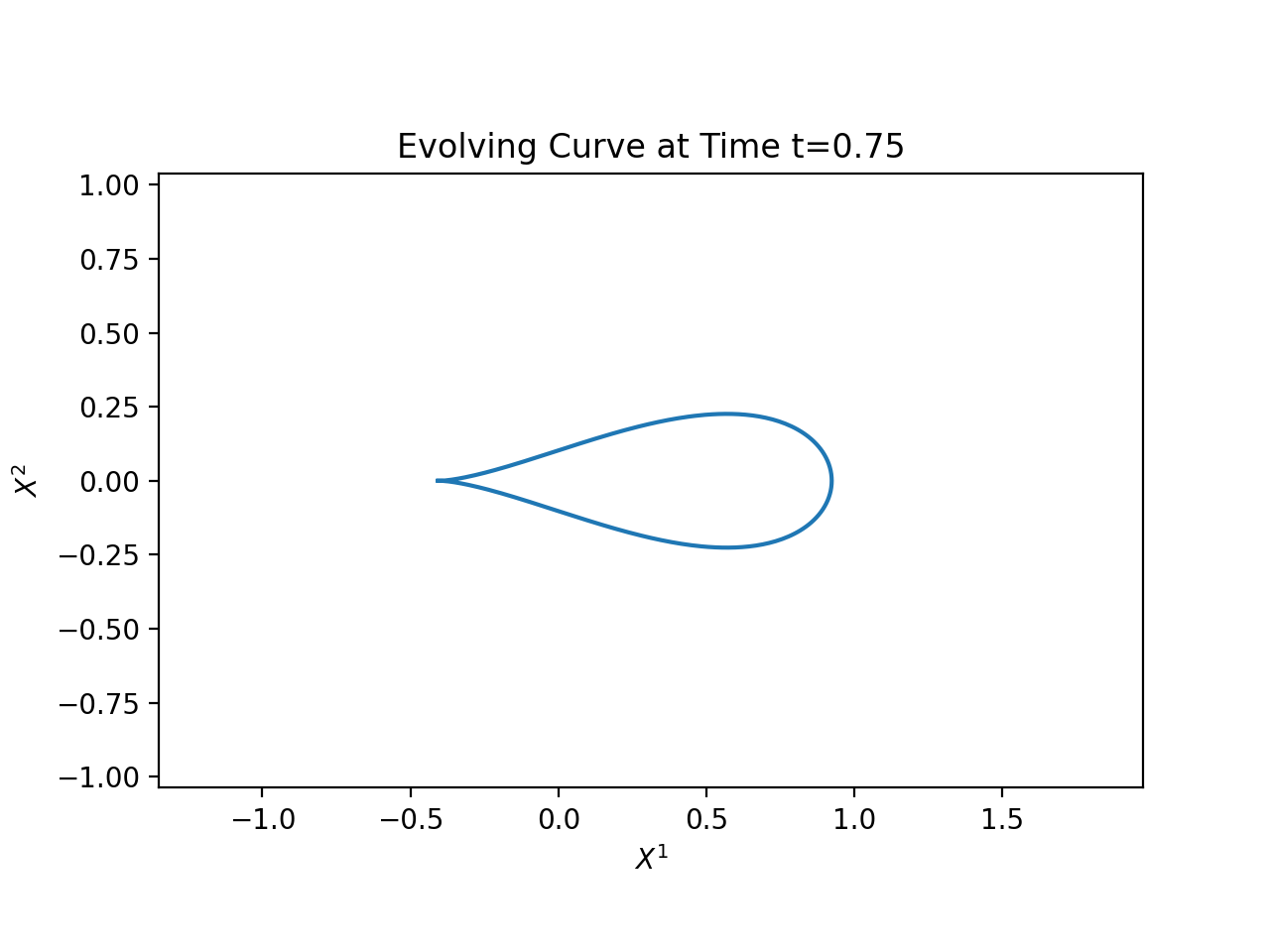}
  \includegraphics[scale=0.5]{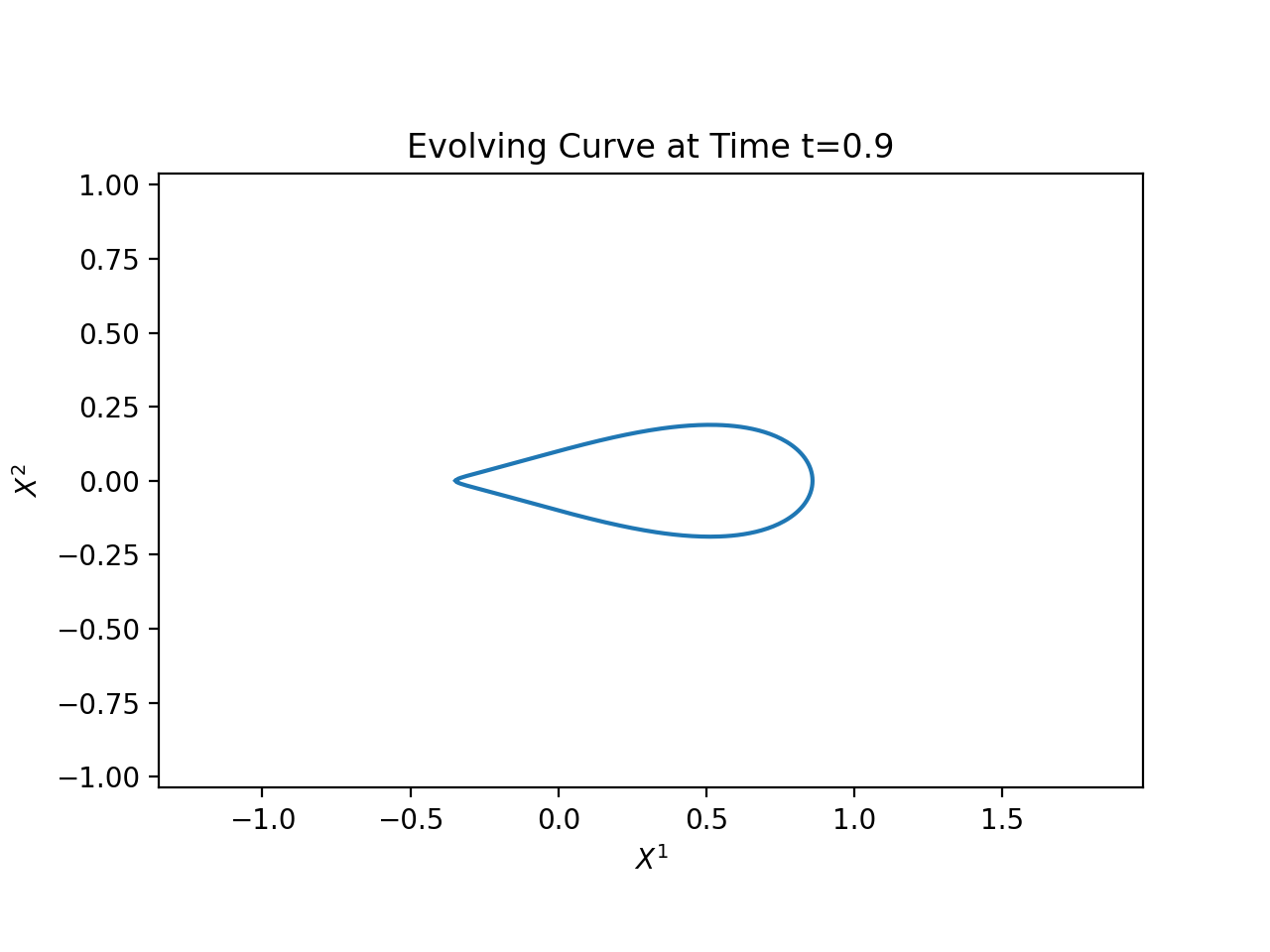}
\caption{Evolution of an infinity like shape consisting of two circles
  with different radii touching in the origin by the simple heat flow \eqref{linDiff}.}
\label{fig:linDiff}
\end{center}
\end{figure}

Here we highlight the fact that loop-shedding is not a consequence of
the nonlinear nature of the CSF as it is observed also when a curve is
evolved by pure diffusion. This is depicted in Figure
\ref{fig:linDiff}, where an initial curve consisting of two adjacent
circles touching in the origin and immersed like an infinity shape
with crossing is evolved by the heat equation with diffusivity
determined by its initial length
\begin{equation}\label{linDiff}
X_t=\frac{1}{L(X_0)^2}X_{ss},\: X(0)=X_0.
\end{equation}
The left circle is taken with radius $\frac{1}{4}$ and right circle
has radius $\frac{3}{4}$. This underscores the insight that
singularity formation at the level of the curve does not stem from a
singularity in any component of its parametrization but rather is due
to certain zero coalescence events as explained earlier in the paper.

\bibliography{csf}

\end{document}